\documentclass{article}[12pt]

\newtheorem{fed}{\textbf{Definition}}[section]
\newtheorem{thm}[fed]{\textbf{Theorem}}
\newtheorem{lemma}[fed]{\textbf{Lemma}}
\newtheorem{ex}[fed]{\textbf{Example}}
\newtheorem{rem}[fed]{\textbf{Remark}}
\newtheorem{prop}[fed]{\textbf{Proposition}}
\newtheorem{cor}[fed]{\textbf{Corollary}}
\usepackage{amssymb,graphicx,epsfig,psfrag,epic,eepic,latexsym}
\usepackage{amsmath}
\usepackage{mathrsfs}
\usepackage[dvips]{color}
\usepackage{xcolor}
\newcommand{\A}{\mathcal{A}}
\newcommand{\B}{\mathcal{B}}
\newcommand{\C}{\mathcal{C}}
\newcommand{\F}{\mathcal{F}}
\newcommand{\G}{\mathcal{G}}
\newcommand{\I}{\mathcal{I}}
\newcommand{\U}{\mathcal{U}}
\newcommand{\N}{\mathbb{N}}

\newcommand{\Z}{\mathbb{Z}}
\newcommand{\R}{\mathbb{R}}

\newcommand{\ddt}{\frac{d}{dt}}
\newcommand{\dds}{\frac{d}{ds}}

\usepackage{verbatim}

\newcommand{\BB}{\mathcal{B}}
\newcommand{\CC}{\mathcal{C}}
\newcommand{\HH}{\mathcal{H}}
\newcommand{\II}{\mathcal{I}}
\newcommand{\LL}{\mathcal{L}}
\newcommand{\NN}{\mathcal{N}}
\newcommand{\QQ}{\mathcal{Q}}
\newcommand{\UU}{\mathcal{U}}
\newcommand{\la}{\langle}
\newcommand{\ra}{\rangle}
\newcommand{\wt}{\widetilde}
\newcommand{\wh}{\widehat}
\newcommand{\eps}{\varepsilon}
\newcommand{\into}{\hookrightarrow}
\newcommand{\ol}{\overline}
\newcommand{\p}{\partial}

\begin{document}

\title{A variational approach to frozen planet orbits in helium}
\author{Kai Cieliebak, Urs Frauenfelder, Evgeny Volkov}
\maketitle
\begin{abstract}
We present variational characterizations of frozen planet orbits for the
helium atom in the Lagrangian and the Hamiltonian picture. They are
based on a nonlocal Levi-Civita regularization \cite{barutello-ortega-verzini} with different time
reparametrizations for the two electrons and lead to nonlocal functionals.
Within this variational setup, we deform the helium problem to one
where the two electrons interact only by their mean values and use
this to deduce the existence of frozen planet orbits. 
\end{abstract}
\tableofcontents

%%%%%%%%%%%%%%%%%%%%%%%%%%%%%%%%%%%%%%%%%%%%%%%%%%%%%%%%%%%%%%%%%%%%%%%%%%
\section{Introduction}
%%%%%%%%%%%%%%%%%%%%%%%%%%%%%%%%%%%%%%%%%%%%%%%%%%%%%%%%%%%%%%%%%%%%%%%%%%

{\em Frozen planet orbits} are periodic orbits in the helium atom
which play an important role in its semiclassical
treatment~\cite{tanner-richter-rost, wintgen-richter-tanner}. In such an
orbit both electrons lie on a line on the same side of the
nucleus. The inner electron undergoes consecutive collisions with the
nucleus while the outer electron (the actual ``frozen planet'') stays
almost stationary at some distance. See the following figure. 
%Despite their physical relevance,
%to our knowledge no rigorous proof of the existence of frozen planet
%orbits has appeared in the literature. One goal of this article is to
%provide such a proof. A second goal is to use this example as an
%illustration that variational methods involving nonlocal functionals
%can be successfully applied to physically relevant problems. 

\bigskip

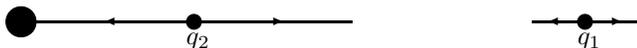
\setlength{\unitlength}{1mm}
\begin{picture}(1,1)
 \put(0,0){\circle*{4}}
 \put(2,0){\line(1,0){9}}
 \put(22,0){\vector(-1,0){11}}
 \put(23,0){\circle*{2}}
 \put(24,0){\vector(1,0){11}}
 \put(35,0){\line(1,0){9}}
 \put(68,0){\line(1,0){2}}
 \put(74,0){\vector(-1,0){4}}
 \put(75,0){\circle*{2}}
 \put(76,0){\vector(1,0){4}}
 \put(80,0){\line(1,0){2}}
 \put(22,-3){$q_2$}
 \put(74,-3){$q_1$}
\end{picture}
\bigskip
\bigskip

An interesting aspect of frozen planet orbits is that they cannot be
obtained using perturbative 
methods. Indeed, if the interaction between the two electrons is
switched off both electrons just fall into the nucleus. In order to
deal with this problem, the second author replaced in~\cite{frauenfelder} 
the instantaneous interaction between the two electrons by a mean
interaction and showed analytically that in this case there exists a
unique nondegenerate frozen planet orbit. 
\smallskip

{\bf Variational setup. }
In this paper we present a variational approach to frozen planet
orbits with instantaneous or mean interaction. One difficulty
lies in the collisions of the inner electron with the nucleus, which
need to be regularized in order to obtain a good functional analytic
setup. A traditional way to regularize two-body collisions is the
Levi-Civita regularization \cite{levi-civita}. In case of mean
interactions our problem has delay and the application of the
traditional Levi-Civita regularization becomes
problematic. Fortunately, in a recent interesting paper by Barutello,
Ortega and Verzini \cite{barutello-ortega-verzini} a new nonlocal
regularization was discovered. This regularization is motivated by the
traditional Levi-Civita regularization, but in sharp contrast to the
latter it is defined on the loop space and therefore fits well with
our problem. Even for loops without collisions this transformation is
quite intriguing. It is not smooth in the usual sense but scale smooth
in the sense of Hofer, Wysocki and Zehnder~\cite{hofer-wysocki-zehnder}.
\smallskip

We study two functionals $\BB_{av}$ and $\BB_{in}$
arising from regularizing frozen planet configurations for the mean
and instantaneous interaction, respectively, as well as their
linear interpolation $\BB_r=r\BB_{in}+(1-r)\BB_{av}$, $r\in[0,1]$.
In general, it is not
clear that critical points of a regularized action functional
correspond precisely to the regularized solutions. It might happen
that new exotic solutions appear as critical points,
see~\cite{barutello-ortega-verzini} for examples of this
phenomenon. Excluding such exotic critical points requires a careful
analysis and this analysis occupies the main part of this paper. In
particular, we prove (see Theorem~\ref{thm:interpol}) 
\smallskip

\textbf{Theorem\,A: }
For each $r\in[0,1]$, critical points of the regularized action
functional $\BB_r$ correspond precisely to frozen planet orbits for
the interpolated interaction.  
\smallskip

{\bf Symmetries. }
There is a special case of frozen planet orbits referred to as
\emph{symmetric frozen planet orbits}. These are frozen planet orbits
in which the outer electron has vanishing velocity whenever the inner
electron collides with the nucleus or is at maximal distance from the
nucleus, see~\cite{frauenfelder2}. We have (see
Theorem~\ref{thm:interpol-symm})
\smallskip

\textbf{Theorem\,B: }
The regularized action functional $\BB_r$ is invariant under an
involution such that the critical points fixed by the involution are
precisely the symmetric frozen planet orbits. 
\smallskip

In view of Theorem\,B one can think of symmetric frozen planet orbits as a
nonlocal generalization of brake orbits.  
\smallskip

{\bf Regularity. }
In order to study critical points of $\BB_r$, we introduce its
$L^2$-gradient as a map $\nabla\BB_r:X\to Y$ from a suitable Hilbert
manifold $X$ to a Hilbert space $Y$. It satisfies (see
Theorem~\ref{thm:Fredholm})
\smallskip

\textbf{Theorem\,C: }
For each $r\in[0,1]$, the $L^2$-gradient $\nabla\BB_r:X\to Y$ is a
$C^1$-Fredholm map of index $0$. 
\smallskip

This result is nontrivial because $\nabla\BB_r$ involves time
reparametrizations depending on points of $X$ as well as singular terms. 
Inspection of its proof shows that $\nabla\BB_r$ is {\em not} of class $C^2$.

Theorem\,C makes the functionals $\BB_r$ amenable to classical
variational methods such as index theory or Morse theory. In this
paper we content ourselves with computing the {\em mod $2$ Euler
number} $\chi(\nabla\BB_r)$, i.e., the count of zeroes modulo $2$
(after suitable perturbation, see Appendix~\ref{sec:Euler}). On
symmetric loops we find 
$$
   \chi(\nabla\BB_{in}) = \chi(\nabla\BB_{av})=1, 
$$
where the first equality follows from homotopy invariance of the mod
$2$ Euler number and the second one from a further deformation and
explicit computation (see Appendix~\ref{sec:Euler-mean-interaction}).
As a consequence, we obtain (see Theorem~\ref{thm:existence}):  

\textbf{Corollary\,C: }
For every $E<0$ there exists a symmetric frozen planet orbit of energy $E$.
\smallskip

%The proof relies on the {\em mod $2$ Euler number} 
%which assigns to each $C^1$-Fredholm map $f:X\to Y$ from a Hilbert
%manifold to a Hilbert space of index $0$ with compact zero set the
%number $\chi(f)$ of its zeroes modulo $2$ 
%Here the main analytical work consists in
%showing that the functional $\BB_{in}$ is of class $C^1$, which is
%nonobvious due to the loss of regularity resulting from the nonlocal
%terms (see Appendix~\ref{sec:diff}). 

{\bf Hamiltonian formulation. }
The regularized action functional $\BB_r$ has an intriguing
structure. It consists of two terms. The first term can be interpreted
as a kinetic energy, but for a nonlocal metric which depends on the
whole loop. The second term can be interpreted as the negative of a
nonlocal potential which is defined on the loop space of the
configuration space. We explain how in this situation a \emph{nonlocal
  Legendre transform} can be carried out which produces a nonlocal
Hamiltonian (see Section~\ref{sec:Ham}):
\smallskip

\textbf{Theorem\,D: }
After applying the nonlocal Legendre transform to the regularized
action functional $\BB_r$, the corresponding Hamiltonian delay
equation reproduces frozen planet orbits.  
\smallskip

Thus there are two nonlocal approaches to frozen planet orbits, one
Lagrangian and one Hamiltonian. This produces food for thought for
many interesting research projects. For instance, the Lagrangian
action functional has a Morse index at each frozen planet orbit. On
the other hand, the Hessian of the Hamiltonian action functional is a
Fredholm operator which gives rise to a nonlocal Conley-Zehnder index
at each frozen planet orbit. Since in the local case the Morse index
corresponds to the Conley-Zehnder index after Legendre transform we
may ask 
\smallskip

\textbf{Question\,1: }
How are the Morse index and the Conley-Zehnder index for frozen planet orbits related? 

We remark that for a simple but yet 
non-trivial delay equation, the regularized
free fall, recently there has been established equality in 
\cite{frauenfelder-weber}. 
\smallskip

The correspondence between these indices \cite{weber} is an important ingredient in
the celebrated adiabatic limit argument by Salamon and Weber
\cite{salamon-weber} relating Floer homology with the heat flow on
chain level. 
\smallskip

\textbf{Question\,2: }
Is there an analogon of the adiabatic limit argument of Salamon and
Weber for frozen planet orbits?

{\bf Acknowledgements. }
Supported by Deutsche Forschungsgemeinschaft (DFG, German Research Foundation) Projects CI 45/8-2
and FR 2637/2-2.

%%%%%%%%%%%%%%%%%%%%%%%%%%%%%%%%%%%%%%%%%%%%%%%%%%%%%%%%%%%%%%%%%%%%%%%%%%
\section{The Kepler problem}\label{sec:Kepler}
%%%%%%%%%%%%%%%%%%%%%%%%%%%%%%%%%%%%%%%%%%%%%%%%%%%%%%%%%%%%%%%%%%%%%%%%%%

When the interaction between the electrons is ignored the system
decouples into two one-electron systems, each of which is equivalent
to the Kepler problem in celestial mechanics. 
%with mass ratio $2$ (corresponding to charge $2$ of the nucleus).
In this section we recall the regularization
procedure of Barutello, Ortega and Verzini~\cite{barutello-ortega-verzini}
for the Kepler problem in the plane. It is based on the Levi-Civita
regularization map $\LL:\C\to\C$, $z\mapsto q=z^2$. Since we are interested
in the case that the position $q$ of the electron remains on the
positive real line, we view the Levi-Civita map as a
map to the non-negative reals
$\R\to\R_{\geq 0}$, $z\mapsto q=z^2$, 
see Equation~\eqref{eq:LC}.

%%%
\subsection{Levi-Civita transformation}\label{ss:LC}
%%%

In this subsection we describe the Levi-Civita transformation
%idea from~\cite{barutello-ortega-verzini}
without worrying about the regularity of the involved maps; precise
statements will be given in the following subsections.

We abbreviate by $S^1=\mathbb{R}/\mathbb{Z}$ the circle. We denote the
$L^2$-inner product of $z_1,z_2\in L^2(S^1,\mathbb{R})$ by
$$
  \langle z_1,z_2\rangle := \int_0^1 z_1(\tau) z_2(\tau) d\tau,
$$
and the $L^2$-norm of $z \in L^2(S^1, \mathbb{R})$ by
$$
  \|z\| := \sqrt{\langle z, z\rangle}.
$$
In the sequel we will work with Sobolev spaces $H^k=W^{k,2}$, but
the only relevant norms and inner products will be the ones from $L^2$.

Consider two maps
$$
   q:S^1\to\R_{\geq 0},\qquad z:S^1\to\R
$$
related by the {\em Levi-Civita transformation}
\begin{equation}\label{eq:LC}
   q(t) = z(\tau)^2
\end{equation}
for a time change $t\longleftrightarrow \tau$ satisfying
$0\longleftrightarrow 0$ and
\begin{equation}\label{eq:t-tau}
   \frac{dt}{q(t)} = \frac{d\tau}{\|z\|^2}. 
\end{equation}
This implies that the mean values of $q$ and $1/q$ are given by
\begin{equation}\label{bov1}
  \overline{q} :=
  \int_0^1 q(t)dt
  = \int_0^1 \frac{z(\tau)^4}{\|z\|^2} d\tau
  = \frac{\|z^2\|^2}{\|z\|^2}
\end{equation}
and
\begin{equation}\label{bov2}
   \int_0^1\frac{dt}{q(t)} = \frac{1}{\|z\|^2}.
\end{equation}
We will denote derivatives with respect to $t$ by a dot and
derivatives with respect to $\tau$ by a prime. Then the first and
second derivatives of $q$ and $z$ (where they are defined) are related by
\begin{equation}\label{bov3}
  \dot q(t)
  = 2z(\tau)z'(\tau)\frac{d\tau}{dt} = \frac{2\|z\|^2z'(\tau)}{z(\tau)}
\end{equation}
and
\begin{equation}\label{eq:ddot-q-z}
  \ddot q(t)
  = 2\|z\|^2\frac{z''(\tau)z(\tau)-z'(\tau)^2}{z(\tau)^2}\frac{d\tau}{dt}
  = \frac{2\|z\|^4}{z(\tau)^4}\bigl(z''(\tau)z(\tau)-z'(\tau)^2\bigr).
\end{equation}
Substituting $z^2$ and $z'^2$ by~\eqref{eq:LC} and~\eqref{bov3} this becomes
\begin{equation}\label{eq:ddot-q}
  \ddot q(t)
  = \frac{1}{q(t)}\Bigl(2\|z\|^4\frac{z''(\tau)}{z(\tau)} - \frac{\dot q(t)^2}{2}\Bigr).
\end{equation}
The $L^2$-norm of the derivative of $q$ is given by 
\begin{eqnarray}\label{bov4}
  \|\dot{q}\|^2
  = \int_0^1 \dot{q}(t)^2 dt %\\ \nonumber
  = \int_0^1
  \frac{4\|z\|^4z'(\tau)^2}{z(\tau)^2}\,\frac{z(\tau)^2}{\|z\|^2}d\tau
  %\\ \nonumber
  = 4\|z\|^2\|z'\|^2.
\end{eqnarray}

%%%
\subsection{Inverting the Levi-Civita transformation}\label{ss:LC-inv}
%%%

In this subsection we prove that, under suitable technical hypotheses,
the Levi-Civita transformation defines a 2-to-1 covering. 

We begin with a precise definition of the Levi-Civita transformation.
Let $z\in C^0(S^1,\R)$ be a continuous function with finite zero set
$$
  Z_z:=z^{-1}(0). 
$$
We associate to $z$ a $C^1$-map $t_z:S^1\to S^1$ by
\begin{equation}\label{eq:tz}
   t_z(\tau) := \frac{1}{\|z\|^2}\int_0^\tau z(\sigma)^2d\sigma.
\end{equation}
Note that $t_z(0)=0$ and
\begin{equation}\label{eq:dert}
   t_z'(\tau) = \frac{z(\tau)^2}{\|z\|^2}.
\end{equation}
Since $z$ has only finitely many zeroes, this shows that $t_z$ is
strictly increasing and we conclude

\begin{lemma}\label{lem:tz-homeo}
If $z\in C^0(S^1,\R)$ has only finitely many zeroes, then
the map $t_z:S^1\to S^1$ defined by~\eqref{eq:tz} is a homeomorphism. 
\hfill$\square$
\end{lemma}

It follows that $t_z \colon S^1 \to S^1$ has a continuous inverse
$$
  \tau_z:=t_z^{-1} \colon S^1 \to S^1.
$$
Since $t_z$ is of class $C^1$, the function $\tau_z$ is also of class $C^1$ on the
complement of the finite set $t_z(Z_z)$ with derivative
\begin{equation}\label{eq:taudot}
   \dot\tau_z(t) = \frac{\|z\|^2}{z(\tau_z(t))^2}.
\end{equation}
We define a continuous map $q:S^1\to \R_{\ge 0}$ by
\begin{equation}\label{eq:LC-def}
  q(t) := z(\tau_z(t))^2.
\end{equation}
Then the two maps $z,q$ are related by the Levi-Civita
transformation~\eqref{eq:LC} with $\tau=\tau_z$. Their zero sets 
$$
  Z_z=z^{-1}(0)\quad \text{and}\quad Z_q:=q^{-1}(0)=t_z(Z_z)
$$ 
are in bijective correspondence via $t_z$ (or equivalently $\tau_z$).
Moreover, by~\eqref{bov2} we have
$$
   \int_0^1\frac{ds}{q(s)} = \frac{1}{\|z\|^2} <\infty.
$$
Conversely, suppose we are given a map $q\in C^0(S^1,\R_{\geq 0})$
with finite zero set $Z_q$ satisfying $\int_0^1\frac{ds}{q(s)}<\infty$. 
We associate to $q$ the time reparametrization $\tau_q:S^1\to S^1$,
\begin{equation}\label{eq:tauq}
   \tau_q(t) :=
   \Bigl(\int_0^1\frac{ds}{q(s)}\Bigr)^{-1} \int_0^t\frac{1}{q(s)}ds. 
\end{equation}
Then $\tau_q(1)=1$, $\tau_q$ is of class $C^1$ outside the zero set $Z_q=q^{-1}(0)$
with derivative
\begin{equation}\label{eq:tauqdot}
   \tau_q'(t) =
   \Bigl(\int_0^1\frac{ds}{q(s)}\Bigr)^{-1}\frac{1}{q(t)}, \qquad t\in
   S^1\setminus Z_q.
\end{equation}
By~\cite[Lemma 2.1]{barutello-ortega-verzini} the map $\tau_q:S^1\to S^1$
is a homeomorphism whose inverse $t_q:=\tau_q^{-1}$ is of class $C^1$
and satisfies $t_q(1)=\tau_q^{-1}(1)=1$ and
\begin{equation}\label{eq:tqdot}
   t_q'(\tau) = \Bigl(\int_0^1\frac{ds}{q(s)}\Bigr)\, q(t_q(\tau)),\qquad
   \tau\in S^1.
\end{equation}
Suppose that $z:S^1\to\R$ is a continuous function satisfying
\begin{equation}\label{eq:LC-inv}
   z(\tau)^2 = q(t_q(\tau)).
\end{equation}
Then $z$ has finite zero set $Z_z=\tau_q(Z_q)$, so we can associate to
$z$ the homeomorphism $t_z:S^1\to S^1$ defined by~\eqref{eq:tz} and
its inverse $\tau_z$. We claim that 
\begin{equation}\label{eq:homeoqz}
   \tau_q=\tau_z\quad \text{and}\quad t_q=t_z
\end{equation}
It is enough to check the second equality. For this we compute 
$$
  \int_0^\tau z(\sigma)^2d\sigma
  =\int_0^\tau q(t_q(\sigma))d\sigma
  \stackrel{(*)}{=}\Bigl(\int_0^1\frac{ds}{q(s)}\Bigr)^{-1}
  \int_0^{t_q(\tau)}ds
  =\Bigl(\int_0^1\frac{ds}{q(s)}\Bigr)^{-1}t_q(\tau),
$$
where $(*)$ follows from the coordinate change $\sigma=\tau_q(s)$
and~\eqref{eq:tauqdot}. Evaluating at $\tau=1$ gives us
\begin{equation}\label{eq:meanqinv}
   \frac{1}{\|z\|^2}=\int_0^1\frac{ds}{q(s)}.
\end{equation}
Therefore,
$$
   t_z(\tau) = \frac{1}{||z||^2}\int_0^\tau z(\sigma)^2d\sigma = t_q(\tau)
$$
and~\eqref{eq:homeoqz} is established. Hence $q$ is the Levi-Civita
transform of $z$ defined by~\eqref{eq:LC-def}. 

Equation~\eqref{eq:LC-inv} does not uniquely determine $z$ for given
$q$ because the sign of $z$ can be arbitrarily chosen on each connected
component of $S^1\setminus Z_z$. If $Z_z$ consists of an {\em even}
number of points, then we can determine $z$ up to a global sign by the
requirement that $z$ changes its sign at each zero. Therefore, the
preceding discussion shows

\begin{lemma}\label{lem:LC}
The Levi-Civita transformation $z\mapsto q$ given by~\eqref{eq:LC-def}
defines for each even integer $m\in2\N$ a surjective 2-to-1 map
\begin{align*}
   \LL\colon &\{z\in C^0(S^1,\R)\mid \text{$z$ has precisely $m$ zeroes and
     switches sign at each zero}\} \\
   &\longrightarrow \{q\in C^0(S^1,\R_{\geq 0})\mid \text{$z$ has precisely $m$
     zeroes and }\int_0^1\frac{ds}{q(s)}<\infty\}.
\end{align*}
\hfill$\square$
\end{lemma}

The following lemma shows how additional regularity properties
translate between $z$ and $q$. Near each zero $t_*$ of $q$ we define
the local sign function
$$
  s_*(t) := \begin{cases}
  -1 & t<t_*, \\+1 & t>t_* \,.\end{cases}
$$
If $q$ is of class $C^1$ outside its zero set, we denote by 
\begin{equation}\label{eq:energy-q}
   E_q(t) := \frac{\dot q^2(t)}{2}-\frac{N}{q(t)},\qquad t\in S^1\setminus Z_q
\end{equation}
its {\em Kepler energy} at time $t$, for some fixed
$N>0$. By~\eqref{bov4} it corresponds under the Levi-Civita
transformation to
\begin{equation}\label{eq:energy-z}
   E_z(\tau) := \frac{2\|z\|^4z'(\tau)^2-N}{z(\tau)^2},\qquad \tau\in S^1\setminus Z_z.
\end{equation}

\begin{lemma}\label{lem:LC-reg}
Let $z,q$ be as in Lemma~\ref{lem:LC} related by the Levi-Civita
transform~\eqref{eq:LC-def}, and let $k$ be a nonnegative integer. Then the following hold. 

(a) $z\in H^1(S^1,\R)$ if and only if $q\in H^1(S^1,\R_{\geq 0})$;

(b) $z$ is of class $C^k$ outside its zeroes if and only if $q$ is of
class $C^k$ outside its zeroes;

(c) $z$ is of class $C^1$ on all of $S^1$ if and only if $q$ is of class
$C^1$ outside $Z_q$ and for each $t_*\in
Z_q$ the following limit exists:
$$
   \lim_{t_*\neq t\to t_*}s_*(t)\sqrt{q(t)}\dot q(t)\,;
$$
(d) $z$ is of class $C^1$ with transverse zeroes if and only if $q$ is
of class $C^1$ outside $Z_q$ and for each $t_*\in Z_q$ the limit in
(c) exists and is positive. 

(e) the energy $E_z:S^1\setminus Z_z\to\R$ is defined and extends to a
continuous function $S^1\to\R$ if and only if the energy $E_q:S^1\setminus 
Z_q\to\R$ is defined and extends to a continuous function $S^1\to\R$;

(f) the conditions in (e) imply those in (d). 
\end{lemma}

{\bf Proof: }
Part (a) follows immediately from formula~\eqref{bov4}. For (b)
just note that if $z$ is of class $C^k$ outside $Z_z$ then $t_z$ is of
class $C^{k+1}$ outside $Z_z$, so $\tau_z$ and therefore also $q$ is of class
$C^{k+1}$ outside $t_z(Z_z)=Z_q$, and the same in the reverse direction.

For (c) and (d) suppose that $z,q$ are of class $C^1$ outside their zero sets.
In the following we will always denote by $\tau,t$ times related by
the time transformation $t=t_z(\tau)=t_q(\tau)$. 
Consider a zero $t_*\in Z_q$ with corresponding $\tau_*\in Z_z$. 
Since $z$ switches sign at $\tau_*$ we can write
$$
   z(\tau) = \eps s_*(t)\sqrt{q(t)}
$$
for $\tau\neq\tau_*$ near $\tau_*$, with some sign $\eps\in\{-1,1\}$. 
Inserting this into~\eqref{bov3} and solving for $z'(\tau)$ yields
\begin{equation}
   z'(\tau) = \frac{\eps s_*(t)\sqrt{q(t)}\dot q(t)}{2\|z\|^2},
\end{equation}
from which (c) and (d) follow. 

Part (e) follows immediately from $E_z(\tau)=E_q(t_z(\tau))$. To see
that (e) implies (d), note first that the existence and continuity of
$E_q:S^1\setminus Z_q\to\R$ implies that $q$ is of class $C^1$ on
$S^1\setminus Z_q$. Suppose now that $E_q$ extends to a continuous
function $S^1\to\R$, so for each $t_*\in Z_q$ the limit
$$
   \lim_{t_8\neq t\to t_*}\Bigl(\frac{\dot q(t)^2}{2}-\frac{N}{q(t)}\Bigr)
$$
exists. This implies that $\dot q(t)^2\to\infty$ as $t\to t^*$, in
particular $\dot q(t)\neq 0$ for all $t\neq t_*$ close to $t_*$. 
Since $q(t)>0$ for $t\neq t_*$ and $q(t_*)=0$, this implies that
$$
  \dot q(t) = s_*(t)\sqrt{\dot q(t)^2} = s_*(t)\sqrt{2E_q(t)+\frac{2N}{q(t)}}
$$
for all $t\neq t_*$ close to $t_*$. It follows that
$$
  \lim_{t_*\neq t\to t_*}s_*(t)\sqrt{q(t)}\dot q(t)
  = \lim_{t_*\neq t\to t_*}\sqrt{2E(t)q(t)+2N} = \sqrt{2N} > 0,
$$
which is the condition in (d). This proves the lemma. 
\hfill$\square$

Note that $z\in H^1$ and the extension of $E_q$ to a continuous
function $S^1\to\R$ implies the existence of the integral
$$
  \int_0^1E_q(t)dt = \frac{\|\dot q\|^2}{2}-\int_0^1\frac{dt}{q(t)},
$$
and therefore $\int_0^1\frac{ds}{q(s)}<\infty$. Hence
Lemma~\ref{lem:LC-reg} implies

\begin{cor}\label{cor:LC}
For each even $m\in 2\N$ the Levi-Civita map $\LL$ of Lemma~\ref{lem:LC}
restricts to a surjective 2-to-1 map
$$
   \LL\colon \CC^1_{ce}(S^1,\R)\to \HH^1_{ce}(S^1,\R_{\geq 0}),\quad\text{where}
$$
\begin{itemize}
\item $\CC^1_{ce}(S^1,\R)$ denotes the set of $z\in C^1(S^1,\R)$ with
precisely $m$ zeroes such that all zeroes are transverse and the energy $E_z$ extends to a 
continuous function $S^1\to\R$, and
\item $\HH^1_{ce}(S^1,\R_{\geq 0})$ denotes the set of $q\in H^1(S^1,\R_{\geq 0})$ with
precisely $m$ zeroes such that $q$ is of class $C^1$ outside its
zeroes and the energy $E_q$ extends to a  continuous function $S^1\to\R$.
\end{itemize}
\hfill$\square$
\end{cor}

%%%
\subsection{Variational characterization of generalized solutions}
%%%

An electron moving in the electric field of a fixed nucleus with
charge $N>0$ is described by Newton's equation
\begin{equation}\label{eq:Kepler}
  \ddot q(t) = -\frac{N}{q(t)^2}.
\end{equation}
Alternatively, it describes the Kepler problem of a body of mass $1$
moving in the gravitational field of a body of mass $N$. Its periodic solutions
avoiding the origin are the critical points of the Lagrangian action functional
$$
  \mathcal{S}(q) := \frac12\int_0^1\dot q(t)^2dt + \int_0^1\frac{N}{q(t)}dt.
$$
Since all periodic solutions 
have collisions there are actually no critical points for this unregularized functional. 
Let $q$ and $z$ be related by the Levi-Civita transformation~\eqref{eq:LC}.
Using the relations~\eqref{eq:t-tau}, \eqref{bov3} and~\eqref{eq:meanqinv}, we rewrite the
Lagrangian action of $q$ in terms of $z$: 
\begin{align*}
  \mathcal{S}(q)
  &= \frac12\int_0^1\frac{4\|z\|^4z'(\tau)^2}{z(\tau)^2}\frac{z(\tau)^2}{\|z\|^2}d\tau
  + \frac{N}{\|z\|^2} \cr
  &= 2\|z\|^2\|z'\|^2 + \frac{N}{\|z\|^2}.
\end{align*}
We denote the resulting action functional of $z$ by
\begin{equation}\label{eq:Q}
  \QQ: H^1(S^1,\R)\setminus\{0\}\to\R,\qquad 
  \QQ(z) := 2\|z\|^2\|z'\|^2 + \frac{N}{\|z\|^2}. 
\end{equation}
Following~\cite{barutello-ortega-verzini} we call $q\in
H^1(S^1,\R_{\geq 0})$ a {\em generalized solution of~\eqref{eq:Kepler}} if 
\begin{enumerate}
\item the zero set $Z=q^{-1}(0)\subset S^1$ is finite and has an even
  number of elements;
\item on $S^1\setminus Z$ the map $q$ is smooth and satisfies~\eqref{eq:Kepler};
\item the energy
  $$
    E(t) := \frac{\dot q(t)^2}{2}-\frac{N}{q(t)},\qquad t\in
    S^1\setminus Z
  $$
extends to a continuous function $E:S^1\to\R$. 
\end{enumerate}
Note that the energy $E$ is then constant (by conservation of energy)
and negative (for $q$ to be bounded). 

\begin{thm}[Barutello, Ortega and Verzini~\cite{barutello-ortega-verzini}]\label{thm:BOV}
Under the Levi-Civita transformation~\eqref{eq:LC} with time
change~\eqref{eq:t-tau}, critical points $z:S^1\to\R$ of the action
functional $\QQ$ defined in~\eqref{eq:Q} are in 2-to-1 correspondence
with generalized solutions $q:S^1\to\R_{\geq 0}$ of~\eqref{eq:Kepler}.
\end{thm}

In the remainder of this section we will spell out the proof of this
theorem because it uses some ingredients that will also be needed in
later sections. 

%%%
\subsection{From critical points to generalized solutions}\label{ss:z-to-q-Kepler}
%%%

The differential of $\QQ$ at $z\in H^1(S^1,\R)\setminus\{0\}$ in
direction $v\in H^1(S^1,\R)$ is given by
\begin{equation}\label{eq:DQ}
  D\QQ(z)v
  = 4\|z\|^2 \langle z',v'\rangle
  + 4\|z'\|^2\langle z,v\rangle
  - \frac{2N}{\|z\|^4}\langle z,v\rangle
\end{equation}
This shows that a critical point $z$ has a weak second derivative and
satisfies the second order ODE with constant coefficient
\begin{equation}\label{eq:z-Kepler}
  z''(\tau) = a\,z(\tau),\qquad a = \frac{\|z'\|^2}{\|z\|^2} - \frac{N}{2\|z\|^6}.
\end{equation}
It follows that $z$ is smooth. Moreover, $z\in H^1(S^1,\R)$ 
implies boundedness of $z$ and thus forces $a<0$.
So $z$ is a shifted sine function. In particular, $z$ has {\em transverse
zeroes} in the sense that $z'(\tau)\neq 0$ whenever $z(\tau)=0$.
In particular, its zero set
$$
   Z:=\{\tau\in S^1\mid z(\tau)=0\}
$$
is finite. We associate to $z$ the smooth map $t_z:S^1\to S^1$ defined
by~\eqref{eq:tz}. By Lemma~\ref{lem:tz-homeo}, the map $t_z$ is a
homeomorphism with continuous inverse $\tau_z=t_z^{-1} \colon S^1 \to
S^1$. Since $t_z$ is smooth, the function $\tau_z$ is also smooth on the
complement of the finite set $t_z(Z)$ with derivative given by
equation~\eqref{eq:taudot}. 
We define a continuous map $q:S^1\to S^1$ by
$$
  q(t) := z(\tau_z(t))^2.
$$
Then the two maps $z,q:S^1\to\R$ are smooth except at finitely many
points and related by the Levi-Civita transformation~\eqref{eq:LC}
with $\tau=\tau_z$. 
Substituting $z''$ by~\eqref{eq:z-Kepler} in equation~\eqref{eq:ddot-q}
we get the following ODE for $q$ at points 
$t\in S^1\setminus t_z(Z)$:
\begin{equation}\label{eq:ddot-q-Kepler}
  \ddot q(t)
  = \frac{1}{q(t)}\Bigl(2\|z\|^4a-\frac{\dot q(t)^2}{2}\Bigr) 
  = \frac{1}{q(t)}\Bigl(c-\frac{\dot q(t)^2}{2}\Bigr)
\end{equation}
with the constant (using~\eqref{bov2} and~\eqref{bov4})
$$
  c := 2\|z\|^4a = 2\|z'\|^2\|z\|^2 - \frac{N}{\|z\|^2}
  = \frac{\|\dot q\|^2}{2} - \int_0^1\frac{N}{q(s)}ds.
$$
Since at a local maximum $t$ of $q$ we must have $\dot q(t)=0$ and
$\ddot q(t)<0$, it follows from~\eqref{eq:z-Kepler} that $c<0$, hence $\ddot q<0$ outside its zeroes.
Consider now two consecutive zeroes $t_-<t_+$ of $q$ and the smooth map
$$
   \beta := \frac{\ddot q}{q}: (t_-,t_+)\to\R_-\,.
$$
Then (omitting the $t$) we have $\beta q^2 = q\ddot q = c-\dot q^2/2$
and taking a time derivative yields
$$
   \dot\beta q^2+2\beta q\dot q = -\dot q\ddot q = -\beta q\dot q,
$$
hence
\begin{equation}\label{eq:betaqu}
  \dot\beta q = -3\beta\dot q.
\end{equation}

\begin{lemma}\label{lem:betamu}
Equation~\eqref{eq:betaqu} for functions $q>0$ and $\beta<0$ on
$(t_-,t_+)$ implies that
\begin{equation}\label{eq:betamu}
  \beta = -\frac{\mu}{q^3}
\end{equation}
on $(t_-,t_+)$ for some constant $\mu>0$.
\end{lemma}

\textbf{Proof: }
Dividing both sides of equation~\eqref{eq:betamu} by $q\beta$ yields
$$
  \ddt \log(-\beta)=-3\ddt \log(q),
$$
which by integration implies the lemma.
\hfill $\square$
\\

The lemma implies that
\begin{equation}\label{eq:Kepler-mu}
   \ddot q = -\mu/q^2
\end{equation}
on $(t_-,t_+)$. Combining this with~\eqref{eq:ddot-q-Kepler} yields 
\begin{equation}\label{eq:mu-q}
  -\frac{\mu}{q} = q\ddot q
  = c-\frac{\dot q^2}{2}
  = \frac{\|\dot q\|^2}{2} - \int_0^1\frac{N}{q(s)}ds -\frac{\dot q^2}{2}
\end{equation}
on $(t_-,t_+)$. Multiplying this equation by $-q$, smoothness of $q$ gives
$$
  \mu = \lim_{t\to t_\pm}-q(t)\Bigl(\frac{\|\dot q\|^2}{2} -
  \int_0^1\frac{N}{q(s)}ds -\frac{\dot q(t)^2}{2}\Bigr).
$$
This shows that the constant $\mu$ is the same for each interval
between consecutive zeroes of $q$, so equation~\eqref{eq:mu-q} holds
on all of $S^1\setminus t_z(Z)$. Now integrating~\eqref{eq:mu-q} over
$S^1$ yields
$$
   \mu = N,
$$
so~\eqref{eq:Kepler-mu} becomes Newton's equation~\eqref{eq:Kepler}.
Inserting $\mu=N$ in~\eqref{eq:mu-q} shows that the energy
$$
  E = \frac{\dot q^2}{2} - -\frac{N}{q}
  = \frac{\|\dot q\|^2}{2} - \int_0^1\frac{N}{q(s)}ds 
$$
is constant, so $q$ is a generalized solution of~\eqref{eq:Kepler}.

%%%
\subsection{From generalized solutions to critical points}\label{ss:q-to-z-Kepler}
%%%

Let now $q\in H^1(S^1,\R_{\geq 0})$ be a generalized solution of
equation~\eqref{eq:Kepler}. Integrating the constant energy yields
$$
    E = \int_0^1\frac{\dot q(t)^2}{2}dt - \int_0^1\frac{N}{q(t)}dt.
$$
Since $q\in H^1$, the first term on the right hand side is finite and
it follows that
$$
   \int_0^1\frac{1}{q(t)}dt < \infty. 
$$
As in Section~\ref{ss:LC-inv}, we associate to $q$ the time
reparametrization $\tau_q:S^1\to S^1$ defined by~\eqref{eq:tauq} and
its inverse $t_q=\tau_q^{-1}$. Recall that $\tau_q$
is smooth outside the zero set $Z_q=q^{-1}(0)$ and $t_q$ is of class $C^1$.
We define a continuous function $z:S^1\to\R$ by the condition
$$
   z(\tau)^2 = q(t_q(\tau))
$$
and the requirement that $z$ changes its sign at each zero. This is
possible because $q$ has an even number of zeroes, and it determines
$z$ up to a global sign. Using the change of variable $\tau=\tau_q(t)$
we find
$$
  \|z\|^2 = \int_0^1z(\tau)^2d\tau 
  = \Bigl(\int_0^1\frac{1}{q(s)}ds\Bigr)^{-1}\int_0^1z(\tau_q(t))^2 \frac{1}{q(t)}dt
  = \Bigl(\int_0^1\frac{1}{q(s)}ds\Bigr)^{-1}.
$$
It follows that $z$ and $q$ are related by the Levi-Civita
transformation~\eqref{eq:LC} with time change
$t\longrightarrow\tau=\tau_q$ satisfying~\eqref{eq:t-tau}. 

In the sequel we will drop the arguments $t$ and $\tau$. 
Combining equations~\eqref{eq:ddot-q-z} and~\eqref{eq:Kepler}
outside $Z_q$ we obtain
$$
  -\frac{N}{z^4} = -\frac{N}{q^2} = \ddot q = 2\frac{\|z\|^4}{z^4}(z''z-z'^2),
$$
hence
\begin{equation}\label{eq:Kepler-z}
  z''(\tau)z(\tau)-z'(\tau)^2 = -\frac{N}{2\|z\|^4},\qquad \tau\in S^1\setminus\tau_q(Z_q).
\end{equation}
Consider the function
$$
  \beta := \frac{z''}{z}: S^1\setminus\tau(Z_q)\to\R. 
$$
Inserting this into equation~\eqref{eq:Kepler-z} and using
equations~\eqref{eq:LC} and~\eqref{bov3} we find
$$
  -\frac{N}{2\|z\|^4} = \beta z^2 - z'^2
  = q\Bigl(\beta - \frac{\dot q^2}{4\|z\|^4}\Bigr),
$$
and solving for $\beta$ yields
$$
  \beta(\tau_q(t))
  = \frac{1}{2\|z\|^4}\Bigl(\frac{\dot q(t)^2}{2} - \frac{N}{q(t)}\Bigr)
  = \frac{1}{2\|z\|^4}E(t),\qquad t\in S^1\setminus Z_q.
$$
Since $q$ is a generalized solution, the energy $E$ is constant and
negative, thus $\beta(\tau) \equiv E/2\|z\|^4 < 0$ and the definition
of $\beta$ implies
$$
  z''(\tau) = \frac{E}{2\|z\|^4}z(\tau),\qquad \tau\in S^1\setminus\tau_q(Z_q).
$$
The solutions of this ODE are shifted sine functions. So the condition
that $z$ switches sign at each zero implies that it defines a smooth
function $z:S^1\to\R$ solving the ODE on all of $S^1$. 
Rewriting the energy via~\eqref{bov2} and~\eqref{bov4} as
$$
  E = \frac{\|\dot q\|^2}{2} - \int_0^1\frac{N}{q(t)}dt
  = 2\|z\|^2\|z'\|^2 - \frac{N}{\|z\|^2} 
$$
and inserting this into the ODE, we see that $z$ satisfies the
ODE~\eqref{eq:z-Kepler} and is therefore a critical point of $\QQ$.

Together with the previous subsection this concludes the proof of
Theorem~\ref{thm:BOV}.

%%%%%%%%%%%%%%%%%%%%%%%%%%%%%%%%%%%%%%%%%%%%%%%%%%%%%%%%%%%%%%%%%%%%%%%%%%
\section{Mean interaction}\label{sec:mean-interaction}
%%%%%%%%%%%%%%%%%%%%%%%%%%%%%%%%%%%%%%%%%%%%%%%%%%%%%%%%%%%%%%%%%%%%%%%%%%

In this section we consider a ``helium atom'' in which the two electrons interact by
the {\em mean values} $\overline{q}_i = \int_0^1q_i(t)dt$ according to
\begin{equation}\label{eq:mean-interaction}
\left\{\;
\begin{aligned}
  \ddot{q}_1(t) &= -\frac{2}{q_1(t)^2}+\frac{1}{(\overline{q}_1-\overline{q}_2)^2}, \cr
  \ddot{q}_2(t) &= -\frac{2}{q_2(t)^2}-\frac{1}{(\overline{q}_1-\overline{q}_2)^2}
\end{aligned}
\right.
\end{equation}
where we impose the condition
\begin{equation}\label{eq:mean-ineq}
  \ol q_1 > \ol q_2. 
\end{equation}

%%%
\subsection{Variational characterization of generalized solutions}\label{ss:var-mean}
%%%

Solutions of~\eqref{eq:mean-interaction} are the critical points of the action functional
$$
  \mathcal{S}_{av}(q_1,q_2)
  := \sum_{i=1}^2\Bigl(\frac12\int_0^1\dot q_i(t)^2dt + \int_0^1\frac{2}{q_i(t)}dt\Bigr)
  - \frac{1}{\ol q_1-\ol q_2}.
$$
For $i=1,2$ let $q_i$ and $z_i$ be related by Levi-Civita transformations
\begin{equation}\label{eq:LC-i}
   q_i(t) = z_i(\tau_i(t))^2
\end{equation}
for time changes $\tau_i(t)$ satisfying $\tau_i(0)=0$ and
\begin{equation}\label{eq:t-tau-i}
   \frac{dt}{q_i(t)} = \frac{d\tau_i(t)}{\|z_i\|^2}. 
\end{equation}
Note that we perform different time changes for the two electrons.
Then all the relations in Section~\eqref{ss:LC} hold with $(q,z,\tau)=(q_i,z_i,\tau_i)$.
In particular, we can use equation~\eqref{bov1} to rewrite the
interaction term in terms of the $z_i$: 
\begin{align*}
  - \frac{1}{\ol q_1-\ol q_2}
  = - \frac{1}{\frac{\|z_1^2\|^2}{\|z_1\|^2} - \frac{\|z_2^2\|^2}{\|z_2\|^2}}
  = - \frac{\|z_1\|^2\|z_2\|^2}{\|z_1^2\|^2\|z_2\|^2 - \|z_2^2\|^2\|z_1\|^2}.
\end{align*}
We denote the resulting {\em mean interaction} of $(z_1,z_2)$ by
\begin{equation}\label{eq:A}
  \mathcal{A}(z_1,z_2) := - \frac{\|z_1\|^2\|z_2\|^2}{\|z_1^2\|^2\|z_2\|^2 - \|z_2^2\|^2\|z_1\|^2}.
\end{equation}
This quantity is naturally defined on the space 
\begin{equation}\label{eq:Hav}
   \mathcal{H}_{av}^1 := \Bigg\{z=(z_1,z_2) \in
   H^1(S^1,\mathbb{R}^2)\;\Bigl|\; ||z_1||>0,\,\,||z_2||>0,\,\, 
   \frac{||z_1^2||^2}{||z_1||^2}>\frac{||z_2^2||^2}{||z_2||^2}\Bigg\}.
\end{equation}
Note that $\mathcal{H}_{av}^1$ is an open subset of the Hilbert space
$H^1(S^1,\mathbb{R}^2)$ and the last condition in its definition
corresponds to condition~\eqref{eq:mean-ineq}. On $\mathcal{H}_{av}^1$
we consider the functional
\begin{equation}\label{eq:Bav}
  \mathcal{B}_{av} \colon \mathcal{H}_{av}^1 \to \mathbb{R}, \qquad
  \mathcal{B}_{av}(z_1,z_2) := \QQ(z_1) + \QQ(z_2) + \mathcal{A}(z_1,z_2),
\end{equation}
with the functionals
\begin{equation*}%\label{eq:Q}
  \QQ(z_i) = 2\|z_i\|^2\|z_i'\|^2 + \frac{2}{\|z_i\|^2}
\end{equation*}
from equation~\eqref{eq:Q} with charge $N=2$. 

We call $(q_1,q_2)\in H^1(S^1,\R_{\geq 0}\times\R_{\geq 0})$ a {\em
  generalized solution of~\eqref{eq:mean-interaction}} if for $i=1,2$
the following holds:
\begin{enumerate}
\item $\ol q_1>\ol q_2$;
\item the zero sets $Z_i=q_i^{-1}(0)\subset S^1$ are finite and each have an even
  number of elements;
\item the restricted maps $q_i:S^1\setminus Z_i\to\R_{\geq 0}$ are
  smooth and satisfy~\eqref{eq:mean-interaction}; 
\item the energies
  $$
    E_i(t) := \frac{\dot q_i(t)^2}{2}-\frac{2}{q_i(t)},\qquad t\in
    S^1\setminus Z_i
  $$
extend to continuous functions $E_i:S^1\to\R$. 
\end{enumerate}
Note that the individual energies $E_i$ need not be constant, but
their sum is constant and negative. 

\begin{thm}[Generalized solutions with mean interaction]\label{thm:mean-interaction}
Under the Levi-Civita transformations~\eqref{eq:LC-i} with time
changes~\eqref{eq:t-tau-i}, critical points $(z_1,z_2)$ of the action
functional $\mathcal{B}_{av}$ defined in~\eqref{eq:Bav} are in 4-to-1 correspondence
with generalized solutions $(q_1,q_2)$ of~\eqref{eq:mean-interaction}.
\end{thm}

The proof of this theorem will take up the remainder of this section.

%%%
\subsection{The differential of $\mathcal{B}_{av}$}\label{ss:diffav}
%%%

The differential of the mean interaction $\mathcal{A}$ at
$(z_1,z_2)\in\mathcal{H}_{av}^1$ in direction $(v_1,v_2)\in
H^1(S^1,\R^2)$ is given by 
\begin{eqnarray}\label{eq:DA}
  D\mathcal{A}[z_1,z_2](v_1,v_2) \nonumber
  &=& -2\frac{||z_2||^2 \cdot \langle z_1,v_1\rangle+||z_1||^2 \cdot
    \langle z_2,v_2 \rangle}{||z_1^2||^2\cdot ||z_2||^2-||z_2^2||^2
    \cdot ||z_1||^2}\\ \nonumber
  & &+2\frac{||z_1||^2 \cdot ||z_2||^2\Big(2||z_2||^2 \cdot \langle
    z_1^3, v_1 \rangle+||z_1^2||^2 \cdot \langle z_2,v_2\rangle\Big) }
    {\big(||z_1^2||^2\cdot ||z_2||^2-||z_2^2||^2 \cdot
      ||z_1||^2\big)^2 }\\ \nonumber
  & &-2\frac{||z_1||^2 \cdot ||z_2||^2\Big(2||z_1||^2 \cdot \langle
      z_2^3,v_2 \rangle+||z_2^2||^2 \cdot \langle z_1,v_1\rangle
      \Big)}{\big(||z_1^2||^2\cdot ||z_2||^2-||z_2^2||^2 \cdot
      ||z_1||^2\big)^2 }\\ \nonumber
%  &=& +2\frac{||z_2||^2 \cdot \big(||z_2^2||^2\cdot
%      ||z_1||^2-||z_1^2||^2\cdot ||z_2||^2\big)-||z_1||^2 \cdot
%      ||z_2||^2 \cdot ||z_2^2||^2}{\big(||z_1^2||^2\cdot ||z_2||^2 
%    -||z_2^2||^2 \cdot ||z_1||^2\big)^2}\langle z_1,v_1 \rangle \\
%  & &+2\frac{||z_1||^2 \cdot \big(||z_2^2||^2\cdot
%      ||z_1||^2-||z_1^2||^2\cdot ||z_2||^2\big)+||z_1||^2 \cdot
%      ||z_2||^2 \cdot ||z_1^2||^2}{\big(||z_1^2||^2\cdot ||z_2||^2 
%    -||z_2^2||^2 \cdot ||z_1||^2\big)^2}\langle z_2,v_2 \rangle\\
%  & &+4\frac{||z_1||^2\cdot ||z_2||^4}{\big(||z_1^2||^2\cdot
%      ||z_2||^2-||z_2^2||^2 \cdot ||z_1||^2\big)^2} \langle z_1^3,v_1\rangle\\
%  & &-4\frac{||z_1||^4\cdot ||z_2||^2}{\big(||z_1^2||^2\cdot
%      ||z_2||^2-||z_2^2||^2 \cdot ||z_1||^2\big)^2} \langle z_2^3,v_2\rangle\\
  &=& -2\frac{||z_2||^4 \cdot ||z_1^2||^2}{\big(||z_1^2||^2\cdot ||z_2||^2
    -||z_2^2||^2 \cdot ||z_1||^2\big)^2}\langle z_1,v_1 \rangle \\ \nonumber
  & &+2\frac{||z_1||^4 \cdot ||z_2^2||^2}{\big(||z_1^2||^2\cdot ||z_2||^2
    -||z_2^2||^2 \cdot ||z_1||^2\big)^2}\langle z_2,v_2 \rangle\\ \nonumber
  & &+4\frac{||z_1||^2\cdot ||z_2||^4}{\big(||z_1^2||^2\cdot
      ||z_2||^2-||z_2^2||^2 \cdot ||z_1||^2\big)^2} \langle
    z_1^3,v_1\rangle\\ \nonumber
  & &-4\frac{||z_1||^4\cdot ||z_2||^2}{\big(||z_1^2||^2\cdot
      ||z_2||^2-||z_2^2||^2 \cdot ||z_1||^2\big)^2} \langle z_2^3,v_2\rangle
\end{eqnarray}
Combined with equation~\eqref{eq:DQ} with $z=z_i$, $v=v_i$ and $N=2$
for the differentials $D\QQ(z_i)v_i$, this yields the differential of $\B_{av}$:
\begin{eqnarray}\label{eq:DBav}
  D\mathcal{B}_{av}[z_1,z_2](v_1,v_2) \nonumber
  &=& 4 \sum_{i=1}^2 \bigg(||z_i||^2 \langle z_i',v_i'\rangle
    + ||z_i'||^2 \cdot \langle z_i,v_i\rangle
    - \frac{\langle z_i,v_i\rangle}{||z_i||^4}\bigg) \\ \nonumber
  & & -2\frac{||z_2||^4 \cdot ||z_1^2||^2}{\big(||z_1^2||^2\cdot ||z_2||^2
    -||z_2^2||^2 \cdot ||z_1||^2\big)^2}\langle z_1,v_1 \rangle \\
  & & +2\frac{||z_1||^4 \cdot ||z_2^2||^2}{\big(||z_1^2||^2\cdot ||z_2||^2
    -||z_2^2||^2 \cdot ||z_1||^2\big)^2}\langle z_2,v_2 \rangle\\ \nonumber
  & &+4\frac{||z_1||^2\cdot ||z_2||^4}{\big(||z_1^2||^2\cdot
      ||z_2||^2-||z_2^2||^2 \cdot ||z_1||^2\big)^2} \langle z_1^3,v_1\rangle\\ \nonumber
  & &-4\frac{||z_1||^4\cdot ||z_2||^2}{\big(||z_1^2||^2\cdot
      ||z_2||^2-||z_2^2||^2 \cdot ||z_1||^2\big)^2} \langle z_2^3,v_2\rangle
\end{eqnarray}

%%%
\subsection{Critical points of $\mathcal{B}_{av}$}
%%%

%The goal of this part is to describe critical points of $\mathcal{B}_{av}$ in terms of a 
%differential equation they must satisfy. Then we use the differential equation to derive certain 
%properties of critical points.
%We say that a function $f\in C^1(S^1,\R)$ {\it has transverse zeros} if for all $t\in S^1$
%with $f(t)=0$ we have $f'(t)\ne 0$.
Equation~\eqref{eq:DBav} leads to the characterization of critical
points of $\mathcal{B}_{av}$:

\begin{prop}\label{prop:critav}
A point $(z_1,z_2)\in \mathcal{H}_{av}^1$ is a critical point of
$\mathcal{B}_{av}$ if and only if $(z_1,z_2)$ is smooth and solves
the system of (uncoupled!) ODEs
\begin{equation}\label{eq:z-av}
\left\{
\begin{aligned}
  z_1'' &= a_1 z_1 + b_1z_1^3 \\
  z_2'' &= a_2 z_2 + b_2z_2^3
\end{aligned}
\right.
\end{equation}
with the constants
\begin{eqnarray*}
a_1&=&\frac{||z_1'||^2}{||z_1||^2}-\frac{1}{||z_1||^6}-\frac{||z_2||^4 \cdot ||z_1^2||^2}{2||z_1||^2 \cdot\big(||z_1^2||^2\cdot ||z_2||^2
-||z_2^2||^2 \cdot ||z_1||^2\big)^2}\\
b_1&=&+\frac{||z_2||^4}{\big(||z_1^2||^2\cdot ||z_2||^2-||z_2^2||^2 \cdot ||z_1||^2\big)^2}\\
a_2&=&\frac{||z_2'||^2}{||z_2||^2}-\frac{1}{||z_2||^6}+\frac{||z_1||^4 \cdot ||z_2^2||^2}{2||z_2||^2 \cdot\big(||z_1^2||^2\cdot ||z_2||^2
-||z_2^2||^2 \cdot ||z_1||^2\big)^2}\\
b_2&=&-\frac{||z_1||^4}{\big(||z_1^2||^2\cdot ||z_2||^2-||z_2^2||^2 \cdot ||z_1||^2\big)^2}.
\end{eqnarray*}
\end{prop}

\textbf{Proof: } 
From equation~\eqref{eq:DBav} we see that $(z_1,z_2)$ is a critical
point of $\mathcal{B}_{av}$ if and only if $z_1$ and $z_2$ have
weak second derivatives and satisfy the system of
ODEs~\eqref{eq:z-av}. Bootstrapping these equations
%using the Sobolev embedding $W^{k,2}\rightarrow C^k(R)$ for $k=1,2\dots$ 
we conclude that $z_1$ and $z_2$ are smooth and the proposition follows.
\hfill $\square$

\begin{cor}\label{cor:mean-transverse-zeroes}
Suppose that $(z_1,z_2)$ is a critical point of
$\mathcal{B}_{av}$. Then $z_1$ and $z_2$ have transverse zeros. In
particular, their zero sets
$$
  Z_i := \big\{\tau \in S^1 \mid z_i(\tau)=0\big\}, \qquad i=1,2
$$
are finite.
\end{cor}

\textbf{Proof: } 
Arguing by contradiction, suppose that there exists a point
$\tau_0\in S^1$ such that $z_1(\tau_0)=z_1'(\tau_0)=0$. 
Then the function $z_1$ and the zero function both solve the first
equation in~\eqref{eq:z-av} with the same initial conditions at
$\tau_0$. By uniqueness of solutions of ODEs we conclude $z_1\equiv
0$, contradicting the condition $\|z_1\|>0$ in the definition of
$\mathcal{H}_{av}$. An analogous argument applies to $z_2$.
\hfill $\square$

%%%
\subsection{From critical points to generalized solutions}\label{ss:z-to-q-mean}
%%%

Let $(z_1,z_2)\in \mathcal{H}_{av}^1$ be a critical point of $\mathcal{B}_{av}$,
so by Proposition~\ref{prop:critav} the maps $z_1,z_2:S^1\to\R$ are
smooth and satisfy~\eqref{eq:z-av}. 
For $i=1,2$ we define the smooth maps $t_{z_i}:S^1\to S^1$ by
\begin{equation}\label{eq:ti}
   t_{z_i}(\tau) := \frac{1}{\|z_i\|^2}\int_0^\tau z_i(\sigma)^2d\sigma.
\end{equation}
Since by Corollary~\ref{cor:mean-transverse-zeroes} the map $z_i$ has
only finitely many zeroes, it follows from Lemma~\ref{lem:tz-homeo}
that $t_{z_i}:S^1\to S^1$ is a homeomorphism with continuous inverse 
$\tau_{z_i} \colon S^1 \to S^1$. We define continuous maps $q_i:S^1\to S^1$ by
$$
  q_i(t) := z_i(\tau_{z_i}(t))^2.
$$
Then for $i=1,2$ the maps $z_i,q_i:S^1\to\R$ are smooth except at finitely many
points and related by the Levi-Civita transformation~\eqref{eq:LC}
with $\tau=\tau_{z_i}$. As explained in Section~\ref{ss:var-mean}, the
last condition in the definition of $\mathcal{H}_{av}^1$ implies
$$
  \ol q_1>\ol q_2. 
$$
Let us now focus on $i=1$. Substituting $z_1''$ by~\eqref{eq:z-av} in
equation~\eqref{eq:ddot-q} with $q=q_1$ and $z=z_1$ 
we compute at points $t\in S^1\setminus t_{z_1}(Z_1)$:
\begin{eqnarray}
  \ddot{q}_1 \nonumber
  &=&\bigg(2||z_1'||^2 \cdot
  ||z_1||^2-\frac{2}{||z_1||^2}-\frac{||z_1||^2 \cdot ||z_2||^4 \cdot
  ||z_1^2||^2}{\big(||z_1^2||^2\cdot ||z_2||^2 
  -||z_2^2||^2 \cdot ||z_1||^2\big)^2}\bigg)\frac{1}{q_1}\\ \nonumber
& &+\frac{2||z_1||^4 \cdot ||z_2||^4}{\big(||z_1^2||^2\cdot ||z_2||^2-||z_2^2||^2 \cdot ||z_1||^2\big)^2}
-\frac{\dot{q}_1^2}{2q_1}\\ \nonumber
&=&\bigg(2||z_1'||^2 \cdot ||z_1||^2-\frac{2}{||z_1||^2}-\frac{||z_1||^4 \cdot ||z_2||^4 \cdot \frac{||z_1^2||^2}{
||z_1||^2}}{||z_1||^4 \cdot ||z_2||^4 \cdot\Big(\frac{||z_1^2||^2}{||z_1||^2}
-\frac{||z_2^2||^2}{||z_2||^2}\Big)^2}-\frac{\dot{q}_1^2}{2}\bigg)\frac{1}{q_1}\\ \nonumber
& &+\frac{2||z_1||^4 \cdot ||z_2||^4}{||z_1||^4 \cdot ||z_2||^4 \cdot\Big(\frac{||z_1^2||^2}{||z_1||^2}
-\frac{||z_2^2||^2}{||z_2||^2}\Big)^2}\\ \nonumber
&=&\bigg(2||z_1'||^2 \cdot ||z_1||^2-\frac{2}{||z_1||^2}-\frac{\frac{||z_1^2||^2}{
||z_1||^2}}{\Big(\frac{||z_1^2||^2}{||z_1||^2}
-\frac{||z_2^2||^2}{||z_2||^2}\Big)^2}-\frac{\dot{q}_1^2}{2}\bigg)\frac{1}{q_1}\\ \nonumber
& &+\frac{2}{\Big(\frac{||z_1^2||^2}{||z_1||^2}
-\frac{||z_2^2||^2}{||z_2||^2}\Big)^2}\\ \nonumber
&=&\bigg(\frac{||\dot{q}_1||^2}{2}-\int_0^1 \frac{2}{q_1(s)}ds-\frac{\overline{q}_1}{(\overline{q}_1-\overline{q}_2)^2}
-\frac{\dot{q}_1^2}{2}\bigg)\frac{1}{q_1}+\frac{2}{(\overline{q}_1-\overline{q}_2)^2}.
\end{eqnarray}
Thus $q_1$ satisfies the ODE
\begin{equation}\label{bov6}
  \ddot{q}_1 = \bigg(c_1
  -\frac{\dot{q}_1^2}{2}\bigg)\frac{1}{q_1}+\frac{2}{(\overline{q}_1-\overline{q}_2)^2}
\end{equation}
with the constant
\begin{equation}\label{eq:c1-mean}
  c_1 = \frac{||\dot{q}_1||^2}{2}-\int_0^1
  \frac{2}{q_1(s)}ds-\frac{\overline{q}_1}{(\overline{q}_1-\overline{q}_2)^2}. 
\end{equation}
At the global maximum $t_{\rm max}$ of $q_1$ equation~\eqref{bov6} becomes 
$$
  \frac{c_1}{q_1(t_{\rm max})} + \frac{2}{(\overline{q}_1-\overline{q}_2)^2} =
  \ddot q_1(t_{\rm max})\leq 0,
$$
hence
\begin{equation}\label{eq:c1-ineq}
   c_1 \leq -\frac{2q_1(t_{\rm max})}{(\overline{q}_1-\overline{q}_2)^2}.
\end{equation}
Let now $t_-<t_+$ be adjacent zeroes of $q_1$ and consider the smooth map
$$
  \beta_1:=\frac{\ddot{q}_1-\frac{1}{(\overline{q}_1-\overline{q}_2)^2}}{q_1}
  \colon (t_-,t_+) \to \mathbb{R}.
$$
From (\ref{bov6}) we obtain
$$
  \beta_1 q_1^2
  = c_1 -\frac{\dot{q}_1^2}{2} + \frac{q_1}{(\overline{q}_1-\overline{q}_2)^2}.
$$
With $q_1\leq q_1(t_{\rm max})$ and inequality~\eqref{eq:c1-ineq} this implies
$$
  \beta_1 q_1^2
  \leq -\frac{\dot{q}_1^2}{2} - \frac{q_1(t_{\rm
      max})}{(\overline{q}_1-\overline{q}_2)^2} < 0,
$$
hence $\beta_1<0$ on $(t_-,t_+)$. Differentiating both sides of the
equation for $\beta q_1^2$ we get
$$\dot{\beta}_1 q_1^2+2 \beta_1 q_1 \dot{q}_1=-\ddot{q}_1 \dot{q}_1+\frac{\dot{q}_1}{(\overline{q}_1-\overline{q}_2)^2}=-\beta_1 q_1 \dot{q}_1$$
and therefore
\begin{equation*}%\label{eq:betaqu}
  \dot{\beta}_1 q_1 = -3\beta_1 \dot{q}_1.
\end{equation*}
By Lemma~\ref{lem:betamu} this implies that 
\begin{equation*}%\label{eq:betamu}
  \beta_1 = -\frac{\mu}{q_1^3}
\end{equation*}
on $(t_-,t_+)$ for some constant $\mu>0$. By definition of $\beta_1$
this yields
\begin{equation}\label{bov7}
  \ddot{q}_1(t) =
  -\frac{\mu}{q_1(t)^2}+\frac{1}{(\overline{q}_1-\overline{q}_2)^2}
\end{equation}
for $t\in (t_-,t_+)$. Plugging this into~\eqref{bov6} we infer
\begin{equation}\label{bov7a}
  \mu = -\bigg(\frac{\|\dot{q}_1\|^2}{2}-\int_0^1 \frac{2}{q_1(s)}ds
  -\frac{\overline{q}_1}{(\overline{q}_1-\overline{q}_2)^2}
  -\frac{\dot{q}_1(t)^2}{2}\bigg)q_1(t)-\frac{q_1(t)^2}{(\overline{q}_1-\overline{q}_2)^2}
\end{equation}
for $t\in(t_-,t_+)$. In particular, use~\eqref{bov3} in the second identity to obtain 
\begin{equation}\label{bov8}
  \mu = \lim_{t \to t_\pm} \frac{\dot{q}_1(t)^2 q_1(t)}{2}
  = 2||z_1||^4 z_1'\big(\tau_{z_1}(t_\pm)\big)^2.
\end{equation}
We deduce from this that equation~\eqref{bov7} holds on $S^1 \setminus
t_{z_1}(Z_1)$ with a fixed $\mu$ independent of 
the connected component in $S^1 \setminus t_{z_1}(Z_1)$. We divide (\ref{bov7a}) by $q_1$ to get
\begin{equation}\label{bov7b}
  \frac{\mu}{q_1(t)}
  = -\frac{\|\dot{q}_1\|^2}{2}+\int_0^1
  \frac{2}{q_1(s)}ds+\frac{\overline{q}_1}{(\overline{q}_1-\overline{q}_2)^2} 
  +\frac{\dot{q}_1(t)^2}{2}-\frac{q_1(t)}{(\overline{q}_1-\overline{q}_2)^2}.
\end{equation}
Integrating this equation yields
$$
  \mu\int \frac{1}{q_1(t)}dt=2\int_0^1 \frac{1}{q_1(s)}ds,
$$
and therefore
$$\mu=2.$$
Thus equation~\eqref{bov7} becomes the first equation in~\eqref{eq:mean-interaction}.
Similarly, one obtains for $q_2$ the equation
$$\ddot{q}_2=\bigg(\frac{||\dot{q}_2||^2}{2}-\int_0^1 \frac{2}{q_2}dt+\frac{\overline{q}_2}{(\overline{q}_1-\overline{q}_2)^2}
-\frac{\dot{q}_2^2}{2}\bigg)\frac{1}{q_2}-\frac{2}{(\overline{q}_1-\overline{q}_2)^2}.$$
Setting
$$\beta_2:=\frac{\ddot{q}_2+\frac{1}{(\overline{q}_1-\overline{q}_2)^2}}{q_2}$$
one gets
$$\dot{\beta}_2 q_2=-3\beta_2 \dot{q}_2$$
implying that there exists $\mu \in \mathbb{R}$ such that
$$\beta_2=-\frac{\mu}{q_2^3}.$$
Arguing as above one deduces from this that $\mu=2$ and thus $q_2$
satisfies the second equation in~\eqref{eq:mean-interaction}.
%\begin{equation}\label{bov10}
%\ddot{q}_2(t)=-\frac{2}{q_2(t)^2}-\frac{1}{(\overline{q}_1-\overline{q}_2)^2}, \quad t \in S^1 \setminus t_{z_2}(Z_2).
%\end{equation}

To see the continuity of $E_1$, we solve equation~\eqref{bov7b} (with
$\mu=2$) for 
$$
  E_1(t) = \frac{\dot q_1(t)^2}{2} - \frac{2}{q_1(t)} 
  = \frac{\|\dot{q}_1\|^2}{2} - \int_0^1 \frac{2}{q_1(s)}ds
  - \frac{\overline{q}_1}{(\overline{q}_1-\overline{q}_2)^2} 
  -\frac{q_1(t)}{(\overline{q}_1-\overline{q}_2)^2}
$$
and note that the right hand side is continuous as a function of $t\in [0,1]$.
Continuity of $E_2$ follows similarly, and we have shown that
$(q_1,q_2)$ is a generalized solution of equation~\eqref{eq:mean-interaction}.

%%%
\subsection{From generalized solutions to critical points}\label{ss:q-to-z-mean}
%%%

Let now $(q_1,q_2)\in H^1(S^1,\R_{\geq 0}\times\R_{\geq 0})$ be a
generalized solution of equation~\eqref{eq:mean-interaction}.
The definition of a generalized solution implies that  
$q_1,q_2\in \HH_{ce}^1(S^1,\R_{\geq 0})$. Corollary~\ref{cor:LC}
implies that the set $\LL^{-1}(q_1)\times \LL^{-1}(q_2)$ consists of
$4$ elements. The goal of this section is to show that each  
$(z_1,z_2)\in \LL^{-1}(q_1)\times \LL^{-1}(q_2)$ is a critical point
of $\B_{av}$. To see this, observe that smoothness of $q_i$ on the
complement of its zero set $Z_{q_i}$ implies smoothness of $z_i$ on the
complement of its zero set $Z_{z_i}$, $i\in \{1,2\}$. In particular,
second derivatives of $z_i$ make sense there and we can  
make the following statement.

\begin{lemma}\label{lem:eqzmean}
Any $(z_1,z_2)\in \LL^{-1}(q_1)\times \LL^{-1}(q_2)$ satisfies the critical point 
equation~\eqref{eq:z-av} on the complement of the set $Z_{z_1}\cup Z_{z_2}$. 
\end{lemma}

Assuming this lemma for the moment, recall that $z_1$ and $z_2$ are of
class $C^1$ and \eqref{eq:z-av} expresses $z_i''$ through
$z_i$. Thus bootstrapping \eqref{eq:z-av} implies that $z_1$ and
$z_2$ are smooth and \eqref{eq:z-av} holds on the whole
$S^1$. Therefore, it remains to prove Lemma~\ref{lem:eqzmean}.
%The proof will occupy the rest of this section. 

{\bf Proof of Lemma~\ref{lem:eqzmean}: }
We will show the desired equation for $z_1$.
%pretending that it has zeros.
A similar argument will do the job for $z_2$. Recall the
equation satisfied by $q_1$,
\begin{equation}\label{eq:q1mean}
\ddot q_1=-\frac{2}{q_1^2}+\frac{1}{(\overline{q}_1-\overline{q}_2)^2}.
\end{equation}
Set
\begin{equation}\label{eq:beta-mean}
\beta_1:=\frac{\ddot{q}_1-\frac{1}{(\overline{q}_1-\overline{q}_2)^2}}{q_1}
\end{equation}
on $S^1\setminus Z_{q_1}$. Then by~\eqref{eq:q1mean} we have
$$
\beta_1=-\frac{2}{q_1^3}, 
$$
and taking time derivative we obtain
$$
\dot\beta_1=\frac{3\cdot 2}{q_1^4}\dot q_1=-\frac{3\beta_1 \dot q_1}{q_1}.
$$
We multiply both sides with $q_1^2$ to get
$$
\dot\beta_1q_1^2=-3\beta_1\dot q_1q_1\,.
$$
We bring $-2\beta_1\dot q_1q_1$ to the other side to continue
$$
\dot\beta_1q_1^2+2\beta_1\dot q_1q_1=-\beta_1\dot q_1q_1\,.
$$
We substitute the original definition~\eqref{eq:beta-mean} 
of $\beta_1$ in the right hand side to get 
$$
\dot\beta_1q_1^2+2\beta_1\dot q_1q_1=-\dot q_1\ddot q_1+
\frac{\dot q_1}{(\overline{q}_1-\overline{q}_2)^2}.
$$
Integrating both sides from $0$ to $t$ we get
$$
  \beta_1 q_1^2
  = C -\frac{\dot{q}_1^2}{2} + \frac{q_1}{(\overline{q}_1-\overline{q}_2)^2}
$$
with some constant $C\in \R$.
We substitute the original definition~\eqref{eq:beta-mean} 
of $\beta_1$ in the left hand side to get 
$$
\ddot{q}_1q_1 = C
-\frac{\dot{q}_1^2}{2}+\frac{2q_1}{(\overline{q}_1-\overline{q}_2)^2}.
$$
Observe that modulo the exact value of $C$ this is exactly equation~\eqref{bov6}, 
which is equivalent the first equation of~\eqref{eq:z-av}.
%Analogous discussion starting with $q_2$
%will lead to the second equation of \eqref{eq:z-av}. 
Therefore, we are left with  computing the constant $C$. 
For this we solve the last equation for $C$ and use~\eqref{eq:q1mean} to obtain
$$
C=\frac{\dot{q}_1^2}{2}-\frac{2}{q_1}-\frac{q_1}{(\overline{q}_1-\overline{q}_2)^2},
$$
and integrating from $0$ to $1$ gives us
$$
  C = \frac{||\dot{q}_1||^2}{2}-\int_0^1
  \frac{2}{q_1(s)}ds-\frac{\overline{q}_1}{(\overline{q}_1-\overline{q}_2)^2},
$$
which matches the constant $c_1$ in~\eqref{eq:c1-mean}. This concludes
the proof of Lemma~\ref{lem:eqzmean}, and therefore of
Theorem~\ref{thm:mean-interaction}. 
\hfill$\square$

%%%%%%%%%%%%%%%%%%%%%%%%%%%%%%%%%%%%%%%%%%%%%%%%%%%%%%%%%%%%%%%%%%%%%%%%%%
\section{Instantaneous interaction}\label{sec:inst-interaction}
%%%%%%%%%%%%%%%%%%%%%%%%%%%%%%%%%%%%%%%%%%%%%%%%%%%%%%%%%%%%%%%%%%%%%%%%%%

In this section we consider the real helium atom in which the two electrons interact by
their Coulomb repulsion according to
\begin{equation}\label{eq:inst-interaction}
\left\{\;
\begin{aligned}
  \ddot{q}_1(t) &= -\frac{2}{q_1(t)^2}+\frac{1}{(q_1(t)-q_2(t))^2}, \cr
  \ddot{q}_2(t) &= -\frac{2}{q_2(t)^2}-\frac{1}{(q_1(t)-q_2(t))^2}
\end{aligned}
\right.
\end{equation}
where we impose the condition
\begin{equation}\label{eq:inst-ineq}
  q_1(t) > q_2(t)\geq 0\qquad\text{for all }t\in S^1. 
\end{equation}

%%%
\subsection{Variational characterization of generalized
  solutions}\label{ss:gen-solutions-inst} 
%%%

Solutions of~\eqref{eq:inst-interaction} are the critical points of the action functional
$$
  \mathcal{S}_{in}(q_1,q_2)
  := \sum_{i=1}^2\Bigl(\frac12\int_0^1\dot q_i(t)^2dt + \int_0^1\frac{2}{q_i(t)}dt\Bigr)
  - \int_0^1\frac{1}{q_1(t)-q_2(t)}dt.
$$
For $i=1,2$ let $q_i$ and $z_i$ be related by Levi-Civita transformations
\begin{equation}\label{eq:LC-zi}
   q_i(t) = z_i(\tau_{z_i}(t))^2,
\end{equation}
where $\tau_{z_i}:S^1\to S^1$ is the inverse of $t_{z_i}:S^1\to S^1$
defined in equation~\eqref{eq:tz} with $z=z_i$. As in Section it
follows that
\begin{equation}\label{eq:t-tau-zi}
   \frac{dt}{q_i(t)} = \frac{d\tau_{z_i}(t)}{\|z_i\|^2},
\end{equation}
so all the relations in Section~\eqref{ss:LC} hold with $(q,z,\tau)=(q_i,z_i,\tau_{z_i})$.
In particular, we can rewrite the last integral in terms of the $z_i$
as the {\em instantaneous interaction}
\begin{align}\label{eq:I}
  \II(z_1,z_2)
  &:= - \frac{1}{\|z_1\|^2}\int_0^1\frac{z_1(\tau)^2}{z_1^2(\tau)-z_2^2
  (\tau_{z_2}(t_{z_1}(\tau)))}d\tau \\
  &= - \int_0^1\frac{1}{z_1^2(\tau_{z_1}(t)) - z_2^2(\tau_{z_2}(t))}dt \cr 
  &= \frac{1}{\|z_2\|^2}\int_0^1\frac{z_2(\tau)^2}{z_2^2(\tau)-z_1^2(\tau_{z_1}(t_{z_2}(\tau)))}d\tau,
\end{align}
where in the last two equalities we have changed the integration
variable to $\tau=\tau_{z_1}(t)$ resp.~$\tau=\tau_{z_2}(t)$ using
equation~\eqref{eq:t-tau-zi}. Note that the functional $\II$ is {\em
  nonlocal} due to the appearance of the time changes, which we have
written as $\tau_{z_i}$ rather than $\tau_i$ to indicate their
dependence on $z_i$.

The instantaneous interaction $\II$ is naturally defined on the space 
\begin{equation}\label{eq:Hin}
\begin{aligned}
  \mathcal{H}_{in}^1 := \Bigl\{
  &z=(z_1,z_2) \in H^1(S^1,\mathbb{R}^2)\;\Bigl|\; 
  \|z_1\|>0,\,\,\|z_2\|>0,\,\, \cr
  &\ \ z_1^2(\tau)-z_2^2(\tau_{z_2}(t_{z_1}(\tau)))>0 \,\,\text{for
    all}\,\,\tau\in S^1\Bigr\}.
\end{aligned}
\end{equation}
Note that $\mathcal{H}_{in}^1$ is an open subset of the Hilbert space
$H^1(S^1,\mathbb{R}^2)$ and the last condition in its definition
corresponds to condition~\eqref{eq:inst-ineq}. Since integrating
condition~\eqref{eq:inst-ineq} leads to the averaged
condition~\eqref{eq:mean-ineq}, which is in turn equivalent to the
last condition in the definition of $\mathcal{H}_{av}^1$ in
Section~\ref{ss:var-mean}, we have
\begin{equation}\label{eq:inavincl}
  \mathcal{H}_{in}^1\subset \mathcal{H}_{av}^1.
\end{equation}
On $\mathcal{H}_{in}^1$ we consider the functional
\begin{equation}\label{eq:Bin}
  \mathcal{B}_{in} \colon \mathcal{H}_{in}^1 \to \mathbb{R}, \qquad
  \mathcal{B}_{in}(z_1,z_2) := \QQ(z_1) + \QQ(z_2) + \mathcal{I}(z_1,z_2),
\end{equation}
with the functionals
\begin{equation*}%\label{eq:Q}
  \QQ(z_i) = 2\|z_i\|^2\|z_i'\|^2 + \frac{2}{\|z_i\|^2}
\end{equation*}
from equation~\eqref{eq:Q} with charge $N=2$. 

We call $(q_1,q_2)\in H^1(S^1,\R_{\geq 0}\times\R_{\geq 0})$ a {\em
  generalized solution of~\eqref{eq:inst-interaction}} if for $i=1,2$
the following holds:
\begin{enumerate}
\item $q_1(t)>q_2(t)$ for all $t\in S^1$;
\item the zero sets $Z_i=q_i^{-1}(0)\subset S^1$ are finite and each have an even
  number of elements;
\item the restricted maps $q_i:S^1\setminus Z_i\to\R_{\geq 0}$ are
  smooth and satisfy~\eqref{eq:inst-interaction}; 
\item the energies
  $$
    E_i(t) := \frac{\dot q_i(t)^2}{2}-\frac{2}{q_i(t)},\qquad t\in
    S^1\setminus Z_i
  $$
extend to continuous functions $E_i:S^1\to\R$. 
\end{enumerate}
Note that the individual energies $E_i$ need not be constant, but
the total energy
$$
   E = E_1(t) + E_2(t) + \frac{1}{q_1(t)-q_2(t)}
$$
is constant and negative. 

\begin{thm}[Generalized solutions with instantaneous
    interaction]\label{thm:inst-interaction} 
Under the Levi-Civita transformations~\eqref{eq:LC-zi} with time
changes~\eqref{eq:t-tau-zi}, critical points $(z_1,z_2)$ of the action
functional $\mathcal{B}_{in}$ defined in~\eqref{eq:Bin} are in 4-to-1 correspondence
with generalized solutions $(q_1,q_2)$ of~\eqref{eq:inst-interaction}.
\end{thm}

The proof of this theorem will take up the remainder of this section.

%%%
\subsection{The differential of $\mathcal{B}_{in}$}\label{ss:diffin}
%%%

In this subsection we compute the differential of $\mathcal{B}_{in}$
at $(z_1,z_2)\in\mathcal{H}_{in}^1$ in direction $(v_1,v_2)\in H^1(S^1,\R^2)$.
For this we will need for $i=1,2$ the derivative of the time change
$t_{z_i}$ with respect to $z_i$ in direction
$v_i$. By~\eqref{eq:dt} it is given by
\begin{equation}\label{eq:dt-i}
   Dt_{z_i}(v_i)(\tau) = \frac{2}{\|z_i\|^2}\int_0^\tau z_i(\sigma) v_i(\sigma) d\sigma
   - \frac{2\langle z_i,v_i\rangle}{\|z_i\|^4}\int_0^\tau z_i(\sigma)^2 d\sigma\,. 
\end{equation}
%where as in Appendix~\ref{} we use the notation
%\begin{eqnarray*}
%  \langle z_i, v_i \rangle(\tau) &=& \int_0^\tau z_i(\sigma) v_i(\sigma) d\sigma\\
%  \|z_i\|^2(\tau) &=& \langle z_i, z_i\rangle(\tau) = \int_0^\tau z_i(\sigma)^2 d\sigma. 
%\end{eqnarray*}
Using this we now compute the differential of the instantaneous interaction
$\mathcal{I}$ at $(z_1,z_2)\in\mathcal{H}_{in}^1$ with respect to
$z_1$ in direction $v_1\in H^1(S^1,\R)$. Using the first expression
for $\I$ in~\eqref{eq:I} we obtain
\begin{eqnarray*}
  & & D_1\I[z_1,z_2](v_1)\\
  &=& \frac{2\langle z_1,v_1\rangle}{||z_1||^4}
  \int_0^1\frac{z_1^2(\tau)}{z_1^2(\tau)-z_2^2(\tau_{z_2}(t_{z_1}(\tau)))}d\tau\\ 
  & & +\frac{2}{||z_1||^2}\int_0^1 \frac{z_1(\tau)
    z_2^2(\tau_{z_2}(t_{z_1}(\tau)))v_1(\tau)}{\big(z_1^2(\tau) -
    z_2^2(\tau_{z_2}(t_{z_1}(\tau)))\big)^2}d\tau\\  
  & & -\frac{2}{||z_1||^2}\int_0^1 \frac{z_1^2(\tau)
    z_2(\tau_{z_2}(t_{z_1}(\tau)))z_2'(\tau_{z_2}(t_{z_1}(\tau)))\dot{\tau}_{z_2}(t_{z_1}(\tau))
    Dt_{z_1}(v_1)(\tau)}{\big(z_1^2(\tau)-z_2^2(\tau_{z_2}(t_{z_1}(\tau)))\big)^2}d\tau
\end{eqnarray*}
We rewrite the third term on the right hand side as
\begin{eqnarray*}
  & & -\frac{2}{||z_1||^2}\int_0^1 \frac{z_1^2(\tau)
    z_2(\tau_{z_2}(t_{z_1}(\tau)))z_2'(\tau_{z_2}(t_{z_1}(\tau)))\dot{\tau}_{z_2}(t_{z_1}(\tau))
    Dt_{z_1}(v_1)(\tau)}{\big(z_1^2(\tau)-z_2^2(\tau_{z_2}(t_{z_1}(\tau)))\big)^2}d\tau\\
  &=& -\frac{2||z_2||^2}{||z_1||^2}\int_0^1 \frac{z_1^2(\tau)
    z_2'(\tau_{z_2}(t_{z_1}(\tau)))Dt_{z_1}(v_1)(\tau)}
    {z_2(\tau_{z_2}(t_{z_1}(\tau)))\big(z_1^2(\tau)-z_2^2(\tau_{z_2}(t_{z_1}(\tau)))\big)^2}d\tau\\ 
  &=& -\frac{4||z_2||^2}{||z_1||^4}\int_0^1 \frac{z_1^2(\tau)
    z_2'(\tau_{z_2}(t_{z_1}(\tau)))}{z_2(\tau_{z_2}(t_{z_1}(\tau)))\big(z_1^2(\tau)
      -z_2^2(\tau_{z_2}(t_{z_1}(\tau)))\big)^2}
    \bigg(\int_0^\tau z_1(\sigma)v_1(\sigma) d\sigma\bigg) d\tau\\ 
  & & +\frac{4||z_2||^2\langle z_1,v_1\rangle}{||z_1||^6}\int_0^1
    \frac{z_1^2(\tau) z_2'(\tau_{z_2}(t_{z_1}(\tau)))\int_0^\tau z_1(\sigma)^2d\sigma
    }{z_2(\tau_{z_2}(t_{z_1}(\tau)))\big(z_1^2(\tau)
      -z_2^2(\tau_{z_2}(t_{z_1}(\tau)))\big)^2} d\tau\\ 
  &=& -\frac{4||z_2||^2}{||z_1||^4}\int_0^1\bigg(\int_\sigma^1
    \frac{z_1^2(\tau)z_2'(\tau_{z_2}(t_{z_1}(\tau)))}{z_2(\tau_{z_2}(t_{z_1}(\tau)))
    \big(z_1^2(\tau)-z_2^2(\tau_{z_2}(t_{z_1}(\tau)))\big)^2}d\tau\bigg)
    z_1(\sigma) v_1(\sigma) d\sigma\\ 
  & & +\frac{4||z_2||^2\langle z_1,v_1\rangle}{||z_1||^4}\int_0^1
    \frac{z_1^2(\tau) z_2'(\tau_{z_2}(t_{z_1}(\tau)))
      t_{z_1}(\tau)}{z_2(\tau_{z_2}(t_{z_1}(\tau)))\big(z_1^2(\tau)
      -z_2^2(\tau_{z_2}(t_{z_1}(\tau)))\big)^2}d\tau.\\ 
\end{eqnarray*}
Here for the first equality we haved used formula~\eqref{eq:taudot} for
$\dot\tau_{z_i}$, and for the second one we have used equation~\eqref{eq:dt-i}.
For the third equality we have switched the order of integration in the
first integral, and we have used equation~\eqref{eq:ti} to replace
$\int_0^\tau z_1(\sigma)^2d\sigma$ by $\|z_1\|^2t_{z_1}(\tau)$. 

By symmetry, the expression for $D_2\I[z_1,z_2](v_2)$ is the same with
a global minus sign and the roles of $(z_1,v_1)$ and $(z_2,v_2)$ reversed.
Putting everything together,
%and renaming the intergration variables $\sigma,\tau$ in the double
%integral,
we obtain the differential of $\B_{in}$:
 \begin{eqnarray}\label{eq:DBin} \nonumber
  & &d\B_{in}[z_1,z_2](v_1,v_2)\\ \nonumber
  &=&4 \sum_{i=1}^2 \Bigg(||z_i'||^2 \cdot \langle z_i,v_i\rangle
+||z_i||^2 \langle z_i',v_i'\rangle-\frac{\langle z_i,v_i\rangle}{||z_i||^4}\Bigg)\\ \nonumber
  & &+\frac{2\langle z_1,v_1\rangle}{||z_1||^4}\int_0^1\frac{z_1^2(\tau)}{z_1^2(\tau)
  -z_2^2(\tau_{z_2}(t_{z_1}(\tau)))}d\tau\\ \nonumber
  & &+\frac{2}{||z_1||^2}\int_0^1 \frac{z_1(\tau)
  z_2^2(\tau_{z_2}(t_{z_1}(\tau)))v_1(\tau)}{\big(z_1^2(\tau)
  -z_2^2(\tau_{z_2}(t_{z_1}(\tau)))\big)^2}d\tau\\ \nonumber
  & &-\frac{4||z_2||^2}{||z_1||^4}\int_0^1\Bigg(\int_\sigma^1 \frac{z_1^2(\tau) z_2'(\tau_{z_2}(t_{z_1}(\tau)))}{z_2(\tau_{z_2}(t_{z_1}(\tau)))\big(z_1^2(\tau)-z_2^2(\tau_{z_2}(t_{z_1}(\tau)))\big)^2}d\tau
\Bigg) z_1(\sigma) v_1(\sigma) d\sigma\\ 
& &+\frac{4||z_2||^2\langle z_1,v_1\rangle}{||z_1||^4}\int_0^1
\frac{z_1^2(\tau) z_2'(\tau_{z_2}(t_{z_1}(\tau)))
  t_{z_1}(\tau)}{z_2(\tau_{z_2}(t_{z_1}(\tau)))\big(z_1^2(\tau)
  -z_2^2(\tau_{z_2}(t_{z_1}(\tau)))\big)^2}d\tau\\ \nonumber
& &-\frac{2\langle
  z_2,v_2\rangle}{||z_2||^4}\int_0^1\frac{z_2^2(\tau)}
     {z_2^2(\tau)-z_1^2(\tau_{z_1}(t_{z_2}(\tau)))}d\tau\\ \nonumber
& &-\frac{2}{||z_2||^2}\int_0^1 \frac{z_2(\tau)
       z_1^2(\tau_{z_1}(t_{z_2}(\tau)))v_2(\tau)}
     {\big(z_2^2(\tau)-z_1^2(\tau_{z_1}(t_{z_2}(\tau)))\big)^2}d\tau\\ \nonumber
& &+\frac{4||z_1||^2}{||z_2||^4}\int_0^1\Bigg(\int_\sigma^1 \frac{z_2^2(\tau) z_1'(\tau_{z_1}(t_{z_2}(\tau)))}{z_1(\tau_{z_1}(t_{z_2}(\tau)))\big(z_2^2(\tau)-z_1^2(\tau_{z_1}(t_{z_2}(\tau)))\big)^2}d\tau\Bigg)
z_2(\sigma) v_2(\sigma)d\sigma\\ \nonumber
& &-\frac{4||z_1||^2\langle z_2,v_2\rangle}{||z_2||^4}\int_0^1
\frac{z_2^2(\tau) z_1'(\tau_{z_1}(t_{z_2}(\tau)))
  t_{z_2}(\tau)}{z_1(\tau_{z_1}(t_{z_2}(\tau)))\big(z_2^2(\tau)
  -z_1^2(\tau_{z_1}(t_{z_2}(\tau)))\big)^2}d\tau\\ \nonumber 
\end{eqnarray} 
 
%%%
\subsection{Critical points of $\mathcal{B}_{in}$}
%%%

%The goal of this part is to describe critical points of $\mathcal{B}_{in}$ in terms of a 
%(nonlocal nonlinear coupled integral-differential) system of equations they must satisfy. Then we use this system to derive certain properties of critical points.
Equation~\eqref{eq:DBin} leads to the characterization of critical
points of $\mathcal{B}_{in}$:

\begin{prop}\label{prop:critin}
A point $(z_1,z_2)\in \mathcal{H}_{in}^1$ is a critical point of $\mathcal{B}_{in}$ if and only if 
$(z_1,z_2)$ is smooth and solves the following system of coupled
nonlocal integral--differential equations:
\begin{eqnarray}\label{eq:critinst}
z_1''(\tau)&=&\frac{||z_1'||^2z_1(\tau)}{||z_1||^2}-\frac{z_1(\tau)}{||z_1||^6}\\ \nonumber
& &+\frac{z_1(\tau)}{2||z_1||^6}\int_0^1\frac{z_1^2(\sigma)}{z_1^2(\sigma)-z_2^2(\tau_{z_2}(t_{z_1}(\sigma)))}d\sigma\\ \nonumber
& &+\frac{z_1(\tau)z_2^2(\tau_{z_2}(t_{z_1}(\tau)))}{2||z_1||^4\big(z_1^2(\tau)-z_2^2(\tau_{z_2}(t_{z_1}(\tau)))\big)^2}\\ \nonumber
& &-\frac{||z_2||^2z_1(\tau)}{||z_1||^6}
\int_\tau^1 \frac{z_1^2(\sigma) z_2'(\tau_{z_2}(t_{z_1}(\sigma)))}{z_2(\tau_{z_2}(t_{z_1}(\sigma)))\big(z_1^2(\sigma)-z_2^2(\tau_{z_2}(t_{z_1}(\sigma)))\big)^2}d\sigma\\ \nonumber
& &+\frac{||z_2||^2 z_1(\tau)}{||z_1||^6}
\int_0^1 \frac{z_1^2(\sigma) z_2'(\tau_{z_2}(t_{z_1}(\sigma))) t_{z_1}(\sigma)}
{z_2(\tau_{z_2}(t_{z_1}(\sigma)))\big(z_1^2(\sigma)-z_2^2(\tau_{z_2}(t_{z_1}(\sigma)))\big)^2}d\sigma\\
\nonumber
z_2''(\tau)&=&\frac{||z_2'||^2z_2(\tau)}{||z_2||^2}-\frac{z_2(\tau)}{||z_2||^6}\\ \nonumber
& &-\frac{z_2(\tau)}{2||z_2||^6}\int_0^1\frac{z_2^2(\sigma)}{z_2^2(\sigma)-z_1^2(\tau_{z_1}(t_{z_2}(\sigma)))}d\sigma\\ \nonumber
& &-\frac{z_2(\tau)z_1^2(\tau_{z_1}(t_{z_2}(\tau)))}{2||z_2||^4\big(z_2^2(\tau)-z_1^2(\tau_{z_1}(t_{z_2}(\tau)))\big)^2}\\ \nonumber
& &+\frac{||z_1||^2z_2(\tau)}{||z_2||^6}
\int_\tau^1 \frac{z_2^2(\sigma) z_1'(\tau_{z_1}(t_{z_2}(\sigma)))}{z_1(\tau_{z_1}(t_{z_2}(\sigma)))\big(z_2^2(\sigma)-z_1^2(\tau_{z_1}(t_{z_2}(\sigma)))\big)^2}d\sigma\\ \nonumber
& &-\frac{||z_1||^2 z_2(\tau)}{||z_2||^6}
\int_0^1 \frac{z_2^2(\sigma) z_1'(\tau_{z_1}(t_{z_2}(\sigma))) t_{z_2}(\sigma)}{z_1(\tau_{z_1}(t_{z_2}(\sigma)))\big(z_2^2(\sigma)-z_1^2(\tau_{z_1}(t_{z_2}(\sigma)))\big)^2}d\sigma
\end{eqnarray}
\end{prop}

\textbf{Proof: } 
From equation~\eqref{eq:DBin} we see that $(z_1,z_2)$ is a critical
point of $\mathcal{B}_{in}$ if and only if $z_1$ and $z_2$ have
weak second derivatives and satisfy the system of
equations~\eqref{eq:critinst}. Bootstrapping these equations
%using the Sobolev embedding $W^{k,2}\rightarrow C^k(R)$ for $k=1,2\dots$ 
we conclude that $z_1$ and $z_2$ are smooth and the proposition follows.
\hfill $\square$

\begin{cor}\label{cor:inst-transverse-zeroes}
Suppose that $(z_1,z_2)$ is a critical point of
$\mathcal{B}_{in}$. Then $z_1$ has no zeroes and $z_2$ has transverse
zeros. In particular, their zero sets
$$
  Z_i := \big\{\tau \in S^1 \mid z_i(\tau)=0\big\}, \qquad i=1,2
$$
are finite.
\end{cor}

\textbf{Proof: } 
Let $(z_1,z_2)\in\mathcal{H}_{in}^1$ be a critical point of
$\mathcal{B}_{in}$. It follows directly from the definition of
$\mathcal{H}_{in}^1$ that $z_1(t)>0$ for all $t\in S^1$. 
Next note that at $(z_1,z_2)$ solves the system of
coupled nonlinear integral--differential
equations~\eqref{eq:critinst} which has the form
\begin{equation}\label{eq:critinst1}
\left\{
\begin{aligned}
  z_1''(\tau) &= f_1(\tau)z_1(\tau)\\
  z_2''(\tau) &= f_2(\tau)z_2(\tau),
\end{aligned}
\right.
\end{equation}
where $f_i:S^1\to \R$, $i=1,2$ are smooth functions depending on 
$(z_1,z_2)$. Disregarding the dependence of the $f_i$ on $(z_1,z_2)$,
we can view $(z_1,z_2)$ as a solution to the system of {\em uncoupled
  linear ODEs}~\eqref{eq:critinst1}.

Arguing by contradiction, suppose now that there exists a point
$\tau_0\in S^1$ such that $z_2(\tau_0)=z_2'(\tau_0)=0$. 
Then the function $z_2$ and the zero function both solve the second
equation in~\eqref{eq:critinst1} with the same initial conditions at
$\tau_0$. By uniqueness of solutions of ODEs we conclude $z_2\equiv
0$, contradicting the condition $\|z_2\|>0$ in the definition of
$\mathcal{H}_{in}$. %An analogous argument applies to $z_2$.
\hfill $\square$

%%%
\subsection{From critical points to generalized solutions}\label{ss:z-to-q-inst}
%%%

Let $(z_1,z_2)\in \mathcal{H}_{in}^1$ be a critical point of $\mathcal{B}_{in}$,
so by Proposition~\ref{prop:critin} the maps $z_1,z_2:S^1\to\R$ are
smooth and satisfy~\eqref{eq:critinst}. 
As in the previous section, for $i=1,2$ we define the smooth maps
$t_{z_i}:S^1\to S^1$ by~\eqref{eq:ti}. 
Since by Corollary~\ref{cor:inst-transverse-zeroes} the map $z_i$ has
only finitely many zeroes, it follows from Lemma~\ref{lem:tz-homeo}
that $t_{z_i}:S^1\to S^1$ is a homeomorphism with continuous inverse 
$\tau_{z_i} \colon S^1 \to S^1$. We define continuous maps $q_i:S^1\to S^1$ by
\begin{equation}\label{eq:q-i-z-i}
  q_i(t) := z_i(\tau_{z_i}(t))^2.
\end{equation}
Then for $i=1,2$ the maps $z_i,q_i:S^1\to\R$ are smooth except at finitely many
points and related by the Levi-Civita transformation~\eqref{eq:LC}
with $\tau=\tau_{z_i}$. The
last condition in the definition of $\mathcal{H}_{in}^1$ implies
$$
  q_1(t) > q_2(t) \quad\text{for all }t\in S^1. 
$$
Let us now focus on $i=1$. Although by Corollary~\ref{cor:inst-transverse-zeroes}
the function $z_1$ has no zeroes, in the following argument we
will allow $z_1$ to have a finite set $Z_1$ of zeroes; this will
ensure that the same argument carries over to $z_2$ (which may have zeroes).
At points $t\in S^1\setminus t_{z_1}(Z_1)$ we compute:
\begin{eqnarray}\label{qinst}
  \ddot{q}_1(t) q_1(t) \nonumber
  &=& 2\|z_1\|^4 \frac{z_1''(\tau_{z_1}(t))}{z_1(\tau_{z_1}(t))}
  -\frac{\dot{q}_1(t)^2}{2}\\ \nonumber
  &=& 2||z_1||^2\cdot ||z_1'||^2
  -\frac{2}{||z_1||^2}
  -\frac{\dot{q}_1^2(t)}{2}\\ \nonumber
  & & +\frac{1}{||z_1||^2}\int_0^1\frac{z_1(\tau_{z_1}(s))^2}
  {z_1^2(\tau_{z_1}(s))-z_2^2(\tau_{z_2}(s))}\dot{\tau}_{z_1}(s)ds\\ \nonumber 
  & &+\frac{z_2^2(\tau_{z_2}(t))}{\big(z_1^2(\tau_{z_1}(t))-z_2^2(\tau_{z_2}(t))\big)^2}\\ \nonumber
  & &-\frac{2||z_2||^2}{||z_1||^2}\int_t^1 \frac{z_1^2(\tau_{z_1}(s))z_2'(\tau_{z_2}(s))}{z_2(\tau_{z_2}(s))
\big(z_1^2(\tau_{z_1}(s))-z_2^2(\tau_{z_2}(s))\big)^2}\dot{\tau}_{z_1}(s)ds\\ \nonumber
  & &+\frac{2||z_2||^2}{||z_1||^2}\int_0^1
  \frac{z_1^2(\tau_{z_1}(s))z_2'(\tau_{z_2}(s))s}{z_2(\tau_{z_2}(s)) 
\big(z_1^2(\tau_{z_1}(s))-z_2^2(\tau_{z_2}(s))\big)^2}\dot{\tau}_{z_1}(s)ds\\ \nonumber
  &=& \frac{||\dot{q}_1||^2}{2}
  -\int_0^1 \frac{2}{q_1(s)}ds
  -\frac{\dot{q}_1^2(t)}{2}
  +\int_0^1\frac{1}{q_1(s)-q_2(s)}ds\\ \nonumber
  & &+\frac{q_2(t)}{\big(q_1(t)-q_2(t)\big)^2}
  -2||z_2||^2 \int_t^1 \frac{\frac{d}{ds}z_2(\tau_{z_2}(s))}{\dot{\tau}_{z_2}(s)
  z_2(\tau_{z_2}(s))\big(q_1(s)-q_2(s)\big)^2}ds\\ \nonumber
  & &+2||z_2||^2 \int_0^1 \frac{s\frac{d}{ds}z_2(\tau_{z_2}(s))}{\dot{\tau}_{z_2}(s)
  z_2(\tau_{z_2}(s))\big(q_1(s)-q_2(s)\big)^2}ds\\ \nonumber
  &=& \frac{||\dot{q}_1||^2}{2}
  -\int_0^1 \frac{2}{q_1(s)}ds
  -\frac{\dot{q}_1^2(t)}{2}
  +\int_0^1\frac{1}{q_1(s)-q_2(s)}ds\\ \nonumber
  & &+\frac{q_2(t)}{\big(q_1(t)-q_2(t)\big)^2}
  -\int_t^1 \frac{2z_2(\tau_{z_2}(s))\frac{d}{ds}z_2(\tau_{z_2}(s))}{\big(q_1(s)-q_2(s)\big)^2}ds\\ \nonumber 
  & &+\int_0^1 \frac{2sz_2(\tau_{z_2}(s))\frac{d}{ds}z_2(\tau_{z_2}(s))}{\big(q_1(s)-q_2(s)\big)^2}ds\\ \nonumber
  &=& \frac{||\dot{q}_1||^2}{2}
  -\int_0^1 \frac{2}{q_1(s)}ds
  -\frac{\dot{q}_1^2(t)}{2}
  +\int_0^1\frac{1}{q_1(s)-q_2(s)}ds\\ \nonumber
  & &+\frac{q_2(t)}{\big(q_1(t)-q_2(t)\big)^2}
  -\int_t^1 \frac{\frac{d}{ds}z_2^2(\tau_{z_2}(s))}
  {\big(q_1(s)-q_2(s)\big)^2}ds\\ \nonumber
  & &+\int_0^1 \frac{s\frac{d}{ds}z_2^2(\tau_{z_2}(s))}{\big(q_1(s)-q_2(s)\big)^2}ds\\ \nonumber
  &=& \frac{||\dot{q}_1||^2}{2}
  -\int_0^1 \frac{2}{q_1(s)}ds
  -\frac{\dot{q}_1^2(t)}{2}
  +\int_0^1\frac{1}{q_1(s)-q_2(s)}ds\\ \nonumber
  & &+\frac{q_2(t)}{\big(q_1(t)-q_2(t)\big)^2}
  -\int_t^1 \frac{\dot{q}_2(s)}{\big(q_1(s)-q_2(s)\big)^2}ds\\ 
  & &+\int_0^1 \frac{s\dot{q}_2(s)}{\big(q_1(s)-q_2(s)\big)^2}ds\,.
\end{eqnarray}
Here the first equality comes from equation~\eqref{eq:ddot-q} with $q=q_1$ and $z=z_1$. 
In the second one we substitute $z_1''$ by~\eqref{eq:critinst} and change
variables $\sigma=\tau_{z_1}(s)$ in the integrals. 
In the third one we use the following replacements from~\eqref{eq:q-i-z-i},
\eqref{bov2}, \eqref{bov4} and~\eqref{eq:taudot}:
\begin{gather*}
  z_i(\tau_{z_i}(s))^2 = q_i(s),\qquad
  \frac{1}{||z_1||^2} = \int_0^1 \frac{1}{q_1(s)}ds,\cr
  4||z_1||^2 \cdot ||z_1'||^2 = ||\dot{q}_1||^2,\qquad
  \dot\tau_{z_1}(s) = \frac{\|z_1\|^2}{z_1(\tau_{z_1}(s))^2}
\end{gather*}
and the chain rule for $\dds z_2(\tau_{z_2}(s))$. 
In the fourth one we use~\eqref{eq:taudot} to replace $\dot\tau_{z_2}(s)$,
in the fifth one the chain rule for $\dds{z_2^2(\tau_{z_2}(s))}$,
and in the sixth one~\eqref{eq:q-i-z-i} to insert $q_2(s)$. 
Thus $q_1$ satisfies the integral--differential equation
\begin{equation}\label{qinst2}
  \ddot{q}_1 = \bigg(c_1
  -\frac{\dot{q}_1^2}{2}
  +\frac{q_2}{\big(q_1-q_2\big)^2}
  -\int_t^1 \frac{\dot{q}_2(s)}{\big(q_1(s)-q_2(s)\big)^2}ds
  \bigg)\frac{1}{q_1}
\end{equation}
with the constant
$$
  c_1 = \frac{||\dot{q}_1||^2}{2}
  -\int_0^1 \frac{2}{q_1(s)}ds
  +\int_0^1\frac{1}{q_1(s)-q_2(s)}ds
  +\int_0^1 \frac{s\,\dot{q}_2(s)}{\big(q_1(s)-q_2(s)\big)^2}ds\,.
$$
Let now $t_-<t_+$ be adjacent zeroes of $q_1$ and consider the smooth map
$$
  \beta_1:=\frac{\ddot{q}_1-\frac{1}{(q_1-q_2)^2}}{q_1}
  \colon (t_-,t_+) \to \mathbb{R}.
$$
From~\eqref{qinst2} we obtain
\begin{eqnarray*}
  \beta_1 q_1^2
  &=& \ddot{q}_1q_1 - \frac{q_1}{(q_1-q_2)^2}\\
%  &=&\frac{||\dot{q}_1||^2}{2}-\int_0^1 \frac{2}{q_1}ds-\frac{\dot{q}_1^2}{2}+\int_0^1\frac{1}{q_1-
%q_2}ds-\int_t^1 \frac{\dot{q}_2}{\big(q_1-q_2\big)^2}ds\\
%& &+\int_0^1 \frac{s\dot{q}_2}{\big(q_1-q_2\big)^2}ds+\frac{q_2-q_1}{(q_1-q_2)^2}\\
%&=&\frac{||\dot{q}_1||^2}{2}-\int_0^1 \frac{2}{q_1}ds-\frac{\dot{q}_1^2}{2}+\int_0^1\frac{1}{q_1-
%q_2}ds-\int_t^1 \frac{\dot{q}_2}{\big(q_1-q_2\big)^2}ds\\
%& &+\int_0^1 \frac{s\dot{q}_2}{\big(q_1-q_2\big)^2}ds+\frac{1}{q_2-q_1}
%&=& c_1 -\frac{\dot{q}_1^2}{2}
% +\frac{q_2-q_1}{\big(q_1-q_2\big)^2}
%  -\int_t^1 \frac{\dot{q}_2(s)}{\big(q_1(s)-q_2(s)\big)^2}ds\\
&=& c_1 -\frac{\dot{q}_1^2}{2}
  +\frac{1}{q_2-q_1}
  -\int_t^1 \frac{\dot{q}_2(s)}{\big(q_1(s)-q_2(s)\big)^2}ds\,.
\end{eqnarray*}
Taking the time derivative of this expression for $\beta_1 q_1^2$ we obtain
\begin{eqnarray*}
  \dot{\beta}_1 q_1^2+2\beta_1 q_1\dot{q}_1
  &=&-\ddot{q}_1 \dot{q}_1 -\frac{\dot{q}_2-\dot{q}_1}{(q_2-q_1)^2}
  +\frac{\dot{q}_2}{(q_1-q_2)^2}\\
&=&-\ddot{q}_1 \dot{q}_1+\frac{\dot{q}_1}{(q_1-q_2)^2}\\
&=&-\beta_1 q_1 \dot{q}_1\,.
\end{eqnarray*}
Therefore, dividing both sides by $q_1$, we get 
$$
\dot\beta_1q_1=-3\beta_1\dot q_1\,.
$$
This is exactly equation~\eqref{eq:betaqu}. We apply Lemma~\ref{lem:betamu} to get
$$
  \beta_1 = -\frac{\mu}{q_1^3} 
$$
on $(t_-,t_+)$ for some constant $\mu>0$. Here positivity of $\mu$ follows again
because by definition $\beta_1<0$ at the maximum of $q_1$ on $(t_-,t_+)$. 
We substitute the definition of $\beta_1$ in the last displayed equation and solve for $\ddot{q}_1$
to get
\begin{equation}\label{eq:ddotqmu}
\ddot{q}_1(t)=-\frac{\mu}{q_1(t)^2}+\frac{1}{(q_1(t)-q_2(t))^2}, \qquad t \in (t_-,t_+).
\end{equation}
We solve the last equation for $\mu/q_1$ and substitute $\ddot q_1 q_1$
by~\eqref{qinst} to get
\begin{eqnarray}\label{eq:muq}
  \frac{\mu}{q_1(t)} \nonumber
  &=& -\ddot q_1(t)q_1(t) + \frac{q_1(t)}{(q_1(t)-q_2(t))^2}\\ \nonumber
  &=& -\frac{||\dot q_1 ||^2}{2}+\int_0^1\frac{2ds}{q_1(s)}+\frac{\dot q_1(t)^2}{2}
  -\int_0^1\frac{ds}{q_1(s)-q_2(s)}+\frac{1}{q_1(t)-q_2(t)}\\ 
  & & +\int_t^1\frac{\dot q_2(s)ds}{(q_1(s)-q_2(s))^2}
  -\int_0^1\frac{s\,\dot q_2(s)ds}{(q_1(s)-q_2(s))^2}\,.
\end{eqnarray}
We multiply both sides of the last equation by $q_1(t)$
take the limit $t\to t_{\pm}$ and recall that $q_1(t_\pm)=0$ to get 
\begin{equation}%\label{eq:muinst}
  \mu = \lim_{t \to t_\pm} \frac{\dot{q}_1(t)^2 q_1(t)}{2}=2||z_1||^4 z_1'\big(\tau_{z_1}(t_\pm)\big)^2.
\end{equation}
Therefore, $\mu$ is a global constant independent of the interval $(t_-,t_+)$.
(This argument is actually only needed for $q_2$ in place of $q_1$;
for $q_1$ itself, since it has no zeroes, we can replace $(t_-,t_+)$
by all of $S^1$ from the outset). 

Now note that by Fubini's theorem for every integrable function
$f:[0,1]\to\R$ we have 
$$
  \int_0^1dt\int_t^1f(s)ds = \int_0^1ds\,f(s)\int_0^s dt = \int_0^1sf(s)ds\,.
$$
Integrating both sides of Equation~\eqref{eq:muq} from $0$ to $1$ and
applying this identity to 
$$
  f(s):=\frac{\dot q_2}{(q_1(s)-q_2(s))^2}
$$
we obtain
$$
\mu=2\,.
$$
Thus equation~\eqref{eq:ddotqmu} becomes the first equation in~\eqref{eq:inst-interaction}.
Similarly, we deduce that $q_2$ satisfies the second equation
in~\eqref{eq:inst-interaction} outside its zero set. 

To see the continuity of the energy $E_1$, we solve equation~\eqref{eq:muq} (with
$\mu=2$) for 
\begin{eqnarray*}\nonumber
  E_1(t)
  &=& \frac{\dot q_1(t)^2}{2} - \frac{2}{q_1(t)}\\ 
  &=& \frac{||\dot q_1 ||^2}{2} - \int_0^1\frac{2\,ds}{q_1(s)}
  + \int_0^1\frac{ds}{q_1(s)-q_2(s)} - \frac{1}{q_1(t)-q_2(t)}\\ \nonumber
  & & - \int_t^1\frac{\dot q_2(s)ds}{(q_1(s)-q_2(s))^2}
  + \int_0^1\frac{s\,\dot q_2(s)ds}{(q_1(s)-q_2(s))^2}
\end{eqnarray*}
and note that the right hand side is continuous as a function of $t\in [0,1]$.
Continuity of $E_2$ follows similarly, and we have shown that
$(q_1,q_2)$ is a generalized solution of equation~\eqref{eq:inst-interaction}.

%%%
\subsection{From generalized solutions to critical points}\label{ss:q-to-z-inst}
%%%

Let now $(q_1,q_2)\in H^1(S^1,\R_{\geq 0}\times\R_{\geq 0})$ be a
generalized solution of equation~\eqref{eq:inst-interaction}. The
definition of a generalized solution implies that $q_1,q_2\in
\HH_{ce}^1(S^1,\R_{\geq 0})$. Corollary~\ref{cor:LC} implies that the
set $\LL^{-1}(q_1)\times \LL^{-1}(q_2)$ 
consists of $4$ elements. The goal of this section is to show that each 
$(z_1,z_2)\in \LL^{-1}(q_1)\times \LL^{-1}(q_2)$ is a critical point
of $\B_{in}$. To see this, recall that $q_1$ has no zeroes, so
smoothness of $q_1$ implies smoothness of $z_1$. Smoothness of $q_2$ on the
complement of its zero set $Z_{q_2}$ implies smoothness of $z_2$ on
the complement of its zero set $Z_{z_2}$. In
particular, second derivatives of $z_2$ make sense there and we can  
make the following statement.

\begin{lemma}\label{lem:eqzinst}
Any $(z_1,z_2)\in \LL^{-1}(q_1)\times \LL^{-1}(q_2)$ satisfies the critical point 
equation \eqref{eq:critinst} on the complement of the set $Z_{z_2}$. 
\end{lemma}

Assuming this lemma for the moment, recall from Corollary~\ref{cor:LC} that $z_2$ is of class $C^1$ and that
\eqref{eq:critinst} expresses $z_2''$ through $z_1,\, z_1',\, z_2,\,z_2'$.
%and similarly $z_2''$. 
Thus bootstrapping \eqref{eq:critinst} implies that $z_2$
is smooth and \eqref{eq:critinst} holds on the whole $S^1$. Therefore,
it remains to prove  Lemma~\ref{lem:eqzinst}. 

{\bf Proof of Lemma~\ref{lem:eqzinst}: }
We will show the desired equation for $z_1$ pretending that it has
a possibly nonempty zero set $Z_{q_1}$. A similar argument will do the
job for $z_2$. Recall the equation satisfied by $q_1$, 
\begin{equation}\label{eq:q1}
\ddot q_1=-\frac{2}{q_1^2}+\frac{1}{(q_1-q_2)^2}.
\end{equation}
Set
\begin{equation}\label{eq:beta-inst}
\beta_1:=\frac{\ddot{q}_1-\frac{1}{(q_1-q_2)^2}}{q_1}
\end{equation}
on $S^1\setminus Z_{q_1}$. Then by~\eqref{eq:q1} we have
$$
\beta_1=-\frac{2}{q_1^3},
$$
and taking a time derivative we obtain
$$
\dot\beta_1=\frac{3\cdot 2}{q_1^4}\dot q_1=-\frac{3\beta_1 \dot q_1}{q_1}.
$$
We multiply both sides with $q_1^2$ to get
$$
\dot\beta_1q_1^2=-3\beta_1\dot q_1q_1\,.
$$
We bring $-2\beta_1\dot q_1q_1$ to the other side to continue
$$
\dot\beta_1q_1^2+2\beta_1\dot q_1q_1=-\beta_1\dot q_1q_1\,.
$$
We substitute the original definition \eqref{eq:beta-inst} 
of $\beta_1$ in the right hand side to get 
$$
\dot\beta_1q_1^2+2\beta_1\dot q_1q_1=-\dot q_1\ddot q_1+\frac{\dot q_2}{(q_1-q_2)^2}
-\frac{\dot q_2-\dot q_1}{(q_2-q_1)^2}.
$$
Integrating both sides from $t$ to $1$ we get
\begin{equation}\label{eq:aux1}
\beta_1q_1^2=-\frac{\dot{q}_1^2}{2}-\int_t^1 \frac{\dot{q}_2(s)}{\big(q_1(s)-q_2(s)\big)^2}ds
+\frac{q_2-q_1}{(q_2-q_1)^2}+C
\end{equation}
for some constant $C\in \R$. We bring $-\frac{q_1}{(q_2-q_1)^2}$ to the left hand side and 
rewrite \eqref{eq:beta-inst} in the form
$$
\beta_1q_1^2+\frac{q_1}{(q_2-q_1)^2}=\ddot q_1 q_1\,.
$$
Therefore \eqref{eq:aux1} transforms to 
\begin{equation}\label{eq:aux2}
\ddot q_1 q_1=-\frac{\dot{q}_1^2}{2}-\int_t^1 \frac{\dot{q}_2(s)}{\big(q_1(s)-q_2(s)\big)^2}ds
+\frac{q_2}{(q_2-q_1)^2}+C.
\end{equation}
Observe that modulo the exact value of $C$ this is exactly equation~\eqref{qinst2}, 
which is equivalent to the first equation of \eqref{eq:critinst}
An analogous discussion applied to $q_2$
will lead to the second equation of \eqref{eq:critinst}. Therefore, we are left with 
computing the constant $C$. 

To compute $C$ we use \eqref{eq:q1} to get rid of $\ddot q$ on the left hand side of 
\eqref{eq:aux2} and obtain
$$
-\frac{2}{q_1}+\frac{q_1}{(q_1-q_2)^2}
=-\frac{\dot{q}_1^2}{2}-\int_t^1
\frac{\dot{q}_2(s)}{\big(q_1(s)-q_2(s)\big)^2}ds 
+\frac{q_2}{(q_2-q_1)^2}+C\,.
$$
We solve the last equation for $C$,
$$
C=-\frac{2}{q_1}+\frac{1}{q_1-q_2}
+\frac{\dot{q}_1^2}{2}+\int_t^1 \frac{\dot{q}_2}{\big(q_1-q_2\big)^2}ds\,,
$$
and integrate both sides from $0$ to $1$. Noting that 
$$
\int_0^1dt\int_t^1 \frac{\dot{q}_2(s)}{\big(q_1(s)-q_2(s)\big)^2}ds=
\int_0^1ds\int_0^s dt\frac{\dot{q}_2(s)}{\big(q_1(s)-q_2(s)\big)^2}=
\int_0^1\frac{s\dot{q}_2(s)ds}{\big(q_1(s)-q_2(s)\big)^2},
$$
this gives us
$$
C=-\int_0^1\frac{2ds}{q_1(s)}+\int_0^1\frac{ds}{q_1(s)-q_2(s)}
+\frac{||\dot{q}_1||^2}{2}+\int_0^1\frac{s\dot{q}_2ds}{\big(q_1-q_2\big)^2},
$$
which matches the constant in equation~\eqref{qinst}.
This concludes the proof of Lemma~\ref{lem:eqzinst}, and therefore of
Theorem~\ref{thm:inst-interaction}. 
\hfill$\square$

%%%%%%%%%%%%%%%%%%%%%%%%%%%%%%%%%%%%%%%%%%%%%%%%%%%%%%%%%%%%%%%%%%%%%%%%%%
\section{Interpolation}\label{sec:interpol}
%%%%%%%%%%%%%%%%%%%%%%%%%%%%%%%%%%%%%%%%%%%%%%%%%%%%%%%%%%%%%%%%%%%%%%%%%%

We now interpolate linearly between the instantaneous and mean
interactions. That is, for $r\in[0,1]$ we consider the system of coupled ODEs
\begin{equation}\label{eq:interpol}
\left\{\;
\begin{aligned}
  \ddot{q}_1(t) &= -\frac{2}{q_1(t)^2} + \frac{r}{(q_1(t)-q_2(t))^2} 
  + \frac{1-r}{(\overline{q}_1-\overline{q}_2)^2}, \cr
  \ddot{q}_2(t) &= -\frac{2}{q_2(t)^2} - \frac{r}{(q_1(t)-q_2(t))^2}
  - \frac{1-r}{(\overline{q}_1-\overline{q}_2)^2}.
\end{aligned}
\right.
\end{equation}
For $r=0$ this agrees with the system~\eqref{eq:mean-interaction} for
mean interaction, and for $r=1$ with the
system~\eqref{eq:inst-interaction} for instantaneous interaction.  
Solutions of~\eqref{eq:interpol} are critical points of the functional
$r\mathcal{S}_{in}+(1-r)\mathcal{S}_{av}$, which under Levi-Civita transformation corresponds to the functional
\begin{equation}\label{eq:Br}
\mathcal{B}_r:= r\BB_{in} + (1-r)\BB_{av}:
\mathcal{H}_{av}^1\stackrel{\eqref{eq:inavincl}}{\supset}\mathcal{H}_{in}^1\rightarrow \R
\end{equation}
that computes as
\begin{eqnarray*}
\mathcal{B}_r(z_1,z_2)
&=& r\BB_{in}(z_1,z_2) + (1-r)\BB_{av}(z_1,z_2) \\ \nonumber
&=&2\sum_{i=1}^2\Bigg(||z_i||^2 \cdot ||z_i'||^2 
+\frac{1}{||z_i||^2}\Bigg)
-\frac{r||z_1||^2 \cdot ||z_2||^2}
{||z_1^2||^2\cdot ||z_2||^2-||z_2^2||^2 \cdot ||z_1||^2}\\ \nonumber
& &-\frac{1-r}{||z_1||^2}\int_0^1\frac{z_1(\tau)^2}{z_1^2(\tau)-z_2^2
(\tau_{z_2}(t_{z_1}(\tau)))}d\tau.
\end{eqnarray*}
We define a {\em generalized solution $(q_1,q_2)$ of~\eqref{eq:interpol}} as in
Section~\ref{ss:gen-solutions-inst}, only that instead
of~\eqref{eq:inst-interaction} it now solves~\eqref{eq:interpol}.
Then we have the following generalization of Theorems~\ref{thm:mean-interaction}
and~\ref{thm:inst-interaction}.  

\begin{thm}[Generalized solutions for interpolated interaction]\label{thm:interpol} 
Under the Levi-Civita transformations~\eqref{eq:LC-zi} with time
changes~\eqref{eq:t-tau-zi}, critical points $(z_1,z_2)$ of the action
functional $\mathcal{B}_r$, $r\in[0,1]$ are in 4-to-1 correspondence
with generalized solutions $(q_1,q_2)$ of~\eqref{eq:interpol}.
\end{thm}

{\bf Proof: }
The proof is very similar to the proofs of Theorems~\ref{thm:mean-interaction}
and~\ref{thm:inst-interaction}. Critical points of $\BB_r$ are solutions of the problem 
\begin{eqnarray}\label{eq:z-interpol}
z_1''(\tau) \nonumber
&=&\frac{||z_1'||^2z_1(\tau)}{||z_1||^2}-\frac{z_1(\tau)}{||z_1||^6}+\frac{(1-r)z_1(\tau)}{2||z_1||^6}\int_0^1\frac{z_1(\sigma)^2}{z_1^2(\sigma)-z_2^2(\tau_{z_2}(t_{z_1}(\sigma)))}d\sigma\\ \nonumber
& &+\frac{(1-r)z_1(\tau)z_2^2(\tau_{z_2}(t_{z_1}(\tau)))}{2||z_1||^4\big(z_1^2(\tau)-z_2^2(\tau_{z_2}(t_{z_1}(\tau)))\big)^2}\\ \nonumber
& &-\frac{(1-r)||z_2||^2z_1(\tau)}{||z_1||^6}
\int_\tau^1 \frac{z_1^2(\sigma) z_2'(\tau_{z_2}(t_{z_1}(\sigma)))}{z_2(\tau_{z_2}(t_{z_1}(\sigma)))\big(z_1^2(\sigma)-z_2^2(\tau_{z_2}(t_{z_1}(\sigma)))\big)^2}d\sigma\\ \nonumber
& &+\frac{(1-r)||z_2||^2 z_1(\tau)}{||z_1||^6}
\int_0^1 \frac{z_1^2(\sigma) z_2'(\tau_{z_2}(t_{z_1}(\sigma))) t_{z_1}(\sigma)}
{z_2(\tau_{z_2}(t_{z_1}(\sigma)))\big(z_1^2(\sigma)-z_2^2(\tau_{z_2}(t_{z_1}(\sigma)))\big)^2}d\sigma\\ \nonumber
& &-\frac{r||z_2||^4 ||z_1^2||^2 z_1(\tau)}{2||z_1||^2 \big(||z_1^2||^2 ||z_2||^2-||z_2^2||^2
||z_1||^2\big)^2}\\ \nonumber
& &+\frac{r||z_2||^4 z_1^3(\tau)}{\big(||z_1^2||^2 ||z_2||^2-||z_2^2||^2||z_1||^2\big)^2}\\ \nonumber
z_2''(\tau)&=&\frac{||z_2'||^2z_2(\tau)}{||z_2||^2}-\frac{z_2(\tau)}{||z_2||^6}-\frac{(1-r)z_2(\tau)}{2||z_2||^6}\int_0^1\frac{z_2(\sigma)^2}{z_2^2(\sigma)-z_1^2(\tau_{z_1}(t_{z_2}(\sigma)))}d\sigma\\ \nonumber
& &-\frac{(1-r)z_2(\tau)z_1^2(\tau_{z_1}(t_{z_2}(\tau)))}{2||z_2||^4\big(z_2^2(\tau)-z_1^2(\tau_{z_1}(t_{z_2}(\tau)))\big)^2}\\ \nonumber
& &+\frac{(1-r)||z_1||^2z_2(\tau)}{||z_2||^6}
\int_\tau^1 \frac{z_2^2(\sigma) z_1'(\tau_{z_1}(t_{z_2}(\sigma)))}{z_1(\tau_{z_1}(t_{z_2}(\sigma)))\big(z_2^2(\sigma)-z_1^2(\tau_{z_1}(t_{z_2}(\sigma)))\big)^2}d\sigma\\ \nonumber
& &-\frac{(1-r)||z_1||^2 z_2(\tau)}{||z_2||^6}
\int_0^1 \frac{z_2^2(\sigma) z_1'(\tau_{z_1}(t_{z_2}(\sigma))) t_{z_2}(\sigma)}{z_1(\tau_{z_1}(t_{z_2}(\sigma)))\big(z_2^2(\sigma)-z_1^2(\tau_{z_1}(t_{z_2}(\sigma)))\big)^2}d\sigma\\ \nonumber
& &+\frac{r||z_1||^4||z_2^2||^2 z_2(\tau)}{2||z_2||^2\big(||z_1^2||^2||z_2||^2-||z_2^2||^2||z_1||^2\big)^2}\\
& &-\frac{r||z_1||^4 z_2^3(\tau)}{\big(||z_1^2||^2||z_2||^2-||z_2^2||^2||z_1||^2\big)^2}
\end{eqnarray}
which interpolates between problems~\eqref{eq:z-av}
and~\eqref{eq:critinst}.

Suppose now that $(z_1,z_2)$ is a critical point of $\BB_r$ and define
$$q_i(t) := z_i^2(\tau_{z_i}(t)), \qquad i \in \{1,2\}.$$
From~\eqref{eq:z-interpol} we obtain for $q_1$ the equation
\begin{eqnarray}\label{qint}
2\ddot{q}_1(t) q_1(t)&=&||\dot{q}_1||^2-\int_0^1 \frac{4}{q_1(s)}ds-\dot{q}_1^2(t)+\int_0^1\frac{2(1-r)}{q_1(s)-
q_2(s)}ds\\ \nonumber
& &+\frac{2(1-r)q_2(t)}{\big(q_1(t)-q_2(t)\big)^2}-\int_t^1 \frac{2(1-r)\dot{q}_2(s)}{\big(q_1(s)-q_2(s)\big)^2}ds\\ \nonumber
& &+\int_0^1 \frac{2(1-r)s\dot{q}_2(s)}{\big(q_1(s)-q_2(s)\big)^2}ds
-\frac{2r\overline{q}_1}{(\overline{q}_1-\overline{q}_2)^2}
+\frac{4r q_1(t)}{(\overline{q}_1-\overline{q_2})^2}
\end{eqnarray}
which interpolates between (\ref{bov6}) and (\ref{qinst}). As before,
outside collisions we define 
$$\beta_1:=\frac{\ddot{q}_1-\frac{r}{(q_1-q_2)^2}-\frac{1-r}{(\overline{q}_1-\overline{q}_2)^2}}{q_1}.$$
Again $\beta_1$ solves the ODE
$$\dot{\beta}_1 q_1=-3\beta_1 \dot{q}_1$$
and it follows that
$$\beta_1=-\frac{\mu}{q_1^3},$$
where $\mu$ is locally constant on the complement of collisions on the
circle. Using the continuity of $z_1'$, it follows again that $\mu$ is
constant and from~\eqref{qint} we conclude that $\mu=2$. 
Therefore, outside of collisions $q_1$ solves the first equation
in~\eqref{eq:interpol}. Similarly, it follows that outside collisions
$q_2$ solves the second equation in~\eqref{eq:interpol}. 
The converse direction is proved similarly as in the previous cases.
\hfill$\square$

%%%%%%%%%%%%%%%%%%%%%%%%%%%%%%%%%%%%%%%%%%%%%%%%%%%%%%%%%%%%%%%%%%%%%%%%%%
\section{Existence of symmetric frozen planet orbits}
%%%%%%%%%%%%%%%%%%%%%%%%%%%%%%%%%%%%%%%%%%%%%%%%%%%%%%%%%%%%%%%%%%%%%%%%%%

A {\em simple frozen planet orbit of period $T>0$} is a map
$q=(q_1,q_2)\in H^1(\R/T\Z,\R^2)$ with the following properties:
\begin{enumerate}
\item $q_1(t)>q_2(t)\geq 0$ for all $t\in\R/T\Z$;
\item $q_2$ has a unique zero at $t=0$;
\item $(q_1,q_2):(0,T)\to\R^2$ is
  smooth and satisfies~\eqref{eq:inst-interaction}; 
\item the energies
  $$
    E_i(t) := \frac{\dot q_i(t)^2}{2}-\frac{2}{q_i(t)},\qquad t\in
    (0,T)
  $$
extend to continuous functions $E_i:\R/T\Z\to\R$. 
\end{enumerate}
Here simplicity corresponds to the second condition, and every frozen
planet orbit is a multiple cover of a simple one. 
Recall that the individual energies $E_i$ need not be constant, but
the total energy
$$
   E = E_1(t) + E_2(t) + \frac{1}{q_1(t)-q_2(t)}
$$
is constant and negative. A simple frozen planet orbit of period $T>0$
is called {\em symmetric} if it satisfies in addition
$$
   q(t) = q(T-t)\quad\text{for all }t\in\R/T\Z. 
$$
In this section we prove the following result which corresponds to
Theorem\,C in the Introduction.

\begin{thm}\label{thm:existence}
For every $E<0$ there exists a symmetric simple frozen planet orbit of energy $E$.
\end{thm}

{\bf Rescaling. }
Let $q=(q_1,q_2):\R\to(\R_{\geq0})^2$ be a generalized solution
of~\eqref{eq:inst-interaction} of period $T$ and energy $E$.
Direct computation shows that for each $c>0$,
$$
   q_c(t):=c^2q(c^{-3}t)
$$
is again a generalized solution
of~\eqref{eq:inst-interaction} of period $c^3T$ and energy $c^{-2}E$. 
Therefore, given a generalized solution $q$ of~\eqref{eq:inst-interaction}
of period $1$ and negative energy, rescaling yields similar
solutions with any prescribed period, or alternatively with any
prescribed negative energy. 
As a consequence, we will from now on restrict our discussion to
generalized solutions of period $1$.

%%%
\subsection{Symmetries}\label{ss:symmetries}
%%%

In this subsection we describe the symmetries of the variational
problems in Sections~\ref{sec:Kepler}, \ref{sec:mean-interaction}
and~\ref{sec:inst-interaction}. 

{\bf The Kepler problem. }
The functional $\QQ$ in~\eqref{eq:Q} is obviously invariant under the
following transformations: 
\begin{itemize}
\item time shift $T_sz(\tau):=z(s+\tau)$, $s\in S^1$;
\item time reversal $Rz(\tau):=z(-\tau)$;
\item sign reversal $Sz(\tau):= -z(\tau)$,
\end{itemize}
where $z\in H^1(S^1,\R)\setminus\{0\}$.
\begin{lemma}\label{lem:time-rev}
The homeomorphism $t_z:S^1\to S^1$ defined in~\eqref{eq:tz} 
and its inverse $\tau_z$ transform under time shift and
reversal as
\begin{equation}\label{eq:shift}
t_{T_sz}(\tau)=t_z(s+\tau)-t_z(s),\quad 
\tau_{T_sz}(t)=\tau_z(t+t_z(s))-s,
\end{equation}
\begin{equation}\label{eq:rev}
   t_{Rz}(\tau) = -t_z(-\tau),\qquad \tau_{Rz}(t) = -\tau_z(-t).
\end{equation}
\end{lemma}

{\bf Proof: }
For the first equation of \eqref{eq:shift} we compute 
\begin{eqnarray*}
t_{T_sz}(\tau)&=& 
\frac{1}{\|T_sz\|^2}\int_0^\tau z(s+\sigma)^2d\sigma\\
&=&
\frac{1}{\|z\|^2}\left(\int_0^{s+\tau} z(\sigma)^2d\sigma-
\int_0^s z(\sigma)^2d\sigma\right)\\
&=&t_z(s+\tau)-t_z(s).
\end{eqnarray*}

For the second equation of~\eqref{eq:shift} note that 
$t_{T_sz}$ writes out as a composition of three homeomorphisms 
of the circle: 
$$t_{T_sz}=T_{-t_z(s)}\circ t_z\circ T_s.$$ 
The inverse $\tau_{T_sz}$ of $t_{T_sz}$ is the composition of the 
three inverse homeomorphisms in the opposite order: 
$$
\tau_{T_sz}=T_{-s}\circ\tau_z\circ T_{t_z(s)}.
$$
The first equation of~\eqref{eq:rev} follows from 
\begin{eqnarray*}
  t_{Rz}(\tau)
  &=& \frac{\int_0^\tau (Rz)^2(\sigma)d\sigma}{\|Rz\|^2}
  = \frac{\int_0^\tau z^2(-\sigma)d\sigma}{\|z\|^2}
  = -\frac{\int_0^{-\tau} z^2(\sigma)d\sigma}{\|z\|^2}
  = -t_{z}(-\tau),
\end{eqnarray*}
and the second equation of~\eqref{eq:rev} follows from the first one by splitting $t_{Rz}$ into three homeomorphisms in analogy
with the above.
%checking that the two maps are inverse to one another:
%\begin{align*}
%  t_{Rz}(\tau_{Rz}(t))
%  = t_{Rz}(-\tau_{z}(-t))
%  = -t_{z}(\tau_{z}(-t))
%  = -(-t) =t.
%\end{align*}
%$$\tau_{Iz_i}(t_{Iz_i}(\tau))=\tau_{Iz_i}(-t_{z_i}(-\tau))=-\tau_{z_i}(t_{z_i}(-\tau))=-(-\tau)=\tau$$
\hfill $\square$

We discuss what these symmetries correspond to 
under the Levi-Civita transformation
$q(t)=z(\tau_z(t))^2$.
Let $q_s(t)$ denote 
the Levi-Civita transformation of $T_sz$. Then
$$
q_s(t)=(T_sz(\tau_{T_sz}(t)))^2=
(T_sz(\tau_z(t_z(s)+t)-s))^2=
z(\tau_z(t+t_z(s)))^2=q(t_z(s)+t).
$$
We see that the time shift $T_s$ applied to $z$ 
corresponds to the time 
shift $T_{t_z(s)}$ applied to the Levi-Civita transformation 
$q$ of $z$. Analogous but simpler arguments show that 
time reversal corresponds to time reversal $Rq(t)=q(-t)$, 
and sign reversal corresponds to the identity $Sq(t)=q(t)$. 
Thus, the first two symmetries of $z$ correspond to actual symmetries of the Kepler problem, while the sign change $Sz=-z$ just expresses the fact the the Levi-Civita transformation defines a 2-to-1 correspondence.
\smallskip

{\bf Mean and instantaneous interactions. }
The problems with interaction have the following symmetries.

\begin{lemma}\label{lem:inv}
The functionals $\B_{av}$ and $\B_{in}$ in~\eqref{eq:Bav}
and~\eqref{eq:Bin} are invariant under the 
following transformations of $z=(z_1,z_2)$:
\begin{itemize}
\item joint time shift $T_sz(\tau):=z(s+\tau)$, $s\in S^1$;
\item joint time reversal $Rz(\tau):=z(-\tau)$;
\item separate sign reversals $S_1(z_1,z_2):=(-z_1,z_2)$ and $S_2(z_1,z_2):=(z_1,-z_2)$.
\end{itemize}
\end{lemma}

\textbf{Proof: }
Invariance of $\B_{av}$ and $\B_{in}$ under joint time shift $T_s$ and
separate sign changes $S_i$ is obvious, and so is invariance of $\B_{av}$ under
joint time reversal $R$ (in fact, $\B_{av}$ is even invariant under
separate time reversals). For invariance of $\B_{in}$ under $R$ we
write the instantaneous interaction from~\eqref{eq:I} in the form
$$
\mathcal{I}(z_1,z_2)=\int_0^1\frac{dt}{z_2^2(\tau_{z_2}(t))-z_1^2(\tau_{z_1}(t))}.
$$
Using Lemma~\ref{lem:time-rev} we compute
$$
  (Rz_i)(\tau_{Rz_i}(t)) = z_i(-\tau_{Rz_i}(t)) = z_i(\tau_{z_i}(-t))
$$
for $i=1,2$, and therefore
%by substitution $t\mapsto -t$ in the integral
\begin{eqnarray*}
  \mathcal{I}(R(z_1,z_2))
  &=& \int_0^1\frac{dt}{z_2^2(\tau_{z_2}(-t))-z_1^2(\tau_{z_1}(-t))}\\
  &=& \int_0^1\frac{dt}{z_2^2(\tau_{z_2}(t))-z_1^2(\tau_{z_1}(t))}
  = \mathcal{I}(z_1,z_2).
\end{eqnarray*}
\hfill $\square$

%%%
\subsection{Twisted loops}\label{ss:twisted}
%%%

Theorems~\ref{thm:BOV}, ~\ref{thm:mean-interaction} and~\ref{thm:inst-interaction}
establish a correspondence between critical points of the functionals
$\QQ$, $\BB_{av}$ and $\BB_{in}$ and generalized solutions with an
even number of zeroes. In this subsection we explain how to deal with
generalized solutions with an {\em odd} number of zeroes. 

{\bf The Kepler problem. }
Consider $q\in H^1(S^1,\R_{\geq 0})$ satisfying all conditions on a
generalized solution of~\eqref{eq:Kepler} except that it has an odd
number of zeroes. By slight abuse of notation we will still refer to
such $q$ as a ``generalized solution''.
Recall that $S^1=\R/\Z$. Then we can view $q$ as a map
$\wt q\in H^1(\R/2\Z,\R_{\geq 0})$ with an even number of zeroes and
associate to it as in Section~\ref{ss:q-to-z-Kepler} a
homeomorphism $\tau_q:\R/2\Z\to \R/2\Z$ and a map $z\in H^1(\R/2\Z,\R)$
such that
$$
  z(\tau_q(t))^2=\wt q(t),\qquad t\in\R/2\Z.
$$
By construction these maps satisfy
$$
   \wt q(t+1) = \wt q(t),\qquad 
   \tau_q(t+1) = \tau_q(t),\qquad 
   z(t+1) = -z(t),\qquad t\in\R/2\Z. 
$$
This leads us to introduce for each $k\in\N_0$ the Hilbert space of
{\em twisted loops}
$$
  H_{twist}^k(S^1,\R) := \{z\in H^k(\R/2\Z,\R)\mid z(\tau+1)=-z(\tau)
  \text{ for all }\tau\}
$$
with the inner product
$$
   \la z,v\ra := \int_0^1z(\tau)v(\tau)d\tau = \frac12\int_0^2z(\tau)v(\tau)d\tau.
$$
Note that $H_{twist}^k(S^1,\R)$ is the fixed point set of the linear involution
$$
   \sigma=S\circ T_1:H^k(\R/2\Z,\R) \to H^k(\R/2\Z,\R),\qquad \sigma z(\tau)=-z(\tau+1),
$$
where $T_1$ and $S$ are the time shift by $1$ and sign reversal from
Section~\ref{ss:symmetries}. Thus $H_{twist}^k(S^1,\R)$ is a closed linear
subspace of $H^k(\R/2\Z,\R)$, with inner product the one induced from
$H^k(\R/2\Z,\R)$ divided by $2$. 
We define
\begin{equation*}%\label{eq:Q}
  \QQ: H_{twist}^1(S^1,\R)\setminus\{0\}\to\R,\qquad 
  \QQ(z) := 2\|z\|^2\|z'\|^2 + \frac{N}{\|z\|^2}
\end{equation*}
by the same formula as in~\eqref{eq:Q}. As in
Section~\ref{ss:z-to-q-Kepler} it follows that $z\in 
H_{twist}^1(S^1,\R)\setminus\{0\}$
is a critical point of $\QQ$ if and only if
\begin{equation}\label{eq:DQ-twisted}
  \langle z',v'\rangle + \langle az,v\rangle = 0\qquad \text{for all }v\in H_{twist}^1(S^1,\R), 
\end{equation}
with the constant 
\begin{equation*}%\label{eq:z-Kepler}
  a = \frac{\|z'\|^2}{\|z\|^2} - \frac{N}{2\|z\|^6}.
\end{equation*}
To proceed we need the following lemma. 

\begin{lemma}\label{lem:twisted-var}
Suppose that $f,g\in H_{twist}^0(S^1,\R)$ satisfy
$$
  \la f,v'\ra + \la g,v\ra = 0 \qquad\text{for all }
  v\in H_{twist}^1(S^1,\R).
$$
Then $f\in H_{twist}^1(S^1,\R)$ with weak derivative $f'=g$. 
\end{lemma}

{\bf Proof: }
By definition of $H_{twist}^0(S^1,\R)$ we have $f,g\in H^0(\R/2\Z,\R)$
with \\ $f(\tau+1)=-f(\tau)$ and $g(\tau+1)=-g(\tau)$, so the hypothesis of
the lemma reads
\begin{equation}\label{eq:weak-der}
  \int_0^2 f(\tau)v'(\tau)d\tau + \int_0^2g(\tau)v(\tau)d\tau = 0 
\end{equation}
for all $v\in H^1(\R/2\Z,\R)$ with $v(\tau+1)=-v(\tau)$. 
By definition of the weak derivative, we need to show that~\eqref{eq:weak-der}
holds for all $v\in H^1(\R/2\Z,\R)$ (not necessarily satisfying $v(\tau+1)=-v(\tau)$).
Multiplying by bump functions and using linearity, it suffices to show
this for $v\in H^1(\R/2\Z,\R)$ with support in an interval
$I\subset\R/2\Z$ of length less than $1$. Given such $v$ we define
$\wt v\in H^1(\R/2\Z,\R)$ by
$$
  \wt v(\tau) := \begin{cases}
    v(\tau) & \tau \in I, \\  
    -v(\tau-1) & \tau-1 \in I, \\  
    0 & \text{otherwise}.
  \end{cases}
$$
Using $f(\tau+1)=-f(\tau)$ and $\wt v'(\tau+1)=-\wt v'(\tau)$ we compute
\begin{align*}
  \int_0^2 f(\tau)\wt v'(\tau)d\tau
  &= \int_I f(\tau)\wt v'(\tau)d\tau + \int_I f(\tau+1)\wt v'(\tau+1)d\tau\cr
  &= \int_I f(\tau)\wt v'(\tau)d\tau + \int_I f(\tau)\wt v'(\tau)d\tau\cr
  &= 2\int_0^2 f(\tau)v'(\tau)d\tau,
\end{align*}
and similarly
$$
  \int_0^2 g(\tau)\wt v(\tau)d\tau
  = 2\int_0^2 g(\tau)v(\tau)d\tau.
$$
Since~\eqref{eq:weak-der} holds with $\wt v$ in place of $v$, we conclude
\begin{align*}
  0
  &= \int_0^2 f(\tau)\wt v'(\tau)d\tau + \int_0^2g(\tau)\wt v(\tau)d\tau
  = 2\bigg(\int_0^2 f(\tau)v'(\tau)d\tau + \int_0^2g(\tau)v(\tau)d\tau\bigg)
\end{align*}
and the lemma is proved. 
\hfill $\square$

Let us now return to the critical point 
$z\in H_{twist}^1(S^1,\R)\setminus\{0\}$ of $\QQ$ satisfying~\eqref{eq:DQ-twisted} above.
Applying Lemma~\ref{lem:twisted-var} with $f=z'$ and $g(z)=az$, we
conclude that $z\in H_{twist}^2(S^1,\R)$ and its second weak derivative satisfies
$$
  z'' = az. 
$$
This is the same ODE~\eqref{eq:z-Kepler} as for critical points in
the untwisted case. Now the arguments in Section~\ref{sec:Kepler}
carry over without further changes to show the following twisted version of
Theorem~\ref{thm:BOV}. 

\begin{thm}\label{thm:BOV-twisted}
Under the Levi-Civita transformation~\eqref{eq:LC} with time
change~\eqref{eq:t-tau}, critical points $z$ of the action
functional $\QQ: H_{twist}^1(S^1,\R)\setminus\{0\}\to\R$ on twisted
loops are in 2-to-1 correspondence
with generalized solutions $q:S^1\to\R_{\geq 0}$ of~\eqref{eq:Kepler}
having an odd number of zeroes. 
\hfill $\square$
\end{thm}

{\bf Mean and instantaneous interactions. }
Consider now the helium atom. Recall that in a frozen planet
configuration the inner electron $q_2$ should undergo repeated
collisions with the nucleus while the outer electron $q_1$ experiences no
collisions. If in period $1$ the inner electron undergoes an odd number of
collisions (for example a single one), the Levi-Civita transformed
maps $z_1,z_2$ will be $2$-periodic and satisfy
$$
  z_1(\tau+1) = z_1(\tau),\qquad z_2(\tau+1) = -z_2(\tau).
$$
This leads us to introduce for each $k\in\N_0$ the Hilbert space of
{\em twisted loops} 
\begin{align*}
  H_{twist}^k(S^1,\R^2) 
  := \bigl\{&z=(z_1,z_2)\in H^k(\R/2\Z,\R)\;\bigl| \\
  &z_1(\tau+1) = z_1(\tau),\ z_2(\tau+1)=-z_2(\tau)\text{ for all }\tau\bigr\}
\end{align*}
with the inner product
$$
  \la z,v\ra := \sum_{i=1}^2\int_0^1z(\tau)v(\tau)d\tau
  = \frac12\sum_{i=1}^2\int_0^2z(\tau)v(\tau)d\tau.
$$
Note that $H_{twist}^k(S^1,\R^2)$ is the fixed point set of the linear involution
$$
   \sigma=S_2\circ T_1:H^k(\R/2\Z,\R^2) \to H^k(\R/2\Z,\R^2),\quad \sigma
   z(\tau)=\bigl(z_1(\tau+1),-z_2(\tau+1)\bigr), 
$$
where $T_1$ and $S_2$ are the joint time shift by $1$ and sign reversal in
the second component from
Section~\ref{ss:symmetries}. Thus $H_{twist}^k(S^1,\R^2)$ is a closed linear
subspace of $H^k(\R/2\Z,\R^2)$, with inner product the one induced from
$H^k(\R/2\Z,\R^2)$ divided by $2$. 
We define open subsets
\begin{equation*}%\label{eq:Hav}
   \wt{\mathcal{H}}_{av}^1 := \Bigg\{z=(z_1,z_2) \in
   H_{twist}^1(S^1,\mathbb{R}^2)\;\Bigl|\; ||z_1||>0,\,\,||z_2||>0,\,\, 
   \frac{||z_1^2||^2}{||z_1||^2}>\frac{||z_2^2||^2}{||z_2||^2}\Bigg\}
\end{equation*}
and
\begin{equation}\label{eq:Hin-twisted}
\begin{aligned}
  \wt{\mathcal{H}}_{in}^1 := \Bigl\{
  &z=(z_1,z_2) \in H_{twist}^1(S^1,\mathbb{R}^2)\;\Bigl|\; 
  \|z_1\|>0,\,\,\|z_2\|>0,\,\, \cr
  &\ \ z_1^2(\tau)-z_2^2(\tau_{z_2}(t_{z_1}(\tau)))>0 \,\,\text{for
    all}\,\,\tau\in S^1\Bigr\}
\end{aligned}
\end{equation}
as in Sections~\ref{sec:mean-interaction} and~\ref{sec:inst-interaction},
and we define the mean and instantaneous interaction functionals and
their interpolation
\begin{equation*}
  \BB_{av}:\wt{\mathcal{H}}_{av}^1\to\R,\qquad 
  \BB_{in},\BB_r:\wt{\mathcal{H}}_{in}^1\to\R
\end{equation*}
by the same formulas as in~\eqref{eq:Bav}, ~\eqref{eq:Bin} and~\eqref{eq:Br}.
Note that in formulas~\eqref{eq:DBav} and~\eqref{eq:DBin} for the
derivatives of $\BB_{av}$ and $\BB_{in}$ the leading order terms $\la
z_i',v_i'\ra$ are the same as in the Kepler problem. Arguing as above
using Lemma~\ref{lem:twisted-var}, we thus conclude that critical
points of $\BB_{av}$ and $\BB_{r}$ on twisted loops
still satisfy the same equations~\eqref{eq:z-av}
and~\eqref{eq:z-interpol} as in the untwisted case. Therefore, in
analogy with Theorem~\ref{thm:BOV-twisted}, we obtain the following
twisted version of Theorems~\ref{thm:mean-interaction}
and~\ref{thm:interpol}. 

\begin{thm}\label{thm:interpol-twisted}
Under the Levi-Civita transformations~\eqref{eq:LC-i} with time
changes~\eqref{eq:t-tau-i}, critical points $(z_1,z_2)$ of the action
functional $\mathcal{B}_r:\wt{\mathcal{H}}_{av}^1\to\R$ on twisted
loops are in 4-to-1 correspondence
with generalized solutions $(q_1,q_2)$ of~\eqref{eq:mean-interaction}
with $q_1$ having an even and $q_2$ an odd number of zeroes.
Similarly, for each $r\in[0,1]$ critical points $(z_1,z_2)$ of the action
functional $\mathcal{B}_{r}:\wt{\mathcal{H}}_{in}^1\to\R$ on twisted
loops are in 4-to-1 correspondence
with generalized solutions $(q_1,q_2)$ of~\eqref{eq:interpol}
with $q_1$ having an even and $q_2$ an odd number of zeroes.
\hfill $\square$
\end{thm}

%%%
\subsection{Symmetric loops}\label{ss:symmetric}
%%%

In order to study symmetric frozen planet orbits, we introduce
for each $k\in\N_0$ the Hilbert space of {\em symmetric loops} 
\begin{align*}
  H_{sym}^k(S^1,\R^2) 
  := \bigl\{&z=(z_1,z_2)\in H_{twist}^k(S^1,\R^2)\;\bigl| \\
  &z(\tau) = z(1-\tau)\text{ for all }\tau\bigr\}
\end{align*}
with the inner product induced from $H_{sym}^k(S^1,\R^2)$. 
Thus $H_{sym}^k(S^1,\R^2)$ is the fixed point set of the linear involution
$$
   \rho=R\circ T_1:H_{twist}^k(S^1,\R^2) \to H_{twist}^k(S^1,\R^2),\quad \rho
   z(\tau)=z(1-\tau), 
$$
where $T_1$ and $R$ are the joint time shift by $1$ and joint time reversal from Section~\ref{ss:symmetries}. 

Elements $(z_1,z_2)\in H_{sym}^k(S^1,\R^2)$ are $2$-periodic loops
satisfying for all $\tau\in\R/2\Z$ the conditions 
$$
   z_1(\tau) = z_1(\tau+1) = z_1(1-\tau),\qquad
   z_2(\tau) = -z_2(\tau+1) = z_2(1-\tau). 
$$
Taking derivatives they imply
$$
   z_1'(\tau) = z_1'(\tau+1) = -z_1'(1-\tau),\qquad
   z_2'(\tau) = -z_2'(\tau+1) = -z_2'(1-\tau). 
$$
In particular, at $\tau=0$ and $\tau=1/2$ they imply
$$
   z_1'(0) = z_2(0) = 0,\qquad z_1'(1/2) = z_2'(1/2) = 0.
$$
Thus $z_1$ is $1$-periodic with critical points at $\tau=0$ and
$\tau=1/2$, while $z_2$ is $2$-periodic with zeroes at $\tau=0$ and
$\tau=1$ and with critical points at $\tau=1/2$ and $\tau=3/2$.  

The images $q_i(t) = z_i(\tau_{z_i}(t))^2$ under the Levi-Civita transformation
are $1$-periodic and {\em symmetric}, i.e.
$$
  q_i(t)=q_i(1-t),\qquad t\in S^1.
$$
In particular, if $(q_1,q_2)$ satisfies in addition the
ODE~\eqref{eq:inst-interaction},
%and $q_1$ has no zero and $q_2$ has a unique zero,
then $(q_1,q_2)$ is a symmetric frozen planet orbit.

Again, we define open subsets
\begin{equation*}%\label{eq:Hav}
   \wh{\mathcal{H}}_{av}^k := \Bigg\{z=(z_1,z_2) \in
   H_{sym}^k(S^1,\mathbb{R}^2)\;\Bigl|\; ||z_1||>0,\,\,||z_2||>0,\,\, 
   \frac{||z_1^2||^2}{||z_1||^2}>\frac{||z_2^2||^2}{||z_2||^2}\Bigg\}
\end{equation*}
and
\begin{equation}\label{eq:Hin-symm}
\begin{aligned}
  \wh{\mathcal{H}}_{in}^k := \Bigl\{
  &z=(z_1,z_2) \in H_{sym}^k(S^1,\mathbb{R}^2)\;\Bigl|\; 
  \|z_1\|>0,\,\,\|z_2\|>0,\,\, \cr
  &\ \ z_1^2(\tau)-z_2^2(\tau_{z_2}(t_{z_1}(\tau)))>0 \,\,\text{for
    all}\,\,\tau\in S^1\Bigr\},
\end{aligned}
\end{equation}
and we define the mean and instantaneous interaction functionals and
their interpolation
\begin{equation*}
  \BB_{av}:\wh{\mathcal{H}}_{av}^1\to\R,\qquad 
  \BB_{in}:\wh{\mathcal{H}}_{in}^1\to\R,\qquad
  \BB_r:\wh{\mathcal{H}}_{av}^1\stackrel{\eqref{eq:inavincl}}
  {\supset}\wh{\mathcal{H}}_{in}^1\to\R
\end{equation*}
by the same formulas as in~\eqref{eq:Bav}, ~\eqref{eq:Bin} and~\eqref{eq:Br}.
Arguing as in the previous subsection, using a variant of
Lemma~\ref{lem:twisted-var}, we conclude that critical
points of $\BB_{av}$ and $\BB_{r}$ on symmetric loops
still satisfy the same equations~\eqref{eq:z-av}
and~\eqref{eq:z-interpol} and we obtain the following
symmetric version of Theorems~\ref{thm:mean-interaction}
and~\ref{thm:interpol}. 

\begin{thm}\label{thm:interpol-symm}
Under the Levi-Civita transformations~\eqref{eq:LC-i} with time
changes~\eqref{eq:t-tau-i}, critical points $(z_1,z_2)$ of the action
functional $\mathcal{B}_{av}:\wh{\mathcal{H}}_{av}^1\to\R$ on symmetric
loops are in 4-to-1 correspondence with symmetric
generalized solutions $(q_1,q_2)$ of~\eqref{eq:mean-interaction}
with $q_1$ having an even and $q_2$ an odd number of zeroes.
Similarly, for each $r\in[0,1]$ critical points $(z_1,z_2)$ of the action
functional $\mathcal{B}_{r}:\wh{\mathcal{H}}_{in}^1\to\R$ on symmetric
loops are in 4-to-1 correspondence with symmetric
generalized solutions $(q_1,q_2)$ of~\eqref{eq:interpol}
with $q_1$ having an even and $q_2$ an odd number of zeroes
(i.e.~with symmetric frozen planet orbits).
\hfill $\square$
\end{thm}
 
%%%
\subsection{Proof of the Existence Theorem~\ref{thm:existence}}\label{ss:existence-proof}
%%%

Now we are ready to prove Theorem~\ref{thm:existence}. By the
rescaling discussion following the theorem, it suffices to consider the case of
period $1$ without prescribing the energy. 

Using the notation of the previous subsection we define
\begin{equation}\label{def:X}
   X := \{z=(z_1,z_2)\in\wh\HH^2_{in}\mid
   z_i(\tau)>0\text{ for all }\tau\in(0,1)\text{ and }i=1,2\}. 
\end{equation}
This is an open subset of the Hilbert space 
$H_{sym}^2(S^1,\R^2)$ and
thus a Hilbert manifold. Note that the twisting conditions 
$z_1(\tau+1)=z_1(\tau)$ and $z_2(\tau+1)=-z_2(\tau)$ imply that
$z_1(\tau)>0$ for all $\tau\in\R$ and $z_2(\tau)<0$ for $\tau\in(1,2)$. This
ensures that $z_2$ is simple of minimal period $2$ and it
removes the symmetries $z_i\mapsto\pm z_i$.
%, so we can hope for the
%existence of a unique critical point of $\BB_{in}$ in $X$.
We consider the Hilbert space 
$$
   Y := H_{sym}^0(S^1,\R^2)
$$
and the $L^2$-gradient (see beginning of Section~\ref{sec:diff}) of the interpolation functional $\BB_r=(1-r)\BB_{in}+r\BB_{av}$,
$$
   \nabla\BB_r=(1-r)\nabla\BB_{in}+r\nabla\BB_{av}:X\to Y,\qquad r\in[0,1]. 
$$
According to Theorem~\ref{thm:Fredholm}, for each $r\in[0,1]$ this is
a $C^1$-Fredholm map of index zero. Thus
$$
   F:[0,1]\times X\to Y,\qquad (r,z)\mapsto \nabla\BB_r(z)
$$
is a $C^1$-Fredholm map of index $1$. 

According to Theorem~\ref{thm:interpol-symm}, for each $r\in[0,1]$
zeroes $z\in X$ of $\nabla\BB_r$ correspond under the Levi-Civita
transformation to symmetric generalized solutions
$q=(q_1,q_2):S^1\to\R^2$ of~\eqref{eq:interpol}. The condition $z_2(\tau)>0$
for all $\tau\in(0,1)$ in the definition of $X$ implies that $q_2$ has
a unique zero at $t=0$. 
By the main result in~\cite{frauenfelder2} there exists a constant
$\kappa$ such that
$$
   \max_{t\in S^1}\bigg\{q_1(t),\frac{1}{q_1(t)-q_2(t)}\bigg\}\leq\kappa
$$
for all such solutions $q$ of~\eqref{eq:interpol} and all $r\in[0,1]$. Thus 
on $(r,z)\in F^{-1}(0)$ the function $z_1^2(t)$ is uniformly bounded
from above and the difference
$z_1^2(\tau_{z_1}(t))-z_2^2(\tau_{z_2}(t))$ is uniformly bounded away
from zero. In view of the ODE~\eqref{eq:z-interpol}, this implies that 
the zero set $F^{-1}(0)\subset[0,1]\times X$ is compact. Hence
$F:[0,1]\times X\to Y$ is a homotopy as in Theorem~\ref{thm:zero-Fredholm}
between $f_0=\nabla\BB_{in}$ and $f_1=\nabla\BB_{av}$, and it follows that
$f_0,f_1$ have well-defined mod $2$ Euler numbers satisfying
$$
   \chi(\nabla\BB_{in}) = \chi(\nabla\BB_{av}).
$$
%Now by the main result in~\cite{frauenfelder}, for the mean
%interaction there is a unique simple symmetric generalized solution
%$q$ of~\eqref{eq:mean-interaction}. By Theorem~\ref{thm:mean-interaction-symm}
%and the discussion following the definition~\eqref{def:X} of the space
%$X$, under the Levi-Civita transformation this corresponds to a unique
%zero $z\in X$ of $\nabla\BB_{av}$. According to Proposition~\ref{prop:triv-kernel},
%the differential $\nabla\BB_{av}(z):\wh H^2(S^1,\R^2)\to\wh H^0(S^1,\R^2)$ is 
%injective, hence surjective (because it is Fredholm of index $0$), so
%the (Transversality) axiom in Theorem~\ref{thm:zero-Fredholm} yields 
By Theorem~\ref{thm:Euler-mean-intersection}, the mod $2$ Euler number
of $\nabla\BB_{av}$ satisfies
$$
   \chi(\nabla\BB_{av}) = 1.
$$
Together with the previous displayed equation this shows that
$\chi(\nabla\BB_{in})=1$, so $\nabla\BB_{in}$ possesses a zero
whose Levi-Civita transform is the desired symmetric frozen planet orbit. 
This concludes the proof of Theorem~\ref{thm:existence}.

%%%%%%%%%%%%%%%%%%%%%%%%%%%%%%%%%%%%%%%%%%%%%%%%%%%%%%%%%%%%%%%%%%%%%%%%%%
\section{Hamiltonian formulation}\label{sec:Ham}
%%%%%%%%%%%%%%%%%%%%%%%%%%%%%%%%%%%%%%%%%%%%%%%%%%%%%%%%%%%%%%%%%%%%%%%%%%

In this section we present the Hamiltonian formulations of the
problems described in the previous sections. They will all be derived
from a general result proved in the first subsection. 

%%%
\subsection{Legendre transform}\label{ss:Legendre-transform}
%%%

In this subsection we describe an abstract Legendre transform which
will be applied to the helium atom in the following subsections. 
For concreteness we restrict to the case that the configuration
space is $\R^n$, but everything could be easily extended to more general
configuration manifolds. 

For $k\in\N_0$ we abbreviate $H^k:=H^k(S^1,\R^n)$. We denote the
derivative of $q\in H^1$ by $\dot q$. Suppose we are given
an open subset $\UU^1\subset H^1$ and a Lagrange function
$$
  \LL: \UU^1\times H^0\to\R,\qquad (q,v)\mapsto \LL(q,v).
$$
We say that $\LL$ possesses a {\em continuous $L^2$-gradient} if $\LL$ 
  is of class $C^1$ and there exists a continuous map
$$
  \nabla\LL = (\nabla_1\LL,\nabla_2\LL):\UU^1\times H^0\to H^0\times H^0
$$
uniquely defined by the conditions
$$
  \la\nabla_i\LL(q,v),w\ra = D_i\LL(q,v)w \quad \text{for all }w\in H^1,
$$
where $\la\ ,\ \ra$ is the $L^2$-inner product and $D_i\LL$ denotes
the derivative with respect to the $i$-th variable.
%Moreover, we make the following regularity assumptions:
%\begin{description}
%\item[(L1)] $\nabla_1\LL:\UU^1\to H^0$ is continuous; 
%\item[(L2)] $\nabla_1\LL:\UU^1\to H^0$ is of class $C^1$. 
%\end{description}
We associate to such $\LL$ its {\em Lagrangian action}
$$
  \mathcal{S}_\LL: \UU^1\to\R,\qquad q\mapsto \LL(q,\dot q).
$$
This is a $C^1$-function whose Fr\'echet derivative at $q\in\UU^1$ in direction $w\in H^1$ is
\begin{align*}
  D\mathcal{S}_\LL(q)w
  &= D_1\LL(q,\dot q)w + D_2\LL(q,\dot q)\dot w 
  = \la\nabla_1\LL(q,\dot q),w\ra + \la\nabla_2\LL(q,\dot q),\dot w\ra,
\end{align*}
%and $\nabla_i\LL$ the corresponding component of the $L^2$-gradient. 
It follows that $q\in\UU^1$ is a critical point of $\mathcal{S}$ if
and only if $\nabla_2\LL(q,\dot q)\in H^1$ and the following
{\em Euler-Lagrange equation} holds:
\begin{equation}\label{eq:EL}
  \frac{d}{dt}\nabla_2\LL(q,\dot q) = \nabla_1\LL(q,\dot q).
\end{equation}
Let us impose the following condition on $\LL$:
\begin{description}
\item[(L)] {\em There exists a differentiable map
$$
  F:\UU^1\times H^0\to H^0,\qquad (q,p)\mapsto F(q,p)
$$
such that for each $q\in\UU^1$ the map $H^0\rightarrow H^0$, 
$v\mapsto \nabla_2\LL(q,v)$ is
a homeomorphism with inverse $p\mapsto F(q,p)$.} 
\end{description}
In particular, the map $F$ then satisfies
\begin{equation}\label{eq:F}
  \nabla_2\LL\bigl(q,F(q,p)\bigr) = p.
\end{equation}
Then we associate to $\LL$ its {\em fibrewise Legendre transform} 
$$
  \HH:\UU^1\times H^0\to\R,\qquad \HH(q,p):=\la p,F(q,p)\ra - \LL(q,F(q,p)).
$$
Using equation~\eqref{eq:F} we compute for $w\in H^1$:
\begin{eqnarray*}
  && D_1\HH(q,p)w \\
  &=& \la p,D_1F(q,p)w\ra - D_1\LL(q,F(q,p))w -
  D_2\LL(q,F(q,p))D_1F(q,p)w \\
  &=& \la p-\nabla_2\LL(q,F(q,p)),D_1F(q,p)w\ra
  - \la\nabla_1\LL(q,F(q,p)),w\ra \\
  &=& - \la\nabla_1\LL(q,F(q,p)),w\ra\,, \\
  && D_2\HH(q,p)w \\
  &=& \la F(q,p),w\ra + \la p,D_2F(q,p)w\ra -
  D_2\LL(q,F(q,p))D_2F(q,p)w \\
  &=& \la F(q,p),w\ra
  + \la p - \nabla_2\LL(q,F(q,p)),D_2F(q,p)w\ra \\
  &=& \la F(q,p),w\ra\,.
\end{eqnarray*}
This shows that $\HH$ has a continuous $L^2$-gradient which is related
to that of $\LL$ by
\begin{equation}\label{eq:H-L}
\begin{aligned}
  \nabla_1\HH(q,p) &= -\nabla_1\LL(q,F(q,p)),\cr
  \nabla_2\HH(q,p) &= F(q,p).
\end{aligned}
\end{equation}
On the other hand, to any Hamilton function $\HH:\UU^1\times
H^0\to\R$ with continuous $L^2$-gradient we can associate its {\em Hamiltonian action}
$$
  \A_\HH:\UU^1\times H^0\to\R,\qquad
  \A_\HH(q,p) := \la p,\dot q\ra - \HH(q,p). 
$$
Its derivatives in direction $w\in H^1$ are given by
\begin{align*}
  D_1\A_\HH(q,p)w
  &= \la p,\dot w\ra - \la\nabla_1\HH(q,p),w\ra,\cr
  D_2\A_\HH(q,p)w
  &= \la \dot q,w\ra - \la\nabla_2\HH(q,p),w\ra. 
\end{align*}
It follows that $(q,p)\in\UU^1\times H^0$ is a critical point of $\A_\HH$ if
and only if $p\in H^1$ and the following {\em Hamilton equations} hold:
\begin{equation}\label{eq:Ham}
\begin{aligned}
  \dot p &= -\nabla_1\HH(q,p),\cr
  \dot q &= \nabla_2\HH(q,p).
\end{aligned}
\end{equation}

\begin{prop}\label{prop:Legendre-transform}
Let $\LL:\UU^1\times H^0\to\R$ be a Lagrange function with continuous gradient
satisfying condition (L) and $\HH:\UU^1\times H^0\to\R$ its fibrewise
Legendre transform. Then the assignments $(q,p)\mapsto q$ and
$q\mapsto (q,p=\nabla_2\LL(q,\dot q))$ define a one-to-one
correspondence between critical points $(q,p)$ of $\A_\HH$ and
critical points $q$ of $\mathcal{S}_\LL$. 
\end{prop}

{\bf Proof: }
Let $q\in\UU^1$ be a critical point of $\mathcal{S}_\LL$, so
$p:=\nabla_2\LL(q,\dot q)\in H^1$ and $q$ solves~\eqref{eq:EL}. 
Then condition (L) and the second equation in~\eqref{eq:H-L} give 
$$
   \dot q = F(q,p) = \nabla_2\HH(q,p), 
$$
and~\eqref{eq:EL} and the first equation in~\eqref{eq:H-L} give
$$
  \dot p = \frac{d}{dt}\nabla_2\LL(q,\dot q) = \nabla_1\LL(q,\dot q) =
  -\nabla_1\LL(q,F(q,p)) = \nabla_1\HH(q,p). 
$$
So $(q,p)$ solves~\eqref{eq:Ham} and is therefore a critical point of $\A_\HH$.

Conversely, let $(q,p)\in\UU^1\times H^0$ be a critical point of $\A_\HH$,
so $p\in H^1$ and $(q,p)$ solves~\eqref{eq:Ham}. 
Then the second equation in~\eqref{eq:Ham} and the second equation in~\eqref{eq:H-L} give 
$$
   \dot q = \nabla_2\HH(q,p) = F(q,p),
$$
which by condition (L) implies $\nabla_2\LL(q,\dot q)=p\in H^1$. 
Now the first equation in~\eqref{eq:Ham} and the first equation
in~\eqref{eq:H-L} give
$$
  \frac{d}{dt}\nabla_2\LL(q,\dot q) = \dot p = \nabla_1\HH(q,p)
  = -\nabla_1\LL(q,F(q,p)) = \nabla_1\LL(q,\dot q). 
$$
So $q$ solves~\eqref{eq:EL} and is therefore a critical point of $\mathcal{S}_\LL$.
\hfill $\square$

\begin{rem}
Inspection of the preceding proof shows that formulae for the
derivatives of $\HH$ do not involve derivatives of $F$. This suggests
that Proposition~\ref{prop:Legendre-transform} should still hold if
in condition (L) we only assume continuity of $F$ rather than
differentiability. 
\end{rem}

\begin{ex}\label{ex:density}
Classically, the Lagrangian $\LL:\UU^1\times H^0\to\R$ has the form
$$
  \LL(q,v) = \int_0^1 L(q(t),v(t))dt
$$
for a smooth Lagrangian density $L:U\times\R^n\to\R$, where $U\subset\R^n$
is an open subset and $\UU^1=H^1(S^1,U)$. Then
$$
  \HH(q,p) = \int_0^1 H(q(t),p(t))dt
$$
with $H:U\times\R^n\to\R$ the fibrewise Legendre transform of $L$, and
the Euler-Lagrange and Hamilton's equations take the traditional form
for $i=1,\dots,n$:
$$
  \frac{d}{dt}\frac{\p L}{\p\dot q_i} = \frac{\p L}{\p q_i}
$$
and
$$
  \dot p_i = -\frac{\p H}{\p q_i},\qquad
  \dot q_i = \frac{\p H}{\p p_i}.
$$
\end{ex}

This example covers the instantaneous interaction Lagrangian
$\LL_{in}$ for helium in the original coordinates $q=(q_1,q_2)\in\R_+^2$
away from collisions. The more general setting in
Proposition~\ref{prop:Legendre-transform} will be needed in the
following subsections to deal with the Lagrangians $\BB_{av}$ and $\BB_{in}$
in the new coordinate $z=(z_1,z_2)$, which do not have the form in
Example~\ref{ex:density}. 

%%%
\subsection{The Kepler problem}
%%%

For the Kepler problem, the function $\QQ$ defined
in~\eqref{eq:Q} is the Lagrangian action $\mathcal{S}_\LL$ associated
to the Lagrange function
\begin{equation}\label{eq:L-Kepler}
  \LL: \UU^1\times H^0(S^1,\R)\to\R,\qquad 
  \LL(z,w) = 2\|z\|^2\|w\|^2 + \frac{N}{\|z\|^2}
\end{equation}
with $\UU^1=H^1(S^1,\R)\setminus\{0\}$.
The computation of the differential in Section~\ref{ss:z-to-q-Kepler}
shows that $\LL$ has a continuous $L^2$-gradient. 
The associated momentum $\eta$ is given by
$$
  \eta = \nabla_2\LL(z,w) = 4\|z\|^2w,
$$
which can be solved for $w$ as
$$
  w = \frac{\eta}{4\|z\|^2} = F(z,\eta). 
$$
Note that the map $F$ is smooth in $(z,\eta)$. 
It follows that $\|\eta\|^2=16\|z\|^4\|w\|^2$ and the associated Hamilton
function becomes
\begin{equation}\label{eq:H-Kepler}
  \HH(z,\eta) = \la\eta,F(z,\eta)\ra - \LL(z,F(z,\eta))
  = \frac{\|\eta\|^2}{8\|z\|^2} - \frac{N}{\|z\|^2}, 
\end{equation}
with Hamiltonian action
\begin{equation}\label{eq:AH-Kepler}
  \A_\HH(z,\eta) = \la\eta,z'\ra
  - \frac{\|\eta\|^2}{8\|z\|^2} + \frac{N}{\|z\|^2}, 
\end{equation}
By Proposition~\ref{prop:Legendre-transform}, critical points of
$\A_\HH$ are in one-to-one correspondence to critical points of $\QQ$.

%%%
\subsection{Mean interaction}
%%%

For the helium atom with mean interaction, the function $\BB_{av}$ defined
in~\eqref{eq:Bav} is the Lagrangian action $\mathcal{S}_{\LL_{av}}$ associated
to the Lagrange function
\begin{equation*}%\label{eq:L-Kepler}
  \LL_{av}: \HH^1_{av}\times H^0(S^1,\R^2)\to\R,\qquad 
  \LL_{av}(z,w) = \LL(z_1,w_1) + \LL(z_2,w_2) + \A(z_1,z_2)
\end{equation*}
where $\HH^1_{av}$ and $\A$ are defined in~\eqref{eq:Hav}
and~\eqref{eq:A}, and $\LL$ is the Kepler Lagrangian
from~\eqref{eq:L-Kepler} with charge $N=2$.
The computation of the differential in Section~\ref{ss:diffav}
shows that $\LL_{av}$ has a continuous $L^2$-gradient. 
Since the interaction term $\A$ does not depend on the $w_i$, 
the associated momenta $\eta_i$ are given as in the Kepler case by
$$
  \eta_i = 4\|z_i\|^2w_i,\qquad 
  w_i = \frac{\eta_i}{4\|z_i\|^2} = F(z_i,\eta_i) 
$$
and the associated Hamilton function becomes
\begin{eqnarray*}%\label{eq:H-Kepler}
  \HH_{av}(z,\eta)
  &=& \HH(z_1,\eta_1) + \HH(z_2,\eta_2) - \A(z_1,z_2) \\
  &=& \sum_{i=1}^2\bigg(\frac{\|\eta_i\|^2}{8\|z_i\|^2} -
  \frac{2}{\|z_i\|^2}\bigg)
  + \frac{\|z_1\|^2\|z_2\|^2}{\|z_1^2\|^2\|z_2\|^2 - \|z_2^2\|^2\|z_1\|^2}, 
\end{eqnarray*}
with Hamiltonian action
\begin{eqnarray*}%\label{eq:AH-Kepler}
  \A_{\HH_{av}}(z,\eta)
  &=& \la\eta_1,z_1'\ra + \la\eta_2,z_2'\ra - \HH_{av}(z,\eta) \\
  &=& \sum_{i=1}^2\bigg(\la\eta_1,z_1'\ra - \frac{\|\eta_i\|^2}{8\|z_i\|^2} +
  \frac{2}{\|z_i\|^2}\bigg)
  - \frac{\|z_1\|^2\|z_2\|^2}{\|z_1^2\|^2\|z_2\|^2 - \|z_2^2\|^2\|z_1\|^2}.
\end{eqnarray*}
By Proposition~\ref{prop:Legendre-transform}, critical points of
$\A_{\HH_{av}}$ are in one-to-one correspondence to critical points of $\BB_{av}$.

%%%
\subsection{Instanteneous interaction}
%%%

For the helium atom with instantaneous interaction, the function $\BB_{in}$ defined
in~\eqref{eq:Bin} is the Lagrangian action $\mathcal{S}_{\LL_{in}}$ associated
to the Lagrange function
\begin{equation*}%\label{eq:L-Kepler}
  \LL_{in}: \HH^1_{in}\times H^0(S^1,\R^2)\to\R,\qquad 
  \LL_{in}(z,w) = \LL(z_1,w_1) + \LL(z_2,w_2) + \II(z_1,z_2)
\end{equation*}
where $\HH^1_{in}$ and $\II$ are defined in~\eqref{eq:Hin}
and~\eqref{eq:I}, and $\LL$ is the Kepler Lagrangian
from~\eqref{eq:L-Kepler} with charge $N=2$.
The computation of the differential in Section~\ref{ss:diffin}
shows that $\LL_{in}$ has a continuous $L^2$-gradient. 
Since the interaction term $\II$ does not depend on the 
$w_i$, 
the associated momenta $\eta_i$ are given as in the Kepler case by
$$
  \eta_i = 4\|z_i\|^2w_i,\qquad 
  w_i = \frac{\eta_i}{4\|z_i\|^2} = F(z_i,\eta_i) 
$$
and the associated Hamilton function becomes
\begin{eqnarray*}%\label{eq:H-Kepler}
  \HH_{in}(z,\eta)
  &=& \HH(z_1,\eta_1) + \HH(z_2,\eta_2) - \II(z_1,z_2) \\
  &=& \sum_{i=1}^2\bigg(\frac{\|\eta_i\|^2}{8\|z_i\|^2} -
  \frac{2}{\|z_i\|^2}\bigg)
  + \int_0^1\frac{1}{z_1^2(\tau_{z_1}(t)) - z_2^2(\tau_{z_2}(t))}dt, 
\end{eqnarray*}
with Hamiltonian action
\begin{eqnarray*}%\label{eq:AH-Kepler}
  \A_{\HH_{in}}(z,\eta)
  &=& \la\eta_1,z_1'\ra + \la\eta_2,z_2'\ra - \HH_{in}(z,\eta) \\
  &=& \sum_{i=1}^2\bigg(\la\eta_1,z_1'\ra - \frac{\|\eta_i\|^2}{8\|z_i\|^2} +
  \frac{2}{\|z_i\|^2}\bigg)
  - \int_0^1\frac{1}{z_1^2(\tau_{z_1}(t)) - z_2^2(\tau_{z_2}(t))}dt\,.
\end{eqnarray*}
By Proposition~\ref{prop:Legendre-transform}, critical points of
$\A_{\HH_{in}}$ are in one-to-one correspondence to critical points of
$\BB_{in}$.

\appendix

%%%%%%%%%%%%%%%%%%%%%%%%%%%%%%%%%%%%%%%%%%%%%%%%%%%%%%%%%%%%%%%%%%%%%%%%%%
\section{Differentiability and Fredholm property}\label{sec:diff}
%%%%%%%%%%%%%%%%%%%%%%%%%%%%%%%%%%%%%%%%%%%%%%%%%%%%%%%%%%%%%%%%%%%%%%%%%%

Let $W\subset H^1(S^1,\mathbb{R}^2)$ be an open subset 
and 
$$
\mathcal{B}:W\longrightarrow\R
$$
be a Fr\'echet differentiable map. We say that $\mathcal{B}$ possesses an {\em $L^2$-gradient} 
if for each $z\in W$ the derivative 
$D\mathcal{B}(z):H^1(S^1,\R^2)\to\R$
extends to a continuous linear functional 
$L^2(S^1,\R^2)\to\R$. In which case the $L^2$-gradient
$$
\nabla\mathcal{B}:W\longrightarrow L^2(S^1,\R^2)
$$ 
uniquely defined by the condition
$$
  \la\nabla\mathcal{B}(z),v\ra = D\mathcal{B}(z)(v)\quad \text{for all }v\in L^2(S^1,\R^2),
$$
where $\la\ ,\ \ra$ is the $L^2$-inner product and 
$D\mathcal{B}$ denotes the Fr\'echet derivative.

The natural domain of definition for such a functional $\mathcal{B}$
is an open subset of $H^1$. 
However, for the applications in this article we are only interested
in functionals that possess $L^2$-gradients of class $C^1$. Since for
the functionals we consider this is only the case when their domain of
definition is restricted to $H^2$, we restrict the following
discussion to subsets of $H^2$.

Recall the setup from Section~\ref{sec:inst-interaction}. The set
\begin{align*}
\mathcal{H}_{in} = 
\Bigg\{&z=(z_1,z_2) \in H^2(S^1,\mathbb{R}^2)| \\
&||z_1||>0,\,\,||z_2||>0,\,\,
 z_1^2(\tau)-z_2^2(\tau_{z_2}(t_{z_1}(\tau)))>0 \,\,\text{for all}\,\,\tau\in S^1\Bigg\}
\end{align*}
is an open subset of the Hilbert space $H^2(S^1,\R^2)=W^{2,2}(S^1,\R^2)$ 
and equipped with the $H^2$-topology. We consider the instantaneous interaction functional
$$
   \mathcal{B}_{in}: \mathcal{H}_{in} \to \R,\qquad z=(z_1,z_2)\mapsto
   Q(z_1,z_2) + I(z_1,z_2)
$$
with the free (non-interacting) term
$$
   Q(z_1,z_2) = 2\sum_{i=1}^2\Bigg(||z_i||^2 \cdot ||z_i'||^2
   +\frac{1}{||z_i||^2}\Bigg) 
$$
and the instantaneous interaction term 
\begin{eqnarray*}
   I(z_1,z_2)
&=&-\frac{1}{||z_1||^2}\int_0^1\frac{z_1(\tau)^2}{z_1^2(\tau)-z_2^2
(\tau_{z_2}(t_{z_1}(\tau)))}d\tau\\
&=&\frac{1}{||z_2||^2}\int_0^1\frac{z_2(\tau)^2}{z_2^2(\tau)-z_1^2(\tau_{z_1}(t_{z_2}(\tau)))}d\tau.
\end{eqnarray*}
For each $z\in\HH_{in}$ the derivative $D\BB_{in}(z):H^2(S^1,\R^2)\to\R$
extends to a continuous linear functional $L^2(S^1,\R^2)\to\R$ (see
the formulas for $DQ(z)$ in the proof of Theorem~\ref{thm:Fredholm}
and for $DI(z)$ in equation~\eqref{eq:derivative} below), so we
can define its $L^2$-gradient as above.
%$\nabla\BB_{in}(z)\in L^2(S^1,\R^2)$ 
%in terms
%of the $L^2$ inner product $\la\ ,\ \ra$ by
%$$
%   \la\nabla\BB_{in}(z),v\ra = D\BB_{in}(z)v\qquad\text{for all }v\in L^2(S^1,\R^2).
%$$
An analogous discussion applies to the mean interaction functional 
$$
   \mathcal{B}_{av}: \mathcal{H}_{av} \to \R,\qquad z=(z_1,z_2)\mapsto
   Q(z_1,z_2) + \mathcal{A}(z_1,z_2)
$$
from Section~\ref{sec:inst-interaction} with the mean interaction term
\begin{equation*}%\label{eq:A}
  \mathcal{A}(z_1,z_2) = - \frac{\|z_1\|^2\|z_2\|^2}{\|z_1^2\|^2\|z_2\|^2 - \|z_2^2\|^2\|z_1\|^2},
\end{equation*}
defined on the open subset
\begin{equation*}%\label{eq:Hav}
   \mathcal{H}_{av} = \Bigg\{z=(z_1,z_2) \in
   H^2(S^1,\mathbb{R}^2)\;\Bigl|\; ||z_1||>0,\,\,||z_2||>0,\,\, 
   \frac{||z_1^2||^2}{||z_1||^2}>\frac{||z_2^2||^2}{||z_2||^2}\Bigg\}
\end{equation*}
of the Hilbert space $H^2(S^1,\R^2)$. Intersecting 
\eqref{eq:inavincl} with $H^2(S^1,\mathbb{R}^2)$ 
gives us 
$$
\HH_{in}\subset\HH_{av}.
$$
The goal of this appendix is the proof of

\begin{thm}\label{thm:Fredholm}
On neighbourhoods of their respective zero sets, the $L^2$-gradients
$\nabla\BB_{in}:\HH_{in}\to L^2(S^1,\R^2)$ and
$\nabla\BB_{av}:\HH_{av}\to L^2(S^1,\R^2)$ as well as their
interpolation $\nabla\BB_r=(1-r)\nabla\BB_{in}+r\nabla\BB_{av}: \HH_{av}\supset\HH_{in}\to L^2(S^1,\R^2)$
are $C^1$-Fredholm maps of index zero. The same holds for their
restrictions to the spaces of symmetric orbits
$\wh\HH^2_{av}\to L_{sym}^2(S^1,\R^2)$ resp.~$\wh\HH^2_{in}\to L_{sym}^2(S^1,\R^2)$
defined in~\eqref{eq:Hin-symm}. 
%are $C^1$-Fredholm maps of index zero on neighbourhoods of their respective zero sets.
\end{thm}

{\bf Proof: }
The derivative of the free term applied to $v=(v_1,v_2)\in
L^2(S^1,\R^2)$ is 
\begin{align*}
  DQ(z)v
  &= 4\sum_{i=1}^2\Bigl(||z_i||^2\la z_i',v_i'\ra + ||z_i'||^2\la z_i,v_i\ra
   - \frac{1}{||z_i||^4}\la z_i,v_i\ra\Bigr) \cr
  &= 4\sum_{i=1}^2\Bigl\la -||z_i||^2z_i'' + ||z_i'||^2 z_i 
   - \frac{1}{||z_i||^4} z_i,v_i\Bigr\ra,
\end{align*}
hence its $L^2$-gradient has components
\begin{equation}\label{eq:gradQ}
   \nabla_i Q(z) = 4\Bigl(-||z_i||^2z_i'' + ||z_i'||^2 z_i 
   - \frac{1}{||z_i||^4} z_i\Bigr),\qquad i=1,2.
\end{equation}
This obviously defines a $C^1$-map $\nabla Q:\HH_{av}\to L^2(S^1,\R^2)$. To see that the derivative of $\nabla Q$
at $w=(w_1,w_2)$ is Fredholm of index zero we write the leading term of $\nabla Q$ near $w$ as 
$z\mapsto -4(||w_1||^2z_1'',||w_2||^2z_2'')$. The latter
map is a restriction of the obvious linear index zero 
Fredholm map $H^2(S^1,\R^2)\to L^2(S^1,\R^2)$ (the kernel is spanned and the image is complemented by constants).
The remaining lower order terms give a compact perturbation, so $\nabla Q$ is a 
nonlinear $C^1$-Fredholm map $\HH_{av}\to L^2(S^1,\R^2)$ of index zero.  

For the mean interaction, the components of the $L^2$-gradient are
read off from~\eqref{eq:DA} to be
\begin{eqnarray}\label{eq:gradA}
  \nabla_1\mathcal{A}[z_1,z_2] \nonumber
  &=& -2\frac{||z_2||^4 \cdot ||z_1^2||^2}{\big(||z_1^2||^2\cdot ||z_2||^2
    -||z_2^2||^2 \cdot ||z_1||^2\big)^2}z_1 \\ \nonumber
  & &+4\frac{||z_1||^2\cdot ||z_2||^4}{\big(||z_1^2||^2\cdot
      ||z_2||^2-||z_2^2||^2 \cdot ||z_1||^2\big)^2}z_1^3 \\ \nonumber
  \nabla_2\mathcal{A}[z_1,z_2]
  &=& +2\frac{||z_1||^4 \cdot ||z_2^2||^2}{\big(||z_1^2||^2\cdot ||z_2||^2
    -||z_2^2||^2 \cdot ||z_1||^2\big)^2} z_2\\ 
  & &-4\frac{||z_1||^4\cdot ||z_2||^2}{\big(||z_1^2||^2\cdot
      ||z_2||^2-||z_2^2||^2 \cdot ||z_1||^2\big)^2} z_2^3.
\end{eqnarray}
Since these are $C^1$-maps whose derivatives at each point 
are compact linear operators, this proves the assertion for 
$\BB_{av}$ (which also follows from Theorem~\ref{thm:Euler-mean-intersection}). 

\quad Existence of the restriction of $\nabla Q$ to 
$\wh\HH^2_{av} \to L_{sym}^2(S^1,\R^2)$ as a Fredholm map of index zero is straightforward. So Theorem~\ref{thm:Fredholm} follows from Proposition~\ref{prop:compact}
below.
\hfill$\square$\medskip

\begin{prop}\label{prop:compact}
On a neighbourhood of the zero set of $\BB$, the gradient $\nabla I:\HH_{in}\to
L^2(S^1,\R^2)$ of the interaction term defines a $C^1$-map whose
derivative at each point is a compact linear operator
$H^2(S^1,\R^2)\to L^2(S^1,\R^2)$. 
\end{prop}

The proof of this proposition will occupy the rest of this appendix. 
It uses some technical lemmas from Appendix~\ref{sec:contdep}. 

%%%
\subsection{Reparametrizations of the circle}
%%%

Here we collect some facts about reparametrizations of the circle that
are used throughout this article. In this subsection we consider
a map $z\in C^1(S^1,\R)$ with finite zero set
$$
   Z:=\{\tau\in S^1\mid z(\tau)=0\}.
$$
As in Section~\ref{ss:LC-inv}, we associate to $z$ a $C^2$-map
$t_z:S^1\to S^1$ by 
\begin{equation}\label{eq:tz-restated}
   t_z(\tau) := \frac{1}{\|z\|^2}\int_0^\tau z(\sigma)^2d\sigma
\end{equation}
with derivative
\begin{equation}\label{eq:dert-restated}
   t_z'(\tau) = \frac{z(\tau)^2}{\|z\|^2}.
\end{equation}
By Lemma~\ref{lem:tz-homeo}, the map $t_z$ is a homeomorphism with
continuous inverse 
\begin{equation}\label{eq:tau}
   \tau_z := t_z^{-1}:S^1\to S^1.
\end{equation}
Since $t_z$ is of class $C^2$, the function $\tau_z$ is also of class $C^2$ on the
complement of the finite set $t_z(Z)$ with derivative
\begin{equation}\label{eq:taudot-restated}
  \dot\tau_z(t) = \frac{\|z\|^2}{z(\tau_z(t))^2}.
\end{equation}
In Section~\ref{sec:inst-interaction} we need the Fr\'echet
derivatives of $t_z$ and $\tau_z$ with respect to $z$.
The derivative of $t_z$ with respect to $z$ in direction $v\in L^2(S^1,\R)$ is given by
\begin{equation}\label{eq:dt}
   Dt_z(v)(\tau) = \frac{2}{\|z\|^2}\int_0^\tau z(\sigma) v(\sigma) d\sigma
   - \frac{2\langle z,v\rangle}{\|z\|^4}\int_0^\tau z(\sigma)^2 d\sigma\,. 
\end{equation}
For future use observe that $Dt_z$ defines a bounded operator 
$L^2(S^1,\R)\to H^1(S^1,\R)$ depending continuously on $z$. Indeed, the first summand is a composition of multiplication with a continuous function and integrating from $0$ to $\tau$. The multiplication is a continuous operator to $L^2$, and integration is a continuous operator to $H^1$. Boundedness of the second summand can be seen analogously. Continuous dependence on $z$ is clear. 

Using equation~\eqref{eq:dt} and equation~\eqref{eq:dert-restated}, we derive a formula
for the derivative of $\tau_z$: 
\begin{eqnarray*}
  0 &=& D(t_z \circ \tau_z)(v)(t)\\
  &=& Dt_z(v)(\tau_z(t))+t'_z(\tau_z(t))D\tau_z (v)(t)\\
  &=& \frac{2}{\|z\|^2}\int_0^{\tau_z(t)}\hspace{-15pt} z(\sigma) v(\sigma) d\sigma
   - \frac{2\langle z,v\rangle}{\|z\|^4}\int_0^{\tau_z(t)}\hspace{-15pt} z(\sigma)^2 d\sigma
   + \frac{z(\tau_z(t))^2}{||z||^2}D\tau_z(v)(t),
\end{eqnarray*}
thus
\begin{equation}\label{eq:dtau}
  D\tau_z(v)(t)
  = \frac{2\langle z,v\rangle}{\|z\|^2z(\tau_z(t))^2}\int_0^{\tau_z(t)}\hspace{-15pt} z(\sigma)^2 d\sigma
  - \frac{2}{z(\tau_z(t))^2}\int_0^{\tau_z(t)}\hspace{-15pt} z(\sigma) v(\sigma) d\sigma\,.
\end{equation}
The instantaneous interaction term in
Section~\ref{sec:inst-interaction} involves the expression
\begin{equation}\label{eq:Delta}
   \Delta(\tau) := z_2^2(\tau)-z_1^2(\tau_{z_1}(t_{z_2}(\tau))),
\end{equation}
whose derivative with respect to $z_1$ in direction $v\in L^2(S^1,\R)$
is given by
\begin{eqnarray}\label{eq:dDeltaz^2}
(D_1\Delta)(v)(\tau) %&=&-D_1z_1^2(\tau_{z_1}(t_{z_2}(\tau)))\\
&=&-2z_1(\tau_{z_1}(t_{z_2}(\tau)))v(\tau_{z_1}(t_{z_2}(\tau)))\nonumber \\
&&-2z_1(\tau_{z_1}(t_{z_2}(\tau)))z_1'(\tau_{z_1}(t_{z_2}(\tau)))D\tau_{z_1}(v)(t_{z_2}(\tau)).
\label{eq:dDelta}
\end{eqnarray}
We also sometimes need the transformation of $\Delta$ under time change
\begin{equation}\label{eq:Delta-transf}
   \Delta(\tau_{z_2}(t_{z_1}(\tau)))=z_2^2(\tau_{z_2}(t_{z_1}(\tau)))-z_1^2(\tau)<0
\end{equation}
and its derivative
\begin{eqnarray}\nonumber
  && D_1\Bigl(\Delta(\tau_{z_2}(t_{z_1}(\cdot)))\Bigr)(v)(\tau)\\ \nonumber
  &=& 2z_2(\tau_{z_2}(t_{z_1}(\tau)))z_2'(\tau_{z_2}(t_{z_1}(\tau)))\dot\tau_{z_2}(t_{z_1}(\tau))Dt_{z_1}(v)(\tau)
  -2z_1(\tau)v(\tau)\\ 
  &=&2||z_2 ||^2\frac{z_2'(\tau_{z_2}(t_{z_1}(\tau)))}{z_2(\tau_{z_2}(t_{z_1}(\tau)))}Dt_{z_1}(v)(\tau) 
  -2z_1(\tau)v(\tau) \label{eq:Deltach}
\end{eqnarray}

%%%
\subsection{$L^2$-gradient of the instantaneous interaction term}\label{ss:der}
%%%

We write the interaction term as
$$
   I(z_1,z_2) = \frac{1}{||z_2||^2}\II(z_1,z_2)
$$
with
$$
   \II(z_1,z_2)
   := \int_0^1\frac{z_2(\tau)^2}{z_2^2(\tau)-z_1^2(\tau_{z_1}(t_{z_2}(\tau)))}d\tau
   = \int_0^1\frac{z_2(\tau)^2}{\Delta(\tau)}d\tau.
$$
Since $\|z_2\|>0$, it is enough to prove Proposition~\ref{prop:compact} with $\II$ in place of $I$. 
In the remainder of this appendix we will prove compactness and continuous
dependence for the $z_1$-derivative of the $z_1$-component of the
$L^2$-gradient $\nabla\II$; the treatments of the $z_2$-component and
the $z_2$-derivatives of both components are analogous and will be omitted. 

Recall that $\Delta(\tau)$ never vanishes; we will use this without
further mention in the computations below.

{\bf Derivative of $\II$. }
Let us compute the derivative of $\II$ with respect to $z_1$ in
direction $v\in L^2(S^1,\R)$: 
\begin{eqnarray}\nonumber
  D_1\mathcal{I}(z_1,z_2)(v)
  &=&-\int_0^1\frac{z_2^2(\tau)}{\Delta(\tau)^2}D_1\Delta(v)(\tau)d\tau\\ \nonumber
  &=&2\int_0^1\frac{z_2^2(\tau)}{\Delta(\tau)^2}z_1(\tau_{z_1}(t_{z_2}(\tau)))v(\tau_{z_1}(t_{z_2}(\tau)))d\tau\\ \nonumber
  &
  &+2\int_0^1\frac{z_2^2(\tau)}{\Delta(\tau)^2}z_1(\tau_{z_1}(t_{z_2}(\tau)))z_1'(\tau_{z_1}(t_{z_2}(\tau)))D\tau_{z_1}(v)(t_{z_2}(\tau))d\tau\\ \nonumber
  &=&2\int_0^1\frac{z_2^2(\tau)}{\Delta(\tau)^2}z_1(\tau_{z_1}(t_{z_2}(\tau)))v(\tau_{z_1}(t_{z_2}(\tau)))d\tau\\ \nonumber
  & &+ \frac{4\langle z_1,v\rangle}{\|z_1\|^2}
  \int_0^1\frac{z_2^2(\tau)}{\Delta(\tau)^2}
  \frac{z_1'(\tau_{z_1}(t_{z_2}(\tau)))}{z_1(\tau_{z_1}(t_{z_2}(\tau)))}
  \int_0^{\tau_{z_1}(t_{z_2}(\tau))}\hspace{-15pt} z_1(\sigma)^2 d\sigma\, d\tau\\ \nonumber
  & & -4\int_0^1\frac{z_2^2(\tau)}{\Delta(\tau)^2}
  \frac{z_1'(\tau_{z_1}(t_{z_2}(\tau)))}{z_1(\tau_{z_1}(t_{z_2}(\tau)))}
  \int_0^{\tau_{z_1}(t_{z_2}(\tau))}\hspace{-15pt} z_1(\sigma) v(\sigma) d\sigma\,d\tau\\ 
%\int_0^{\tau_{z_1}(t_{z_2}(\tau))}z_1(\sigma)v(\sigma)d\sigma
&=:&\sum_{i=1}^3D_1^i\mathcal{I}(z_1,z_2)(v). \label{eq:derivative}
\end{eqnarray}
Here in the second equality we have used equation~\eqref{eq:dDelta},
in the third equality we have substituted $D\tau_{z_1}(v)$ using 
equation~\eqref{eq:dtau} with $z=z_1$, and we denote the resulting three summands by
$D_1^i\mathcal{I}(z_1,z_2)(v)$, $i=1,2,3$. 

Our next goal is to rewrite $D_1\mathcal{I}(z_1,z_2)(v)$ as the
$L^2$-inner product of $v$ with the first component of the
$L^2$-gradient of $\II$,
$$
  D_1\II(z_1,z_2))v = \la\nabla_1\II(z_1,z_2),v\ra\qquad\text{for all }v\in L^2(S^1,\R).
$$
{\bf Coordinate change in the integrals. } 
In order to write the first term in~\eqref{eq:derivative} as an
$L^2$-inner product with $v$, we perform the following coordinate
change that will also be used later. For $\sigma\in S^1$ we set 
$$
   \xi := \tau_{z_1}(t_{z_2}(\sigma)) \in S^1,
$$ 
so that
$$
   \sigma = \tau_{z_2}(t_{z_1}(\xi))
$$
and from equations~\eqref{eq:dert-restated} and~\eqref{eq:taudot-restated} we get
\begin{equation}\label{eq:dsigma}
  d\sigma = \dot\tau_{z_2}(t_{z_1}(\xi))t_{z_1}'(\xi)d\xi
  = \frac{||z_2||^2}{||z_1||^2}\frac{z_1^2(\xi)}{z_2^2(\sigma)}d\xi.
\end{equation}
So after renaming the integration variable from $\tau$ to $\sigma$ the
first term in~\eqref{eq:derivative} becomes
\begin{align*}
  D_1^1\mathcal{I}(z_1,z_2)(v)
  &= 2\int_0^1\frac{z_2^2(\tau)}{\Delta(\tau)^2}z_1(\tau_{z_1}(t_{z_2}(\sigma)))
  v(\tau_{z_1}(t_{z_2}(\sigma)))d\sigma \cr
  &= 2\frac{||z_2||^2}{||z_1||^2}\int_0^1\frac{z_1^3(\xi)v(\xi)}{\Delta(\tau_{z_2}(t_{z_1}(\xi)))^2}d\xi.
\end{align*}
The second term in~\eqref{eq:derivative} has already the form of an
$L^2$-inner product with $v$. Using equation~\eqref{eq:tz-restated} to insert
$t_z(\tau)$ and renaming the integration variable $\tau$ to $\sigma$
it becomes
\begin{align*}
  D_1^2\mathcal{I}(z_1,z_2)(v)
  &= 4\la z_1,v\ra\int_0^1\frac{z_2^2(\sigma)}{\Delta(\sigma)^2}
  \frac{z_1'(\tau_{z_1}(t_{z_2}(\sigma)))}{z_1(\tau_{z_1}(t_{z_2}(\sigma)))}
  t_{z_1}(\tau_{z_1}(t_{z_2}(\sigma)))d\sigma
\end{align*}
{\bf Switching the order of integration. }
To write the third term as an inner product, we need to switch the
order of integration in the double integral.
The general setup for this is the following. 
Let $f,F:S^1\to\R$ be continuous functions. Let $\tau_1:S^1\to S^1$
be a $C^1$-homeomorphism with $\tau(0)=0$ and finitely many critical
points. Then
\begin{align*}
  \int_0^1d\tau F(\tau)\int_0^{\tau_1(\tau)}f(\sigma)d\sigma
  &= \int_0^1d\sigma f(\sigma)\int_{\tau_1^{-1}(\sigma)}^1F(\tau)d\tau \cr
  &= \int_0^1d\tau f(\tau)\int_{\tau_1^{-1}(\tau)}^1F(\sigma)d\sigma.
\end{align*}
We apply this formula with 
$$
   \tau_1(\tau):=\tau_{z_1}(t_{z_2}(\tau)),\qquad f(\sigma):=z_1(\sigma)v(\sigma)
$$ 
to the third term in~\eqref{eq:derivative}. To deal with the
integration limits observe that the inverse of $\tau_1$ is given by
$\tau_1^{-1}(\tau)=\tau_{z_2}(t_{z_1}(\tau))$. Thus we find
\begin{align*}
  D_1^3\mathcal{I}(z_1,z_2)(v)
  &= -4\int_0^1d\tau\frac{z_2^2(\tau)}{\Delta(\tau)^2}
  \frac{z_1'(\tau_{z_1}(t_{z_2}(\tau)))}{z_1(\tau_{z_1}(t_{z_2}(\tau)))}
  \int_0^{\tau_{z_1}(t_{z_2}(\tau))}\hspace{-15pt}z_1(\sigma)v(\sigma)d\sigma \cr
  &= -4\int_0^1z_1(\tau)v(\tau)\int_{\tau_{z_2}(t_{z_1}(\tau))}^1\frac{z_2^2(\sigma)}{\Delta(\sigma)^2}
\frac{z_1'(\tau_{z_1}(t_{z_2}(\sigma)))}{z_1(\tau_{z_1}(t_{z_2}(\sigma)))}d\sigma.
\end{align*}
Altogether this gives us the formula for the first component of the $L^2$-gradient
\begin{eqnarray} \nonumber
  \nabla_1\mathcal{I}(z_1,z_2)
  &=&2\frac{||z_2||^2}{||z_1||^2}\frac{z_1^3(\tau)}{\Delta(\tau_{z_2}(t_{z_1}(\tau)))^2}\\ \nonumber
  & &+4z_1(\tau)\int_0^1\frac{z_2^2(\sigma)}{\Delta(\sigma)^2}
  \frac{z_1'(\tau_{z_1}(t_{z_2}(\sigma)))}{z_1(\tau_{z_1}(t_{z_2}(\sigma)))}
  t_{z_1}(\tau_{z_1}(t_{z_2}(\sigma)))d\sigma\\ \nonumber
  & &-4z_1(\tau)\int_{\tau_{z_2}(t_{z_1}(\tau))}^1\frac{z_2^2(\sigma)}{\Delta(\sigma)^2}
  \frac{z_1'(\tau_{z_1}(t_{z_2}(\sigma)))}{z_1(\tau_{z_1}(t_{z_2}(\sigma)))}d\sigma\\ 
  &=:&\sum_{i=1}^3 \mathcal{V}^i_1(z_1,z_2)(\tau). \label{eq:gradient-I}
\end{eqnarray}
Note that, since $\Delta$ and $z_1$ never vanish,
$\nabla_1\mathcal{I}(z_1,z_2)$ exists as an $L^2$-function and depends
continuously on $(z_1,z_2)\in\HH_{in}$. 

%%%
\subsection{Hessian part 1}
%%%

In this subsection we consider the first part of the gradient
in \eqref{eq:gradient-I},
$$
   \mathcal{V}^1_1(z_1,z_2)(\tau) =
   2\frac{||z_2||^2}{||z_1||^2}\frac{z_1^3(\tau)}{\Delta(\tau_{z_2}(t_{z_1}(\tau)))^2}.
$$
Differentiating it with respect to $z_1$ in direction $v\in
H^2(S^1,\R)$ we obtain
\begin{eqnarray}\nonumber
  \frac{||z_1||^2}{2||z_2||^2}D_1\mathcal{V}^1_1(z_1,z_2)(v)(\tau)
  &=&-\frac{2}{||z_1||^2}\frac{z_1^3(\tau)}{\Delta(\tau_{z_2}(t_{z_1}(\tau)))^2}\langle z_1,v\rangle
\\ \nonumber
  &&+3\frac{z_1^2(\tau)v(\tau)}{\Delta(\tau_{z_2}(t_{z_1}(\tau)))^2}\\ \nonumber 
  &&-2\frac{z_1^3(\tau)}{\Delta(\tau_{z_2}(t_{z_1}(\tau)))^3}
  D_1\Bigl(\Delta(\tau_{z_2}(t_{z_1}(\cdot)))\Bigr)(v)(\tau)\\ \nonumber
&=&-\frac{2}{||z_1||^2}\frac{z_1^3(\tau)}{\Delta(\tau_{z_2}(t_{z_1}(\tau)))^2}\langle z_1,v\rangle
\\ \nonumber
  &&+3\frac{z_1^2(\tau)v(\tau)}{\Delta(\tau_{z_2}(t_{z_1}(\tau)))^2}\\ \nonumber
  &&-4\frac{z_1^3(\tau)||z_2 ||^2z_2'(\tau_{z_2}(t_{z_1}(\tau)))}{\Delta(\tau_{z_2}(t_{z_1}(\tau)))^3z_2(\tau_{z_2}(t_{z_1}(\tau)))}Dt_{z_1}(v)(\tau)
\\ 
  &&+4\frac{z_1^3(\tau)z_1(\tau)v(\tau)}{\Delta(\tau_{z_2}(t_{z_1}(\tau)))^3},
\end{eqnarray}
where we have used equation~\eqref{eq:Deltach} to replace
$D_1\Bigl(\Delta(\tau_{z_2}(t_{z_1}(\cdot)))\Bigr)(v)(\tau)$. 
%After expanding the bracket we get $4$ summands.
We need to show that each of the four summands on the right hand side
as a function of $v$ defines a compact linear operator $H^2(S^1,\R)\to
L^2(S^1,\R)$ that depends continuously on $(z_1,z_2)\in\HH_{in}$
with respect to the operator norm.

The first summand is a $1$-dimensional operator (hence compact) whose
image is spanned by a function that lies
in $H^1(S^1,\R)$ and depends continuously on $(z_1,z_2)$ by
Lemmas~\ref{lem:H^1cont} and~\ref{lem:Banach-alg}. \\
The second and fourth summands are multiplication operators with
functions that lie in $H^1(S^1,\R)$ and depend continuously on $(z_1,z_2)$
by Lemmas~\ref{lem:H^1cont} and~\ref{lem:Banach-alg}.
%and nonvanishing of $\Delta(\tau_{z_2}(t_{z_1}(\tau)))$. 
They are compact because they are compositions
$$
   H^2(S^1,\R)\into H^1(S^1,\R)\to H^1(S^1,\R)\into L^2(S^1,\R),
$$
where the middle map is multiplication with a fixed $H^1$-function and
the two inclusions are compact. \\
For the third summand first note that by 
formula~\eqref{eq:dt} with
$z=z_1$ the map $v\mapsto Dt_{z_1}(v)$ defines a bounded linear operator
$H^2(S^1,\R)\to H^1(S^1,\R)$ depending continuously on $z_1$. 
The third term is the composition of this operator and several 
multiplication operators. The functions with which we multiply lie in
$H^1(S^1,\R)$ except for $\frac{1}{z_2(\tau_{z_2}(t_{z_1}(\tau)))}$,
which lies in $L^2(S^1,\R)$. They depend continuously on $(z_1,z_2)$
by Lemmas~\ref{lem:H^1cont},~\ref{lem:Banach-alg} and~\ref{lem:z_1z_2L^2cont}.  
To show compactness of this operator we write it as the composition
of continuous linear maps
$$
   H^2(S^1,\R)\rightarrow H^1(S^1,\R)
   %\stackrel{A}{\rightarrow}H^{2/3}(S^1,\R)\stackrel{B}{\rightarrow}
   \into C^0(S^1,\R)\rightarrow L^2(S^1,\R),
$$
where the first map sends $v\mapsto Dt_{z_1}(v)$, the third map is
multiplication with $\frac{1}{z_2(\tau_{z_2}(t_{z_1}(\tau)))}$, 
and the canonical inclusion in the middle is compact by the Rellich
embedding theorem. 
This concludes the discussion of $\mathcal{V}^1_1$. 

%%%
\subsection{Hessian part 2}
%%%

The second part of the gradient in~\eqref{eq:gradient-I} has the form
$$
   \mathcal{V}^2_1(z_1,z_2)(\tau) =
   4z_1(\tau)\int_0^1g(z_1,z_2)(\sigma)t_{z_2}(\sigma)d\sigma
$$
with
$$
   g(z_1,z_2)(\sigma) := \frac{z_2^2(\sigma)}{\Delta(\sigma)^2}
  \frac{z_1'(\tau_{z_1}(t_{z_2}(\sigma)))}{z_1(\tau_{z_1}(t_{z_2}(\sigma)))}.
$$
Differentiating it with respect to $z_1$ in direction $v\in
H^2(S^1,\R)$ we obtain
\begin{eqnarray*}
  \frac{1}{4}D_1\mathcal{V}^2_1(z_1,z_2)(v)(\tau)
  &=& v(\tau)\int_0^1g(z_1,z_2)(\sigma)t_{z_2}(\sigma)d\sigma \cr
  && + z_1(\tau)\int_0^1D_1g(z_1,z_2)(v)(\sigma)t_{z_2}(\sigma)d\sigma \cr
  && + z_1(\tau)\int_0^1g(z_1,z_2)(\sigma)Dt_{z_2}(v)(\sigma)d\sigma.
\end{eqnarray*}
Note that $g(z_1,z_2)(\sigma)$ agrees with the integrand of the third
part of the gradient in~\eqref{eq:gradient-I}. It is shown in the next subsection that $g(z_1,z_2)\in
C^0(S^1,\R)$ depends continuously on $(z_1,z_2)$, and $v\mapsto
D_1g(z_1,z_2)(v)$ defines a bounded linear operator
$H^2(S^1,\R)\to L^2(S^1,\R)$ that depends continuously on $(z_1,z_2)$.
Equations~\eqref{eq:tz-restated} and~\eqref{eq:dt} show that $t_{z_2}\in
C^0(S^1,\R)$ depends continuously on $z_2$, and $v\mapsto
Dt_{z_2}(v)$ defines a bounded linear operator
$H^2(S^1,\R)\to L^2(S^1,\R)$ that depends continuously on $z_2$.
This shows that the right hand side as a function of $v$ defines a bounded
linear operator $H^2(S^1,\R)\to L^2(S^1,\R)$ that depends continuously
on $(z_1,z_2)\in\HH_{in}$. It is compact because the first term is a
scalar multiplication operator composed with the compact inclusion
$H^2(S^1,\R)\into L^2(S^1,\R)$, and the other two terms have
$1$-dimensional images. 

%%%
\subsection{Hessian part 3}
%%%

In this subsection we consider the third part of the gradient
in~\eqref{eq:gradient-I},
$$
   \mathcal{V}^3_1(z_1,z_2)(\tau) =
   -4z_1(\tau)\int_{\tau_{z_2}(t_{z_1}(\tau))}^1\frac{z_2^2(\sigma)}{\Delta(\sigma)^2}
   \frac{z_1'(\tau_{z_1}(t_{z_2}(\sigma)))}{z_1(\tau_{z_1}(t_{z_2}(\sigma)))}d\sigma.
$$
Differentiating it with respect to $z_1$ in direction $v\in
H^2(S^1,\R)$ we obtain
\begin{eqnarray*}
  -\frac{1}{4}D_1\mathcal{V}^3_1(z_1,z_2)(v)(\tau)
  &=&v(\tau)\int_{\tau_{z_2}(t_{z_1}(\tau))}^1\frac{z_2^2(\sigma)}{\Delta(\sigma)^2}\frac{z_1'(\tau_{z_1}(t_{z_2}(\sigma)))}{z_1(\tau_{z_1}(t_{z_2}(\sigma)))}d\sigma
  \\& &
  -z_1(\tau)\frac{z_2^2(\tau_{z_2}(t_{z_1}(\tau)) )}{\Delta(\tau_{z_2}(t_{z_1}(\tau)))^2}
  \frac{z_1'(\tau)}{z_1(\tau)}
%  \frac{||z_2||^2}{z_2^2(\tau_{z_2}(t_{z_1}(\tau)))}
  \dot\tau_{z_2}(t_{z_1}(\tau))
  Dt_{z_1}(v)(\tau)
 \\ 
\nonumber
& &
+z_1(\tau)\int_{\tau_{z_2}(t_{z_1}(\tau))}^1z_2^2(\sigma)D_1\left(\frac{1}{\Delta(\sigma)^2}\right)(v)
\frac{z_1'(\tau_{z_1}(t_{z_2}(\sigma)))}{z_1(\tau_{z_1}(t_{z_2}(\sigma)))}d\sigma\\ 
\nonumber
& &
+z_1(\tau)\int_{\tau_{z_2}(t_{z_1}(\tau))}^1\frac{z_2^2(\sigma)}{\Delta(\sigma)^2}D_1\left(\frac{z_1'(\tau_{z_1}(t_{z_2}(\sigma)))}{z_1(\tau_{z_1}(t_{z_2}(\sigma)))}\right)(v)d\sigma\\ 
\nonumber
&=:&\sum_{i=1}^4T_iv(\tau).
\end{eqnarray*}
Again, we need to show that each $T_i$ defines 
a compact linear operator $H^2(S^1,\R)\to
L^2(S^1,\R)$ that depends continuously on $(z_1,z_2)\in\HH_{in}$. 
In the integrals we will use the change of variables $\xi := \tau_{z_1}(t_{z_2}(\sigma))$
from Section~\ref{ss:der}. The integrand transforms according to
formula~\eqref{eq:dsigma}. To change the limits of integration note that 
$\sigma\in [\tau_{z_2}(t_{z_1}(\tau)),1]$ corresponds to $\xi\in [\tau,1]$.

Now we discuss the four terms one by one. We will omit the arguments
$(S^1,\R)$ and simply write $H^2$ instead of $H^2(S^1,\R)$ etc. 

{\bf The first term. }
Change of integration variable turns the first term into
\begin{align*}
  T_1v(\tau)
  &= v(\tau)\int_{\tau_{z_2}(t_{z_1}(\tau))}^1\frac{z_2^2(\sigma)}{\Delta(\sigma)^2}
  \frac{z_1'(\tau_{z_1}(t_{z_2}(\sigma)))}{z_1(\tau_{z_1}(t_{z_2}(\sigma)))}d\sigma \cr
  &= v(\tau)\frac{||z_2||^2}{||z_1||^2}\int_{\tau}^1\frac{z_1'(\xi)z_1(\xi)}
  {\Delta(\tau_{z_2}(t_{z_1}(\xi)))^2}d\xi.
\end{align*}
By Lemmas~\ref{lem:H^1cont} and~\ref{lem:Banach-alg} the integrand is
continuous and depends continuously on $(z_1,z_2)$ as a map
$\HH_{in}\to C^0$. Therefore, $T_1$ defines a bounded linear
operator $H^2\to H^1$ depending continuously on $(z_1,z_2)$,
and composition with the inclusion $H^1\into L^2$ makes it compact. 

{\bf The second term. }
Using equation~\eqref{eq:tz-restated} to replace $\dot\tau_{z_2}(t_{z_1}(\tau))$
turns the seconds term into
\begin{align*}
  T_2v(\tau)
  &= -z_1(\tau)\frac{z_2^2(\tau_{z_2}(t_{z_1}(\tau)) )}{\Delta(\tau_{z_2}(t_{z_1}(\tau)))^2}
  \frac{z_1'(\tau)}{z_1(\tau)}
  \frac{||z_2||^2}{z_2^2(\tau_{z_2}(t_{z_1}(\tau)))}
  Dt_{z_1}(v)(\tau) \cr
  &= -\frac{z_1'(\tau)||z_2||^2}{\Delta(\tau_{z_2}(t_{z_1}(\tau)))^2}Dt_{z_1}(v)(\tau).
\end{align*}
Formula~\eqref{eq:dt} with $z=z_1$ shows that $v\mapsto Dt_{z_1}(v)$
defines a bounded operator to $H^2\to H^1$ depending continuously on $(z_1,z_2)$.
Together with Lemmas~\ref{lem:H^1cont} and~\ref{lem:Banach-alg} this
implies that $T_2$ defines a bounded operator $H^2\to H^1$
depending continuously on $(z_1,z_2)$, 
and composition with the inclusion $H^1\into L^2$ makes it compact. 

{\bf The third term. }
We rewrite the third term as
\begin{align*}
  T_3v(\tau)
  &= z_1(\tau)\int_{\tau_{z_2}(t_{z_1}(\tau))}^1z_2^2(\sigma)
  \frac{-2(D_1\Delta)(v)(\sigma)}{\Delta(\sigma)^3} 
  \frac{z_1'(\tau_{z_1}(t_{z_2}(\sigma)))}{z_1(\tau_{z_1}(t_{z_2}(\sigma)))}d\sigma\cr
  &= -2z_1(\tau)\int_{\tau_{z_2}(t_{z_1}(\tau))}^1\frac{z_2^2(\sigma)}{\Delta(\sigma)^3}
  \frac{z_1'(\tau_{z_1}(t_{z_2}(\sigma)))}{z_1(\tau_{z_1}(t_{z_2}(\sigma)))}
  \Bigl[-2z_1(\tau_{z_1}(t_{z_2}(\tau)))v(\tau_{z_1}(t_{z_2}(\tau)))\cr
  &
  \ \ \ \ \ \ -2z_1(\tau_{z_1}(t_{z_2}(\tau)))z_1'(\tau_{z_1}(t_{z_2}(\tau)))
  D\tau_{z_1}(v)(t_{z_2}(\sigma))\Bigr]d\sigma\cr
  &= 4\frac{||z_2||^2}{||z_1||^2}z_1(\tau)\int_\tau^1
  \frac{z_1(\xi)z_1'(\xi)}{\Delta(\tau_{z_2}(t_{z_1}(\xi)))^3}
  \Bigl[z_1(\xi)v(\xi)+z_1(\xi)z_1'(\xi)Xv(\xi)\Bigr]d\xi.
\end{align*}
Here in the second equality we have used equation~\eqref{eq:dDeltaz^2}
to replace $(D_1\Delta)(v)(\sigma)$, in the third equation we change
the integration variable, and we have abbreviated (replacing
$D\tau_{z_1}(v)$ via equation~\eqref{eq:dtau})
\begin{equation}\label{eq:X}
  Xv(\xi)
  := D\tau_{z_1}(v)(t_{z_1}(\xi))
  =\frac{2}{z_1^2(\xi)}\left(\frac{||z_1||^2(\xi)}{||z_1||^2}
  \langle z_1,v\rangle-\langle z_1,v\rangle(\xi)\right).
\end{equation}
The map $v\mapsto Xv$ defines a bounded operator $H^2\to H^1$
depending continuously on $(z_1,z_2)$.
%This observation will also be used in the discussion of the fourth term.
Together with Lemmas~\ref{lem:H^1cont} and~\ref{lem:Banach-alg} this
implies that $T_3$ defines a bounded operator $H^2\to H^1$
depending continuously on $(z_1,z_2)$, 
and composition with the inclusion $H^1\into L^2$ makes it compact. 

{\bf The fourth term. }
%We abbreviate $\tau_1:=\tau_{z_1}\circ t_{z_2}$. 
We rewrite the fourth term as
\begin{align*}
  T_3v(\tau)
  &= z_1(\tau)\int_{\tau_{z_2}(t_{z_1}(\tau))}^1\frac{z_2^2(\sigma)}{\Delta(\sigma)^2}
  D_1\Bigl((\log z_1)'(\tau_{z_1}(t_{z_2}(\sigma)))\Bigr)(v)d\sigma \cr
  &= z_1(\tau)\int_{\tau_{z_2}(t_{z_1}(\tau))}^1\frac{z_2^2(\sigma)}{\Delta(\sigma)^2}
  \biggl[\left(\frac{v'(\xi)}{z_1(\xi)}-\frac{z_1'(\xi)}{z_1^2(\xi)}v(\xi)
  \right)\bigg|_{\xi=\tau_{z_1}(t_{z_2}(\sigma))} \cr
  &\ \ \ \ \ \  + (\log z_1)''(\tau_{z_1}(t_{z_2}(\sigma)))
  D\tau_{z_1}(v)(t_{z_2}(\sigma))\biggr]d\sigma \cr
  &= \frac{||z_2||^2}{||z_1||^2}z_1(\tau)\int_{\tau}^1
  \frac{z_1^2(\xi)}{\Delta( \tau_{z_2}(t_{z_1}(\xi)))^2}
  \biggl[\left(\frac{v'(\xi)}{z_1(\xi)}-\frac{z_1'(\xi)}{z_1^2(\xi)}v(\xi)\right)\cr
  & \ \ \ \ \ \ + (\log z_1)''(\xi)Xv(\xi)\biggr]d\xi,
\end{align*}
with $Xv(\xi)$ from equation~\eqref{eq:X}. By
Lemma~\ref{lem:Banach-alg} the function
$$
  \xi\mapsto (\log z_1)''(\xi) = \frac{z_1''(\xi)z_1(\xi)-z_1'(\xi)^2}{z_1(\xi)^2}
$$
lies in $L^2$ and depends continuously on $z_1\in H^2$.
All other terms in the integrand are continuous functions of
$\xi$ that depend continuously on $(z_1,z_2)$ by
Lemmas~\ref{lem:H^1cont} and~\ref{lem:Banach-alg} together with
equation~\eqref{eq:X}. So the integrand belongs to 
$L^2(S^1,\R)$, hence the integral belongs to $H^1(S^1,\R)$, thus
$T_4$ defines a bounded operator $H^2\to H^1$ depending continuously on $(z_1,z_2)$, 
and composition with the inclusion $H^1\into L^2$ makes it compact. 

This finishes the discussion of $\mathcal{V}^3_1$, and thus of the
$z_1$-derivative of the $z_1$-component of the $L^2$-gradient
$\nabla\II$. The treatments of the $z_2$-component and 
the $z_2$-derivatives of both components are analogous and will be omitted. 
This concludes the proof of Proposition~\ref{prop:compact}.

%%%%%%%%%%%%%%%%%%%%%%%%%%%%%%%%%%%%%%%%%%%%%%%%%%%%%%%%%%%%%%%%%%%%%%%%%%
\section{Some lemmas on continuous dependence}\label{sec:contdep}
%%%%%%%%%%%%%%%%%%%%%%%%%%%%%%%%%%%%%%%%%%%%%%%%%%%%%%%%%%%%%%%%%%%%%%%%%%

In this appendix we prove some technical lemmas on continuous
dependence that were used in Appendix~\ref{sec:diff}. We will freely use
the notation from Appendix~\ref{sec:diff}. 

%%%
\subsection{The basic lemma on continuous dependence}\label{ss:contbasic}
%%%

We say that a function $f\in C^1(S^1,\R)$ {\it has transverse zeros} if for all $t\in S^1$
with $f(t)=0$ we have $f'(t)\ne 0$. We define the following open subset of $H^2(S^1,\R)$:
$$
  \mathcal{H}_0^2 := \{z\in H^2(S^1,\R)\mid z \text{ has transverse zeros}\}.
$$
Note that by Proposition~\ref{prop:critin} and
Corollary~\ref{cor:inst-transverse-zeroes}, 
for a critical point $(z_1,z_2)$ of the
functional $\BB_{in}$ from Section~\ref{sec:inst-interaction} both components
$z_1,z_2$ belong to $\mathcal{H}_0^2$.
Similarly, we define 
$$
  \C_0^1 := \{f\in C^1(S^1,\R)\mid f \text{ has transverse zeros}\}.
$$
We introduce the maps 
$$
  \F:\mathcal{H}_0^2\longrightarrow \C_0^1,\qquad
  \F(z):=z^3\circ\tau_z,
$$
with $\tau_z$ defined in equation~\eqref{eq:tau}, and
$$
  \G:\C_0^1\longrightarrow L^2(S^1,\R),\qquad
  \G(f):=\frac{1}{f^{1/3}}
$$
(defined outside the zero set of $f$). 
Our goal in this section is to prove the following

\begin{lemma}\label{lem:contin}
The map 
$$
  \G\circ \F: \mathcal{H}_0^2\longrightarrow L^2(S^1,\R),\qquad  
  z\mapsto \frac{1}{z\circ \tau_z}
$$
is continuous.
\end{lemma}

The statement naturally splits in two --- continuity of $\F$ and continuity of $\G$. 

\begin{lemma}\label{lem:continF}
The map $\F$ is continuous.
\end{lemma}

\begin{lemma}\label{lem:continG}
The map $\G$ is continuous.
\end{lemma}

In the proofs we will use the following standard fact whose
easy proof we omit. 

\begin{lemma}\label{lem:homeo}
Denote by ${\rm Homeo}(S^1)$ the space of homeomorphisms of $S^1$
equipped with the $C^0$-topology. Then the inversion $h\mapsto h^{-1}$
defines a continuous map ${\rm Homeo}(S^1)\to {\rm Homeo}(S^1)$.
\hfill$\square$
\end{lemma}

\textbf{Proof of Lemma~\ref{lem:continF}: }
Formula~\eqref{eq:tz-restated} for the homeomorphism $t_z$ shows that it depends continuously
on $z$. Therefore, by Lemma~\ref{lem:homeo}, its inverse $\tau_z$
depends continuously on $z$ as well. This shows that  
$\F$ lands in $C^0(S^1,\R)$ and is continuous as a map to $C^0(S^1,\R)$.
Next we set $f:=\F(z)$ and write out its derivative with respect to $t$:
$$
f'(t)=3z^2(\tau_z(t))z'(\tau_z(t))\dot \tau_z(t)=3z^2(\tau_z(t))z'(\tau_z(t))\frac{||z||^2}{z(\tau_z(t))^2}=
||z||^2z'(\tau_z(t))
$$
Since $z'\in C^0(S^1,\R)$ and $\tau_z$ depend continuously on $z$, we see that $f'$
lands in $C^0(S^1,\R)$ and depends continuously on $z$. Altogether this shows that $\F$
is a continous map to $C^1(S^1,\R)$. 
Transversality of zeros for $z$ and the above formula for $f'(t)$ imply
transversality of zeros for $f$. This completes the proof of Lemma~\ref{lem:continF}.
\hfill $\square$

%%%
\subsection{Proof of Lemma~\ref{lem:continG}}
%%%

{\bf Desingularization procedure. }
The key point in the proof of Lemma~\ref{lem:continG} is the question
of how to deal with integrals of the type 
$$
J=\int_{S^1}\frac{dt}{f^{2/3}(t)}
$$
for a function $f\in \C_0^1$. The idea is to apply a coordinate change that 
turns the integrand into a continuous one. For this consider a $C^1$-homeomorphism 
$\rho:S^1\mapsto S^1$ that restricts to a $C^1$-diffeomorphism 
$S^1\setminus f^{-1}(0)\mapsto S^1\setminus f^{-1}(0)$ and satisfies
$$
  \rho(\tau) = t_i+(\tau-t_i)^3\quad\text{near each zero $t_i$ of $f$}.
$$
%Later, we are going to specify $\delta(\tau_*)>0$ more precisely, see equation~\eqref{eq:quant}.
We perform the coordinate change $t=\rho(\tau)$ in the integral
$J$. Then near a zero $t_i$ of $f$, using $t-t_i=(\tau-t_i)^3$, the
integrand in $J$ becomes
\begin{equation}\label{eq:pullback}
  \frac{dt}{f^{2/3}(t)}
  = \frac{\rho'(\tau)d\tau}{f^{2/3}(t)}
  = \frac{3(\tau-t_i)^2d\tau}{f^{2/3}(t)}
  = 3\left(\frac{(t-t_i)^2}{f^2(t)}\right)^{1/3}\hspace{-12pt}d\tau
  = 3\left(\frac{f(t)}{t-t_i}\right)^{-2/3}\hspace{-15pt}d\tau.
\end{equation}
Since $f$ is a $C^1$-function, the quotient $f(t)/(t-t_i)$ extends
continuously over $t=t_i$ by the derivative $f'(t_i)$, which is
nonzero because $f\in\C_0^1$. Therefore, the function $t\mapsto
3(\frac{f(t)}{t-t_i})^{-2/3}$ extends continuously over
$t=t_i$ by $3f'(t_i)^{-2/3}$. Composing this with the continuous function
$\rho$, we conclude 
%\begin{equation}\label{eq:Lagrange}
%   \frac{f(t)}{t-t_i}=\frac{f(t)-f(t_i)}{t-t_i}=f'(\xi)
%\end{equation}
%for some $\xi$ between $t_i$ and $t$.

\begin{lemma}\label{lem:g}
For $f,\rho$ as above the coefficient in front of $d\tau$ in the
pullback $\rho^*(\frac{dt}{f^{2/3}(t)})$
%Equation~\eqref{eq:pullback}
extends uniquely to a continuous function
\begin{equation}
  g:S^1\to\R,\qquad g(\tau) = \begin{cases}
    \frac{\rho'(\tau)}{f^{2/3}(\rho(\tau))} & f(\tau)\neq 0, \cr
    3f'(t_i)^{-2/3} & \tau = t_i \text{ zero of }f. 
  \end{cases}
\end{equation}
\hfill$\square$
\end{lemma}

We need to globalize this procedure, assigning to each $f\in\C_0^1$ a
map $\rho=\rho_f$ with the properties above in a continuous fashion.
For this, we introduce some notation. 
%{\bf Proof of (C1) and (C2)}
We write $\C_0^1$ as the disjoint union
$$
  \C_0^1 = \coprod_{m\in\N_0} \U_m,
$$
where $\U_m$ is the set of $f\in\C_0^1$ with precisely $m$ zeroes.
Let 
$$
  X_m := \{(t_1,\dots,t_m)\in (S^1)^m\mid t_1<t_2<\cdots<t_m<t_1\}\Bigl/\Z_m
$$
be the configuration space of $m$ cyclically ordered points on $S^1$
modulo cyclic permutations (with the quotient topology). 
Assigning to a function its cyclically ordered zero set defines a
canonical continuous map
$$
  Z:\U_m\longrightarrow X_m.
$$
Let $\G_m$ be the set of $C^1$-homeomorphisms $\rho:S^1\rightarrow
S^1$ with the following properties:
\begin{itemize}
% \item The map $\rho$ is a smooth homeomorphism.
 \item $\rho$ has precisely $m$ critical points (i.e.~zeroes of
   $\rho'$) $t_1,\dots,t_m$.
 \item For each $l=1,\dots,m$ let 
 \begin{equation}\label{eq:quant}
    \delta_i := \frac{1}{4}min\{t_i-t_{i-1},t_{i+1}-t_i\}
 \end{equation}
 be the distance of $t_i$ to the nearest zero, where we stipulate 
 $t_0:=t_m\in S^1$. Then we require that
 $$
    \rho(\tau) = (\tau-t_i)^3+t_i
 $$
 for all $\tau\in (t_i-\delta_i,t_i+\delta_i)$.
\end{itemize}
We equip $\G_m$ with the $C^1$-topology. 
Assigning to a function its cyclically ordered set of critical points defines a
canonical continuous map 
$$
   \pi:\G_m\longrightarrow X_m. 
$$
Since $\pi$ is a fibration with contractible fibres, there exists a
continuous section 
$$
   s:X_m\longrightarrow \G_m.
$$
So we obtain a continuous map
$$
   s\circ Z:\U_m\longrightarrow\G_m,\qquad f\mapsto\rho_f:=s\circ Z(f).
$$
For different $m\in\N_0$ these maps together give a continuous map
\begin{equation}\label{eq:rhof}
   \C_0^1=\coprod_{m\in\N_0}\U_m\longrightarrow\G:=\coprod_{m\in\N_0}\G_m,\qquad f\mapsto\rho_f
\end{equation}
with the property that the critical points of $\rho_f$ are precisely
the zeroes of $f$. 
%$$
%  \rho_f(\tau) = t_i+(\tau-t_i)^3\quad\text{near each zero $t_i$ of $f$}.
%$$
%completes the proof of (C1).
It follows that $f,\rho_f$ satisfy the hypotheses of Lemma~\ref{lem:g};
we denote the resulting continuous function by $g_f:S^1\to\R$.

\begin{lemma}\label{lem:gf}
The map $\C_0^1\to C^0(S^1,\R)$, $f\mapsto g_f$ is continuous. 
\end{lemma}

{\bf Proof: }
We will use the following criterion for uniform convergence of a
sequence of functions $g_n:S^1\to\R$ to a function $g:S^1\to\R$
(which holds more generally for functions on any compact metric space):
\begin{equation}
  g_n\to g \text{ uniformly} \Longleftrightarrow g_n(\tau_n)\to
  g(\tau) \text{ for every sequence } \tau_n\to\tau.
\end{equation}
Consider now a converging sequence $f_n\to f$ in $\C_0^1$ and denote
$g_n:=g_{f_n}$, $g:=g_f$, $\rho_n:=\rho_{f_n}$, $\rho:=\rho_f$. Let
$\tau_n\to\tau$ be a converging sequence in $S^1$. Then by the
criterion we need to show that $g_n(\tau_n)\to g(\tau)$.  
We distinguish two cases.

Case 1: $f(\tau)\neq 0$.\\
Then uniform convergence $f_n\to f$ implies $f_n(\tau_n)\to f(\tau)$,
so $f_n(\tau_n)\neq 0$ for all sufficiently large $n$. Hence by
Lemma~\ref{lem:g} we have
$$
   g_n(\tau_n) = \frac{\rho_n'(\tau_n)}{f_n^{2/3}(\rho_n(\tau_n))}
   \qquad\text{and}\qquad g(\tau) = \frac{\rho'(\tau)}{f^{2/3}(\rho(\tau))}.
$$
Note that $f(\rho(\tau))\neq 0$ and $f_n(\rho_n(\tau_n))\neq 0$ for large $n$. 
Continuity of the map $f\mapsto\rho_f$ implies that $\rho_n\to\rho$ in
$C^1(S^1,S^1)$, hence $\rho_n\to\rho$ and $\rho_n'\to\rho'$ uniformly. 
Applying the above criterion repeatedly it follows that
$\rho_n'(\tau_n)\to \rho'(\tau)$, $\rho_n(\tau_n)\to \rho(\tau)$,
$f_n(\rho_n(\tau_n))\to f(\rho(\tau))\neq 0$, and therefore
$g_n(\tau_n)\to g(\tau)$. 

Case 2: $f(\tau)=0$. \\
In this case let $t_1<\cdots<t_m<t_1$ be the zeroes of $f$. Then for
large $n$ the function $f_n$ also has $m$ zeroes $t_{1,n}<\cdots<t_{m,n}<t_{1,n}$
such that $t_{i,n}\to t_i$ as $n\to\infty$ for each $i$. 
Hence the positive numbers $\delta_i$ and $\delta_{i,n}$ defined via
equation~\eqref{eq:quant} (the latter using the $t_{i,n}$) also satisfy
$\delta_{i,n}\to\delta_i$ as $n\to\infty$. 
We have $\tau=t_i$ for some $i$. Since $\tau_n\to\tau=t_i$ and
$\rho_n\to\rho$, it follows that for large $n$ both $\tau_n$ and
$t_n:=\rho_n(\tau_n)$ lie in the interval $(t_{i,n}-\delta_{i,n},t_{i,n}+\delta_{i,n})$. \\
Let us assume first that $\tau_n\neq t_{i,n}$, hence also $t_n\neq
t_{i,n}$, for all sufficiently large $n$. 
Then $f_n(\tau_n)\neq 0$ for all sufficiently large $n$, and by
Lemma~\ref{lem:g} and equation~\eqref{eq:pullback} we have
$$
   g_n(\tau_n) = \frac{\rho_n'(\tau_n)}{f_n^{2/3}(t_n)} = 3\left(\frac{f_n(t)}{t_n-t_{i,n}}\right)^{-2/3}
   \qquad\text{and}\qquad g(\tau) = 3f'(t_i)^{-2/3}.
$$
By the mean value theorem we have $f_n(t_n)/(t_n-t_{i,n}) = f'(\xi_n)$
for some $\xi_n$ between $t_{i,n}$ and $t_n$. Then $t_{i,n}\to t_i$
and $t_n\to t_i$ implies $\xi_n\to t_i$, so uniform convergence
$f_n'\to f$ yields $f_n'(\xi_n)\to f'(t_i)\neq 0$ and thus
$g_n(\tau_n)\to g(\tau)$. \\
If $\tau_n = t_{i,n}$ for some arbitrarily large $n$, then for these
$n$ by Lemma~\ref{lem:g} we have $g_n(\tau_n)=3f_n'(t_{i,n})^{-2/3}$,
which also converges to $g(\tau) = 3f'(t_i)^{-2/3}$ as $n\to\infty$.
This concludes the proof of Lemma~\ref{lem:gf}.
\hfill$\square$

%We claim that
%\begin{itemize}
%\item [(C1)] there exists a continuous choice of $\rho=\rho_f$ as above for each $f\in \C_0^1$ and
%\item [(C2)] given such a choice, the assignment $f\mapsto \frac{\rho'(\tau)}{f^{2/3}(\rho(\tau))}$ defines a continuous map from $\C_0^1$ to $C^0(S^1,\R)$.
%\end{itemize}
%We finish the proof of Lemma~\eqref{lem:continG} modulo this claim.
%{\bf The proof modulo (C1) and (C2)}

{\bf Proof of Lemma~\ref{lem:continG}. }
Consider a sequence $\{f_n\}_{n\in \N}\subset \C_0^1$ converging to
$f\in \C_0^1$ in $C^1(S^1,\R)$. We need to show $f_n^{-1/3}\to
f^{-1/3}$ in $L^2(S^1,\R)$. We write the squared $L^2$-distance as
\begin{equation}\label{L2-dist}
   ||f^{-1/3}-f_n^{-1/3}||^2_{L^2} = A+A_n-2B_n
\end{equation}
with
$$
   A := \int_0^1\frac{dt}{f^{2/3}(t)},\qquad 
   A_n = \int_0^1\frac{dt}{f_n^{2/3}(t)},\qquad
   B_n = \int_0^1\frac{dt}{f_n^{1/3}(t)f^{1/3}(t)}.
$$
By equation~\eqref{eq:pullback} we have
$$
   A = \int_0^1g_f(\tau)d\tau,\qquad 
   A_n = \int_0^1g_{f_n}(\tau)d\tau
$$
where $g_f,g_{f_n}:S^1\to\R$ are the continuous functions
assigned to $f,f_n$ in Lemma~\ref{lem:g}. From the convergence $f_n\to
f$ in $\C_0^1$ and Lemma~\ref{lem:gf} we obtain $g_{f_n}\to g_f$ in
$C^0(S^1,\R)$ and thus $A_n\to A$. So in view of
equation~\eqref{L2-dist} we are done if we can show $B_n\to A$. 

To prove this we introduce some notation. For $h\in\C_0^1$ and a
measurable subset $I\subset S^1$ we denote
$$
  A(h,I) := \int_{I}\frac{dt}{h^{2/3}(t)}dt.
$$
By equation~\eqref{eq:pullback} we have 
$$
   A(h,I) = \int_{\rho_h^{-1}(I)}g_h(\tau)d\tau 
   \leq |\rho_h^{-1}(I)|\,\|g_h\|_{C^0(S^1,\R)},
$$
where $h\mapsto\rho_h$ and $h\mapsto g_h$ are the continuous maps from 
equation~\eqref{eq:rhof} and Lemma~\ref{lem:gf}, respectively, and
$|I|$ denotes the Lebesgue measure of $I$. It follows from the
definition of $\rho_h$ that $|\rho_h^{-1}(I)|$ depends continuously on
$h\in\C_0^1$ and can be made arbitrarily small by making $|I|$ small.

Let now $\eps>0$ be given and consider the compact subset
$$
   Z := \{f_n\}_{n\in \N}\cup \{f\} \subset \C_0^1.
$$
Since the zero set of $f_n$ converges to that of $f$, there exists for
each $\gamma>0$ a union $I\subset S^1$ of open intervals containing
the zero set of $f$ with the following properties:

\hspace{1em}(i)  there exist $\delta>0$ and $N\in\N$ such that $|f_n|\geq\delta$
  on $S^1\setminus I$ for all $n\geq N$;

\hspace{1em}(ii) $|I|\leq\gamma$. 

By the preceding discussion we can choose $\gamma>0$ so small that
\begin{equation}\label{eq:AhI}
   \max_{h\in Z}A(h,I)
   \leq \max_{h\in Z}|\rho_h^{-1}(I)|\,\max_{h\in Z}\|g_h\|_{C^0(S^1,\R)} 
   < \eps/3. 
\end{equation}
We set
$$
   B_n(I) := \int_I\frac{dt}{f_n^{1/3}(t)f^{1/3}(t)}.
$$
Using the Cauchy-Schwarz inequality and~\eqref{eq:AhI} we estimate
\begin{equation}\label{eq:Bn}
  |B_n(I)|^2
  \leq \int_I\frac{dt}{f_n^{2/3}(t)} \int_{I_k}\frac{dt}{f^{2/3}(t)}
  = A(f_n,I) A(f,I)
  \leq \eps^2/9
\end{equation}
for all $n\ge N$. Now we split $B_n-A$ into summands
$$
  B_n - A
  = B_n(I) - A(f,I) + \int_{S^1\setminus
    I}\left(\frac{1}{f_n^{1/3}(t)f^{1/3}(t)}-\frac{1}{f^{2/3}(t)}\right)dt
$$
and estimate it by the triangle inequality:
$$
  |B_n-A|
  \leq |B_n(I)|+ |A(f,I)| + \int_{S^1\setminus
    I}\left|\frac{1}{f_n^{1/3}(t)f^{1/3}(t)}-\frac{1}{f^{2/3}(t)}\right|dt.
$$
By property (i) above we have $|f_n|\geq\delta$ and $|f|\geq\delta>0$
on the compact set $S^1\setminus I$, so the integrand in the last
integral converges uniformly to zero on $S^1\setminus I$. Therefore,
there exists an integer $N_1\geq N$ such that the last integral is
smaller than $\eps/3$ for all $n\geq N_1$. Together with
equations~\eqref{eq:AhI} and~\eqref{eq:Bn} this implies
$$
   |B_n-A| < \eps/3 + \eps/3 + \eps/3 = \eps
$$
for all $n\geq N_1$. This proves $B_n\to A$, which concludes the proof
of Lemma~\ref{lem:continG} and thus of Lemma~\ref{lem:contin}. 
\hfill$\square$

%%%
\subsection{Further lemmas on continuous dependence.}
%%%

We will frequently use the following standard result for 
which we denote 
$$
H_{\ne 0}^1(S^1,\R):=\{f\in H^1(S^1,\R)\mid f(t)\neq 0
\text{ for all } t\}
$$
\begin{lemma}\label{lem:Banach-alg}
There are continuous maps
$$
   H^1(S^1,\R)\times H^1(S^1,\R)\to H^1(S^1,\R),\qquad (f,g)\mapsto fg
$$
and
$$
H_{\ne 0}^1(S^1,\R)\to H^1(S^1,\R),\qquad f\mapsto 1/f. 
$$
\end{lemma}

\textbf{Proof: } 
The first assertion is just the well-known fact that 
$H^1(S^1,\R)$ is
a Banach algebra. For the second assertion we abbreviate
$H_{\ne 0}^1:=H_{\ne 0}^1(S^1,\R)$ etc. 
The map $H_{\ne 0}^1\to L^2$, $f\mapsto 1/f$ is the
composition of continuous maps
$$
   H_{\ne 0}^1\into C_{\ne 0}^0 \to C_{\ne 0}^0 \into L^2,
$$
where the middle map sends $f\mapsto 1/f$ and the other two maps
are the canonical inclusions. The map $H_{\ne 0}^1\to L^2$,
$f\mapsto (1/f)'=-f'/f^2$ 
is the composition of continuous maps
$$
   H_{\ne 0}^1\to C_{\ne 0}^0\oplus L^2 \to L^2,
$$
where the first map sends $f\mapsto (-1/f^2,f')$ and the second
map is multiplication. Together this proves the lemma. 
\hfill $\square$

\begin{lemma}\label{lem:H^1cont}
The following maps are continuous:
\begin{gather*}
   \mathcal{H}_0^2\longrightarrow H^1(S^1,\R),\qquad
   z\mapsto z^2\circ \tau_z, \cr 
   (H^2(S^1,\R)\times\mathcal{H}_0^2)\cap
   \mathcal{H}_{in}\longrightarrow H^1(S^1,\R), \qquad
   (z_1,z_2)\mapsto \left(\tau\mapsto z^2_2(\tau_{z_2}(t_{z_1}(\tau)))\right), \cr
   (H^2(S^1,\R)\times\mathcal{H}_0^2)\cap
   \mathcal{H}_{in}\longrightarrow H^1(S^1,\R), \qquad
   (z_1,z_2)\mapsto\left(\tau\mapsto\Delta(\tau_{z_2}(t_{z_1}(\tau)))\right),
\end{gather*}
where $\Delta$ is defined in~\eqref{eq:Delta}.
\end{lemma}

\textbf{Proof: }
For the first map continuous dependence of $\tau_z$ on $z$ (which
follows from Lemma~\ref{lem:homeo}) implies that 
$z^2\circ \tau_z$ depends continuously on $z$ as an element of $L^2(S^1,\R)$.
For its time derivative
$$
\frac{d}{dt}z^2(\tau_z(t))=2z(\tau_z(t))\dot\tau_z(t)=2z(\tau_z(t))\frac{||z||^2}{z(\tau_z(t))^2}=
2\frac{||z||^2}{z(\tau_z(t))}
$$
continuous dependence follows from Lemma~\ref{lem:contin}.

For the second map observe that $z_1$ never vanishes. Therefore,
$t_{z_1}$ is a $C^1$-diffeomorphism of $S^1$ depending continuously on
$z_1$. Together with the statement about the first map this
concludes the argument. Continuity of the third map now follows
directly from equation~\eqref{eq:Delta-transf}. 
\hfill $\square$

\begin{lemma}\label{lem:z_1z_2L^2cont}
The map 
$$
(H^2(S^1,\R)\times\mathcal{H}_0^2)\cap \mathcal{H}_{in}\longrightarrow L^2(S^1,\R)
$$
defined by
$$
(z_1,z_2)\mapsto \left(\tau\mapsto\frac{1}{z_2(\tau_{z_2}(t_{z_1}(\tau)))}\right)
$$
is continuous.
\end{lemma}

\textbf{Proof: }
Recall from the previous proof that $t_{z_1}$ is a $C^1$-diffeomorphism 
of $S^1$ depending continuously on $z_1$. This together with
Lemma~\ref{lem:contin} completes the proof. 
\hfill $\square$

%%%%%%%%%%%%%%%%%%%%%%%%%%%%%%%%%%%%%%%%%%%%%%%%%%%%%%%%%%%%%%%%%%%%%%%%%%
\section{The mod $2$ Euler number}\label{sec:Euler}
%%%%%%%%%%%%%%%%%%%%%%%%%%%%%%%%%%%%%%%%%%%%%%%%%%%%%%%%%%%%%%%%%%%%%%%%%%

Throughout this appendix $X$ denotes a Hilbert manifold (an open
subset of a Hilbert space will suffice for our purposes), 
$Y$ an open neighbourhood of $0$ in a Hilbert space,
and $k\in\N_0$.
We say that a $C^1$-map $f:X\to Y$ is {\em transverse to $0$} if $0$
is a regular value of $f$, i.e.~$Df(x):T_xX\to T_0Y$ is surjective for
all $x\in f^{-1}(0)$.  
Our goal is to prove

\begin{thm}\label{thm:zero-Fredholm}
To each $C^1$-Fredholm map $f:X\to Y$ of index $0$ with compact zero
set $f^{-1}(0)$ we can associate its {\em mod $2$ Euler number} $\chi(f)\in\Z/2\Z$
which is uniquely characterized by the following axioms:

(Transversality) 
If $f$ is transverse to $0$, then $\chi(f)=|f^{-1}(0)|$ mod $2$. 

(Excision)
For any open neighbourhoods $\wt X\subset X$ of $f^{-1}(0)$ and $\wt
Y\subset Y$ of $0$ such that $f(\wt X)\subset\wt Y$ we have $\chi(f) =
\chi(f|_{\wt X}:\wt X\to\wt Y)$. 

(Cobordism)
If $W$ is a Hilbert manifold with boundary and $F:W\to Y$ a
$C^1$-Fredholm map of index $1$ with compact zero set $F^{-1}(0)$,
then $\chi(F|_{\p W}:\p W\to Y)=0$. 

(Homotopy)
Let $f_0,f_1:X\to Y$ be $C^1$-Fredholm maps of index $0$ with compact zero
sets. If there exists a $C^1$-Fredholm map $F:[0,1]\times X\to Y$ of
index $1$ with compact zero set $F^{-1}(0)$ such that $F|_{\{i\}\times
  X}=f_i$ for $i=0,1$, then $\chi(f_0)=\chi(f_1)$.
\end{thm}

Note that the (Homotopy) axiom is just a special case of the (Cobordism) axiom.

The invariant $\chi(f)$ in Theorem~\ref{thm:zero-Fredholm} can be
viewed as a special case of either the {\em Caccioppoli--Smale degree}
defined in~\cite{smale}, or of the {\em Euler class of a $G$-moduli
  problem} defined in~\cite{cieliebak-mundet-salamon} (with trivial group $G$). 
The main improvement of Theorem~\ref{thm:zero-Fredholm} over those
results is the fact that in the (Cobordism) we require only regularity
$C^1$ instead of $C^2$. This is not entirely obvious because the
Sard--Smale theorem~\cite{smale} for a Fredholm map of index $1$ requires
regularity $C^2$. While this improvement may seem boring from the viewpoint
of general theory, it is crucial for the application in this paper
because the functional $\nabla\BB_{in}$ is of class $C^1$ but not $C^2$. 

The main idea of the proof of Theorem~\ref{thm:zero-Fredholm}
is the following. Given a $C^1$-Fredholm map $f:X\to Y$
of index $0$ with compact zero set $f^{-1}(0)$ 
we use the Sard-Smale theorem~\cite{smale} to $C^1$-approximate  
$f$ by a $C^1$-Fredholm map $g:X\to Y$
of index $0$ transverse to $0$. We define 
$\chi(f):=\chi(g):=|g^{-1}(0)|$.
To see that this is well-defined we join $g_0$ to $g_1$
by a convex linear combination $G:[0,1]\times X\to Y$.
The map $G$ is a $C^1$-Fredholm map of index $1$. We want 
to approximate $G$ by a $C^1$-Fredholm map $\tilde G$
transverse to $0$ and such that $\tilde G(i,x)=G(i,x)$  
for $i=0,1$. The direct application of the 
Sard-Smale theorem~\cite{smale} is not sufficient, since 
Fredholm maps of index $1$ produce a loss of regularity by $1$. On the
other hand, Lemma~\ref{lem:Fredholm-approx} 
below applies and we conclude by the standard argument that
$M:=\tilde G^{-1}(0)$ is a compact $C^1$-manifold with 
boundary $\p M=g_0^{-1}(0)\amalg g_1^{-1}(0)$. Therefore,
$\chi(g_0)=|g_0^{-1}(0)|=|g_1^{-1}(0)|=\chi(g_1)$.

The proof of Lemma~\ref{lem:Fredholm-approx} below uses
a three-step approximation, where the key application of the Sard-Smale theorem~\cite{smale} occurs on step 3) 
to a $C^\infty$-map defined on an open subset of the domain, so the loss of regularity does not happen. 
To get to this nice situation of
step 3) we use the fact that the restriction of $F$ 
to $\p W$ (notation of Lemma~\ref{lem:Fredholm-approx})
has index zero and can therefore be approximated by 
the Sard-Smale theorem~\cite{smale} without the loss of regularity.

%The proof of Theorem~\ref{thm:zero-Fredholm} uses the following lemma.

\begin{lemma}\label{lem:Fredholm-approx}
Let $W$ be a Hilbert manifold with boundary and $F:W\to Y$ a
$C^k$-Fredholm map of index $k\in\N$ with compact zero set $F^{-1}(0)$
such that $F|_{\p W}$ is transverse to $0$. 
Then $F$ can be $C^k$-approximated by a $C^k$-Fredholm map $\wt F:W\to Y$ of index $k$
with compact zero set $\wt F^{-1}(0)$ such that $F|_{\p W}=\wt F|_{\p W}$
and $\wt F:W\to Y$ is transverse to $0$. 
In particular, $\wt F^{-1}(0)$ is a compact $C^k$-manifold of
dimension $k$ with boundary $\p\wt F^{-1}(0)=(F|_{\p W})^{-1}(0)$. 
\end{lemma}

{\bf Proof: }
1) Let $g:\p W\to Y$ be a $C^\infty$-map which is sufficiently
$C^k$-close to $f:=F|_{\p W}$ so that the map $(1-t)f+tg$ is
transverse to zero for each $t\in[0,1]$ and the map
$$
  [0,1]\times\p W\to Y,\qquad (t,x)\mapsto (1-t)f(x)+tg(x)
$$
is Fredholm of index $k$ with compact zero set. 
Let $\NN\cong[0,2)\times\p W\subset W$ be a collar neighbourhood of
$\p W\cong \{0\}\times\p W$. 
%and set $f_s:=f|_{\{s\}\times\p W}$ for $s\in[0,2)$. 
Pick a smooth cutoff
function $\varphi:\R\to[0,1]$ with compact support in $(0,2)$ which
equals $1$ in a neighbourhood of $1$. Define $G:W\to Y$ by
$$
   G(t,x) := (1-\varphi(t))F(t,x) + \varphi(t)g(x)
$$
for $(t,x)\in\NN$ and $G:=F$ on $W\setminus\NN$. Then $G:W\to Y$ is a
$C^k$-Fredholm map of index $k$ such that $G|_{\p W}=F|_{\p W}$ and
$G$ is $C^\infty$ in a neighbourhood of $\{1\}\times\p W$. 
By choosing $\NN$ sufficiently small we can ensure that $G$ has
compact zero set and $G|_\NN$ is transverse to $0$. 

2) Set $\wt\NN:=[0,1)\times\p W\subset\NN$. We $C^k$-approximate $G$ by
a $C^k$-Fredholm map $H:W\to Y$ of index $k$ with compact zero set
such that $H=G$ on $\wt\NN$ and $H$ is $C^\infty$ on $W\setminus\wt\NN$.

3) Using the Sard--Smale theorem~\cite{smale}, we $C^k$-approximate the restriction
$H|_{W\setminus\wt\NN}$ by a $C^\infty$-Fredholm map $\wt
F:W\setminus\wt\NN\to Y$ with compact zero set which agrees with $H$
near $\{1\}\times\p W$ and is transverse to zero. This map extends via
$H$ over $\wt\NN$ to the desired map $\wt F:W\to Y$. 
\hfill$\square$
\medskip

{\bf Proof of Theorem~\ref{thm:zero-Fredholm}: }

\underline{Uniqueness}:
Let $f:X\to Y$ be a $C^1$-Fredholm map of index $0$ with compact zero
set $f^{-1}(0)$. By the Sard--Smale theorem~\cite{smale}, we can $C^1$-approximate 
$f$ by a $C^1$-Fredholm map $g:X\to Y$ of index $0$ which is
transverse to $0$. By choosing $g$ sufficiently $C^1$-close to $f$ we
can ensure that
$$
  F:[0,1]\times X\to Y,\qquad (t,x)\mapsto (1-t)f(x)+tg(x)
$$
is a $C^1$-Fredholm map of index $1$ with compact zero set. So by the
(Homotopy) and (Transversality) axioms $\chi(f)$ is uniquely
determined by
$$
   \chi(f) = \chi(g) = |g^{-1}(0)| \mod 2.
$$
\underline{Existence}:
On maps $f$ as in the theorem which are transverse to $0$ we define
$\chi(f):=|f^{-1}(0)|$ mod $2$ by the (Transversality) axiom. We claim
that then the (Cobordism) axiom holds under the additional assumption
that $f:=F|_{\p W}$ is transverse to $0$. To see this, we apply
Lemma~\ref{lem:Fredholm-approx} to find a $C^1$-Fredholm map $\wt
F:W\to Y$ of index $1$ whose zero set $M:=\wt F^{-1}(0)$ is a compact
$C^1$-manifold of dimension $1$ with boundary $\p M=f^{-1}(0)$. Since
$M$ has an even number of boundary points, we conclude $|f^{-1}(0)|=0$
mod $2$ and the claim is proved. 
As a special case, the (Homotopy) axiom also holds under the
additional assumption that $f_0$ and $f_1$ are transverse to $0$.

Let now $f:X\to Y$ be any $C^1$-Fredholm map of index $0$ with compact zero
set $f^{-1}(0)$. We choose $g,F$ as in the proof of uniqueness and
define $\chi(f):=\chi(g)=|g^{-1}(0)|$. To see that this is well-defined, let
$g_i,F_i$, $i=0,1$ be two such choices. Then the maps $F_0,F_1$ can be
joined by a cutoff construction to a $C^1$-Fredholm map $G:[0,1]\times X\to Y$ of
index $1$ with compact zero set such that $G|_{\{i\}\times
  X}=g_i$ for $i=0,1$, and the special case of the (Homotopy) axiom
implies $\chi(f_0)=\chi(f_1)$. So $\chi(f)$ is well-defined. By
construction, it satisfies the (Transversality) and (Excision) axioms. 
The (Cobordism) axiom follows from the special case above, and the
(Homotopy) axion is a special case of this. 
\hfill$\square$

%%%%%%%%%%%%%%%%%%%%%%%%%%%%%%%%%%%%%%%%%%%%%%%%%%%%%%%%%%%%%%%%%%%%%%%%%%
\section{The mod $2$ Euler number of the mean interaction
  functional}\label{sec:Euler-mean-interaction} 
%%%%%%%%%%%%%%%%%%%%%%%%%%%%%%%%%%%%%%%%%%%%%%%%%%%%%%%%%%%%%%%%%%%%%%%%%%

In this section we prove the following result which is used in the
proof of the Existence Theorem~\ref{thm:existence}, but which may also
be of independent interest. 

\begin{thm}\label{thm:Euler-mean-intersection}
The $L^2$-gradient of the mean interaction functional on symmetric loops,
$$
  \nabla\BB_{av}:\wh\HH^2_{in}\to H_{sym}^0(S^1,\R^2), 
$$
is a $C^1$-Fredholm map whose mod $2$ Euler number equals $1$. 
\end{thm}

It follows from the main result in~\cite{frauenfelder} that $\nabla\BB_{av}$
has a unique zero $(z_1,z_2)$. To conclude that its mod $2$ Euler number
equals $1$, we would need to prove invertibility of the Hessian
$D\nabla\BB_{av}(z_1,z_2)$, or equivalently (since the Fredholm index
is zero) triviality of its kernel. This is still complicated because
the Hessian has many terms and the equations for an element $(v_1,v_2)$
in its kernel are coupled. Therefore, we will instead further deform
the mean interaction functional to one for which the equations
decouple, and then compute the mod $2$ Euler number of the latter.

%%%
\subsection{Decoupling the mean interaction}\label{ss:decoupling}
%%%

In this subsection we describe the deformation of the mean interaction
functional on symmetric loops to one for which the equations decouple.
We will first phrase it in terms of the original (physical) coordinates $q(t)$;
the Levi-Civita transformation to the new coordinates $z(\tau)$ will
be considered in the next subsection. 

In the following discussion we adapt some arguments from~\cite{frauenfelder2}
to our situation. In order to be consistent with the notation in that paper,
we replace the period $1$ used in this paper before
by period $2$ as in \cite{frauenfelder2}, 
and consider the restriction of a symmetric loop to half a period. 
Thus we consider
$$
   q_1\in C^\infty\bigl([0,1],(0,\infty)\bigr),\qquad
   q_2\in C^0\bigl([0,1],[0,\infty)\bigr)\cap C^\infty\bigl([0,1),(0,\infty)\bigr)
$$
satisfying
\begin{equation}\label{eq:q-symmetric}
   \dot q_1(0)=\dot q_1(1)=\dot q_2(0)=q_2(1) = 0\quad\text{and}\quad
   q_1(t)>q_2(t)\geq 0\;\forall t\in[0,1].
\end{equation}
We assume that $(q_1,q_2)$ solves a system of second order ODEs
coupled through their means
\begin{equation}\label{eq:mean-interaction-f}
\left\{\;
\begin{aligned}
  \ddot{q}_1(t) &= -\frac{2}{q_1(t)^2}+f_1(\ol q_1,\ol q_2), \cr
  \ddot{q}_2(t) &= -\frac{2}{q_2(t)^2}-f_2(\ol q_1,\ol q_2), 
\end{aligned}
\right.
\end{equation}
where $f_i(\ol q_1,\ol q_2)$ are continuous functions defined for $\ol
q_1>\ol q_2\geq 0$ satisfying 
\begin{equation}\label{eq:f-positive}
   f_1(\ol q_1,\ol q_2)>0\quad\text{and}\quad f_2(\ol q_1,\ol q_2)\geq
   0\quad\text{ for all }\ol q_1>\ol q_2\geq 0.
\end{equation}

\begin{lemma}\label{lem:q-estimates}
Under the above assumptions the following hold.

(a) The map $q_1$ is constant, $q_1(t)\equiv\ol q_1$, where $\ol
q_1>0$ solves the equation
\begin{equation}\label{eq:fixed-point}
   \ol q_1^2f_1(\ol q_1,\ol q_2) = 2.
\end{equation}
(b) The map $q_2$ is strictly concave and strictly decreasing with
maximum $q_2^{\rm max}=q_2(0)$ satisfying the estimates
\begin{equation}\label{eq:q2-estimates}
   \frac12\leq\frac{q_2^{\rm max}}{2}\leq\ol q_2\leq q_2^{\rm max}\leq
   2+\frac{f_2(\ol q_1,\ol q_2)}{2}.
\end{equation}
\end{lemma}

{\bf Proof: }
We abbreviate $f_1=f_1(\ol q_1,\ol q_2)$.

(a) Note first that, since $q_1$ solves the first equation
in~\eqref{eq:mean-interaction-f}, it is actually smooth and $1$-periodic.
It attains its maximum at some time $t_{\rm
  max}\in[0,1]$ satisfying $\dot q(t_{\rm max})=0$ and 
$$
   \ddot q_1(t_{\rm max}) = -\frac{2}{q_1(t_{\rm max})^2} + f_1 \leq 0.
$$
It follows that for all $t\in[0,1]$ we have
$$
  \ddot q_1(t) = -\frac{2}{q_1(t)^2} + f_1
  \leq -\frac{2}{q_1(t_{\rm max})^2} + f_1 \leq 0.
$$
By periodicity this implies $\ddot q_1\equiv 0$, so $\dot q_1$ is
constant. Again by periodicity this implies $\dot q_1\equiv 0$, so
$q_1\equiv\ol q_1$ is constant and the first equation
in~\eqref{eq:mean-interaction-f} becomes equation~\eqref{eq:fixed-point}.
%\smallskip

(b) Since $f_2\geq 0$, the second equation
in~\eqref{eq:mean-interaction-f} implies $\ddot q_2(t)<0$ for all
$t\in[0,1)$, so $q_2$ is strictly concave. Together with $\dot
q_2(0)=0$ this implies $\dot q_2(t)<0$ for all $t\in(0,1)$, so $q_2$
is strictly decreasing with maximum $q_2^{\rm max}=q_2(0)>0$.

For the first inequality in~\eqref{eq:q2-estimates}, we use the second
equation in~\eqref{eq:mean-interaction-f} and $f_2\geq 0$ to estimate
for all $t\in[0,1)$:
\begin{align*}
  \dot q_2(t)
  &= \int_0^t\ddot q_2(s)ds 
  \leq -\int_0^t\frac{2}{q_2(s)^2}ds 
  \leq -\int_0^t\frac{2}{(q_2^{\rm max})^2}ds 
  = -\frac{2}{(q_2^{\rm max})^2}t.
\end{align*}
Together with $q_2(1)=0$ this implies
\begin{align*}
  -q_2^{\rm max}
  &= q_2(1)-q_2(0)
  = \int_0^1\dot q_2(t)dt 
  \leq -\frac{2}{(q_2^{\rm max})^2}\int_0^1tdt
  =-\frac{1}{(q_2^{\rm max})^2},
\end{align*}
hence $(q_2^{\rm max})^3\geq 1$ and thus $q_2^{\rm max}\geq 1$. 

The second inequality in~\eqref{eq:q2-estimates} follows from
concavity of $q_2$: For all $t\in[0,1]$ we have
$$
   q_2(t) \geq (1-t)q_2(0)+tq_2(1) = (1-t)q_2^{\rm max},
$$
and it follows that
$$
   \ol q_2 = \int_0^1q_2(t)dt \geq q_2^{\rm max}\int_0^1(1-t)dt =
   \frac{q_2^{\rm max}}{2}. 
$$
The third inequality in~\eqref{eq:q2-estimates} is clear, so it
remains to prove the fourth one. Since $q_2$ is strictly decreasing
and $q_2(0)\geq 1$, there exists a unique $t_0\in[0,1)$ with $q_2(t_0)=1$.
Then for all $t\in[0,t_0]$ we have $q_2(t)\geq 1$, and therefore
\begin{align*}
  \ddot q_2(t) &= -\frac{2}{q_2(t)^2} - f_2 \geq -2-f_2,\cr
  \dot q_2(t) &= \int_0^t\ddot q_2(s)ds \geq -(2+f_2)t
\end{align*}
for all $t\in[0,t_0]$. This implies
$$
  1-q_2^{\rm max} = q_2(t_0)-q_2(0) = \int_0^{t_0}\dot q_2(t)dt \geq
  -(2+f_2)\int_0^{t_0}tdt = -\frac{t_0^2}{2}(2+f_2), 
$$
and with $t_0\leq 1$ we obtain
$$
   q_2^{\max} \leq 1 + \frac{t_0^2}{2}(2+f_2) \leq 1 +
   \frac{1}{2}(2+f_2) = 2 + \frac{f_2}{2}.
$$
\hfill$\square$

{\bf Compactness. }
We wish to consider families of problems~\eqref{eq:mean-interaction-f}
as above parametrized by pairs of functions $(f_1,f_2)$ satisfying~\eqref{eq:f-positive}.
For compactness of the corresponding space of solutions $(q_1,q_2)$ we need
\begin{itemize}
\item a uniform lower bound $\ol q_1-q_2^{max}\geq\delta>0$, and 
\item a uniform upper bound $\ol q_1\leq c<\infty$,
\end{itemize}
where ``uniform'' means independent of the parameter. In view of
Lemma~\ref{lem:q-estimates}, this is ensured by the following
{\em sufficient condition: For all solutions $(\ol q_1,\ol q_2)$ of 
\begin{equation}\label{eq:q-input}
  \ol q_1^2f_1(\ol q_1,\ol q_2) = 2,\qquad
  \frac12\leq\ol q_2\leq 2+\frac{f_2(\ol q_1,\ol q_2)}{2},\qquad
  \ol q_1>\ol q_2
\end{equation}
we have uniform lower and upper bounds
\begin{equation}\label{eq:q-output}
   \ol q_1-2\ol q_2\geq\delta>0,\quad\text{and}\quad \ol q_1\leq c<\infty.
\end{equation}
}
Indeed, by Lemma~\ref{lem:q-estimates} the averages $(\ol q_1,\ol
q_2)$ of a solution $(q_1,q_2)$ of
problem~\eqref{eq:mean-interaction-f} satisfy conditions~\eqref{eq:q-input},
and in view of $q_2^{\rm max}\leq2\ol q_2$ the first inequality in~\eqref{eq:q-output}
implies the uniform lower bound $\ol q_1-q_2^{max}\geq\delta>0$. The
following lemma describes a situation where this sufficient condition
is satisfied. 

\begin{lemma}\label{lem:remove-f2}
Suppose that
$$
   f_1(\ol q_1,\ol q_2) = \frac{1}{(\ol q_1-\ol q_2)^2}\quad\text{and}\quad
   0\leq f_2(\ol q_1,\ol q_2) \leq \frac{1}{(\ol q_1-\ol q_2)^2}.
$$
Then each solution $(\ol q_1,\ol q_2)$ of equation~\eqref{eq:q-input}
satisfies
$$
   \ol q_1 = (2+\sqrt{2})\ol q_2
$$
as well as the lower and upper bounds
\begin{equation*}%\label{eq:q-output}
  \ol q_1-2\ol q_2\geq\frac{1}{\sqrt{2}},\quad\text{and}\quad
  \ol q_1\leq(2+\sqrt{2})\Bigl(2+\frac{2}{(1+\sqrt{2})^2}\Bigr).
\end{equation*}
\end{lemma}

{\bf Proof: }
In this case the first equation in~\eqref{eq:q-input} becomes the
homogeneous quadratic equation
$$
  \ol q_1^2 = \frac{2}{f_1(\ol q_1,\ol q_2)} = 2(\ol q_1-\ol q_2)^2, 
$$
which has the solutions $\ol q_1 = (2\pm\sqrt{2})\ol q_2$. The
condition $\ol q_1>\ol q_2$ enforces $\ol q_1 = (2+\sqrt{2})\ol q_2$.
Together with $\ol q_2\geq 1/2$ this implies the lower bound
$$
   \ol q_1-2\ol q_2 = \sqrt{2}\,\ol q_2 \geq \frac{\sqrt{2}}{2}
$$
as well as
$$
   \ol q_1-\ol q_2 = (1+\sqrt{2})\ol q_2 \geq \frac{1+\sqrt{2}}{2}. 
$$
With the condition on $f_2$ this yields an upper bound on $\ol q_2$,
$$
  \ol q_2
  \leq 2+\frac{f_2(\ol q_1,\ol q_2)}{2}
  \leq 2+\frac{1}{2(\ol q_1-\ol q_2)^2}
  \leq 2 + \frac{2}{(1+\sqrt{2})^2},
$$
and thus on $\ol q_1$,
$$
  \ol q_1 = (2+\sqrt{2})\ol q_2
  \leq (2+\sqrt{2})\Bigl(2 + \frac{2}{(1+\sqrt{2})^2}\Bigr).
$$
\hfill$\square$

Fixing $f_1(\ol q_1,\ol q_2) = \frac{1}{(\ol q_1-\ol q_2)^2}$,
Lemma~\ref{lem:remove-f2} allows us to linearly interpolate between 
$f_2(\ol q_1,\ol q_2) = \frac{1}{(\ol q_1-\ol q_2)^2}$ and $f_2=0$.
We are thus led to consider the {\em decoupled mean interaction problem}
\begin{equation}\label{eq:mean-interaction-decoupled}
\left\{\;
\begin{aligned}
  \ddot{q}_1(t) &= -\frac{2}{q_1(t)^2}+\frac{1}{(\ol q_1-\ol q_2)^2}, \cr
  \ddot{q}_2(t) &= -\frac{2}{q_2(t)^2}. 
\end{aligned}
\right.
\end{equation}
Note that the second equation is a pure Kepler problem which is
not coupled to the first one. It has a unique solution
$q_2:[0,1]\to[0,\infty)$ with $\dot q_2(0)=q_2(1) = 0$ 
and we denote by
$$
   a := \ol q_2
$$
its average. Note that by Lemma~\ref{lem:q-estimates} it satisfies
$$
   1/2 \leq a \leq 2.
$$
Inserting $\ol q_2=a$ into the first equation,
Lemma~\ref{lem:remove-f2} shows that it has a unique solution $q_1$,
which is constant and given by
$$
  q_t(t)\equiv\ol q_1 = (2+\sqrt{2})a.
$$
This concludes our discussion of compactness. 
In the next subsection we will consider its Levi-Civita transformation
and use it to prove Theorem~\ref{thm:Euler-mean-intersection}. 

%%%
\subsection{A Fredholm homotopy}\label{ss:Fredholm-homotopy}
%%%

%{\bf Levi-Civita transformation. }
%Under the Levi-Civita transformation $q(t)=z(\tau(t))^2$ the decoupled
%mean interaction problem~\eqref{eq:mean-interaction-decoupled} becomes

Let $X$ and $Y$ be as in Section~\ref{ss:existence-proof}. 
For $r\in[0,1]$ we consider the map
$$
   F^r = (F_1,F_2^r):X\to Y
$$
given by
\begin{equation*}%\label{eq:z-av}
%\left\{
\begin{aligned}
  F_1(z_1,z_2) &:= -z_1'' + a_1(z_1,z_2) z_1 + b_1(z_1,z_2)z_1^3, \\
  F_2^r(z_1,z_2) &:= -z_2'' + a_2^r(z_1,z_2) z_2 + b_2^r(z_1,z_2)z_2^3
\end{aligned}
%\right.
\end{equation*}
with the functions
\begin{eqnarray*}
  a_1(z_1,z_2)
  &=& \frac{||z_1'||^2}{||z_1||^2}-\frac{1}{||z_1||^6}-\frac{||z_2||^4
    \cdot ||z_1^2||^2}{2||z_1||^2 \cdot\big(||z_1^2||^2\cdot ||z_2||^2 
  -||z_2^2||^2 \cdot ||z_1||^2\big)^2},\\
  b_1(z_1,z_2)
  &=& +\frac{||z_2||^4}{\big(||z_1^2||^2\cdot ||z_2||^2-||z_2^2||^2 \cdot ||z_1||^2\big)^2},\\
  a_2^r(z_1,z_2)
  &:=& \frac{||z_2'||^2}{||z_2||^2}-\frac{1}{||z_2||^6} +
  r\frac{||z_1||^4 \cdot ||z_2^2||^2}{2||z_2||^2
    \cdot\big(||z_1^2||^2\cdot ||z_2||^2 
  -||z_2^2||^2 \cdot ||z_1||^2\big)^2},\\
  b_2^r(z_1,z_2)
  &:=& -r\frac{||z_1||^4}{\big(||z_1^2||^2\cdot ||z_2||^2-||z_2^2||^2 \cdot ||z_1||^2\big)^2}.
\end{eqnarray*}
For $r=1$ comparison with equations~\eqref{eq:gradQ}
and~\eqref{eq:gradA} shows that
$$
   \nabla_1\BB_{av}(z_1,z_2) = 4\|z_1\|^2 F_1(z_1,z_2),\qquad
   \nabla_2\BB_{av}(z_1,z_2) = 4\|z_2\|^2 F_2^1(z_1,z_2).
$$
Thus, up to the irrelevant positive factors $4\|z_i\|^2$, $F^1$
agrees with $\nabla\BB_{av}$. In particular, the zeroes of $F^1$
satisfy the coupled ODEs~\eqref{eq:z-av}. 

For $r=0$ the first component remains unchanged and comparison with
equation~\eqref{eq:gradQ} shows that 
\begin{equation}\label{eq:Q-F2}
   \nabla\QQ(z_2) = 4\|z_2\|^2 F_2^0(z_1,z_2).
\end{equation}
So the second component of $F^0$ is decoupled from the first one and
corresponds to a pure Kepler problem. 

The discussion in Section~\ref{sec:mean-interaction} shows that
under the Levi-Civita transformations $q_i(t)=z_i(\tau{z_i}(t))^2$ zeroes of
$F^r$ correspond to generalized solutions of the coupled ODEs
\begin{equation}\label{eq:mean-interaction-r}
\left\{\;
\begin{aligned}
  \ddot{q}_1(t) &= -\frac{2}{q_1(t)^2}+\frac{1}{(\overline{q}_1-\overline{q}_2)^2}, \cr
  \ddot{q}_2(t) &= -\frac{2}{q_2(t)^2}-\frac{r}{(\overline{q}_1-\overline{q}_2)^2}.
\end{aligned}
\right.
\end{equation}
By definition of the space $X$, the $q_i$ are symmetric and therefore, after
replacing their period $1$ by $2$, they satisfy conditions~\eqref{eq:q-symmetric}
in the previous subsection. By the discussion in that subsection, the
$q_i$ satisfy a lower bound $q_1(t)\geq q_2(t)\geq\delta>0$
and an upper bound $q_1(t)\leq c<\infty$, uniform in $r\in[0,1]$. This
implies that the zero set of the $C^1$-Fredholm homotopy
$$
   F:[0,1]\times X\to Y,\qquad (r,z_1,z_2)\mapsto F^r(z_1,z_2)
$$
is compact, so by the (Homotopy) axiom in
Theorem~\ref{thm:zero-Fredholm} the mod $2$ Euler numbers satisfy
$$
   \chi(\nabla\BB_{av}) = \chi(F^1) = \chi(F^0). 
$$
To prove Theorem~\ref{thm:Euler-mean-intersection}, it thus remains to
compute $\chi(F^0)$. By the discussion in the previous subsection,
$F^0$ has a unique zero $(z_1,z_2)$ whose components correspond under
the Levi-Civita transformation $q_i(t)=z_i(\tau{z_i}(t))^2$ to the
unique solution $q_2$ of the pure Kepler problem and the constant
solution $q_1(t)\equiv\ol q_2=(2+\sqrt{2})a$, where $a=\ol q_2$. In
particular, the first component is constant and given by
\begin{equation}\label{eq:z1bar}
  z_1(t)\equiv\ol z_1 = \sqrt{(2+\sqrt{2})a}.
\end{equation}
It thus remains to prove that the derivative $DF^0(z_1,z_2)$ at its unique
zero $(z_1,z_2)$ has trivial kernel. Suppose $(v_1,v_2)\in\ker
DF^0(z_1,z_2)$. Since the second component $F_2^0(z_1,z_2)=F_2^0(z_2)$
is independent of $z_1$, this implies $DF_2^0(z_2)v_2=0$
%or in view of~\eqref{eq:Q-F2} equivalently $D\nabla\QQ(z_2)v_2=0$.
In the next subsection we will show: 

\begin{prop}\label{prop:Kepler-Hessian}
The derivative $DF_2^0(z_2)$ at the Kepler orbit $z_2$ has trivial kernel. 
\end{prop}

It follows that $v_2=0$ and $v_1$ satisfies $D_1F_1(z_1,z_2)v_1=0$. In
Section~\ref{ss:Kepler-Hessian-force} we will show:

\begin{prop}\label{prop:Kepler-Hessian-force}
For the Kepler orbit $z_2$, the derivative of the map $z_1\mapsto
F_1(z_1,z_2)$ at its unique zero has trivial kernel. 
\end{prop}

This implies $v_1=0$ and thus concludes the proof of
Theorem~\ref{thm:Euler-mean-intersection}. 

%%%
\subsection{Hessian of the Kepler problem}\label{ss:Kepler-Hessian}
%%%

In this subsection we prove Proposition~\ref{prop:Kepler-Hessian}. 
Consider the map 
\begin{equation}
  F:Z\to H_{sym}^0(S^1,\R),\qquad
  F(z) = -z''+a(z)z,\qquad
  a(z) = \frac{\|z'\|^2}{\|z\|^2}-\frac{1}{\|z\|^6}
\end{equation}
defined on the space
$$
  Z := \{z\in H_{sym}^2(S^1,\R)\mid z(\tau)>0\text{ for all }\tau\in(0,1)\}.
$$
Thus $F$ corresponds to the map $F_2^0$ of the previous subsection
describing simple symmetric solutions of the Kepler problem, where we
have renamed $z_2$ to $z$ and $a_2^0$ to $a$. 
The unique zero of $F$ is given by
$$
   z(\tau)=\zeta\,\sin(\pi\tau), 
$$
where $\zeta>0$ is uniquely determined by the equation $F(z)=0$, or
equivalently
$$
   a(z) = -\pi^2. 
$$
We need to show that the derivative of $F$ at its zero $z$ has trivial
kernel. In direction $v\in H_{sym}^2(S^1,\R)$ it is given by
$$
  DF(z)v
  = -v'' + a(z)v + \bigl(Da(z)v\bigr)z
$$
with
\begin{eqnarray*}
  Da(z)v 
  &=& \frac{2\la z',v'\ra}{\|z\|^2} - \frac{2\|z'\|^2\la
    z,v\ra}{\|z\|^4} + \frac{6\la z,v\ra}{\|z\|^8} \\
  &=& \frac{2}{\|z\|^2}\Bigl\la -z'' - \frac{\|z'\|^2}{\|z\|^2}z +
  \frac{3}{\|z\|^6}z,v\Bigl\ra \\ 
  &=& \frac{2}{\|z\|^2}\Bigl(-2a(z) + \frac{2}{\|z\|^6}\Bigr)\la z,v\ra, \\ 
\end{eqnarray*}
where for the last equality we have used $F(z)=0$. Using
$a(z)=-\pi^2$, it follows that an element $v$ in the kernel of $DF(z)$ satisfies
\begin{equation}\label{eq:v-Kepler}
   -v'' -\pi^2v + b\la z,v\ra z = 0
\end{equation}
with the constant
$$
   b = \frac{4}{\|z\|^2}\Bigl(\pi^2 + \frac{1}{\|z\|^6}\Bigr) > 0.
$$
It follows that $v$ is smooth. Multiplying~\eqref{eq:v-Kepler} by $v$
and integrating from $0$ to $1$ yields 
\begin{equation}\label{eq:v-Kepler-integrated}
   \la -v'',v\ra -\pi^2\|v\|^2 + b\la z,v\ra^2 = 0. 
\end{equation}
Since $v$ extends to an odd $2$-periodic function, it has a Fourier expansion 
$$
   v(\tau) = \sum_{k=1}^\infty c_k\sin(\pi k\tau),\qquad c_k\in\R. 
$$
We deduce $v'(\tau) = \sum_{k=1}^\infty \pi kc_k\cos(\pi k\tau)$
and thus the Poincar\'e inequality
$$
  \la -v'',v\ra = \|v'\|^2
  = \frac12 \sum_{k=1}^\infty \pi^2k^2c_k^2
  \geq \frac12 \sum_{k=1}^\infty \pi^2c_k^2 = \pi^2\|v\|^2.
$$
Hence~\eqref{eq:v-Kepler-integrated} can only hold if $b\la
z,v\ra^2\leq 0$, i.e.~$\la z,v\ra=0$, and equality holds in the
Poincar\'e inequality, i.e.~$c_k=0$ for all $k\geq 2$. Thus
$v(\tau)=c_1\sin(\pi\tau)$, and $\la v,z\ra=0$ implies $v=0$.
This concludes the proof of Proposition~\ref{prop:Kepler-Hessian}.

%%%
\subsection{Hessian of the Kepler problem with constant force}\label{ss:Kepler-Hessian-force}
%%%

In this subsection we prove Proposition~\ref{prop:Kepler-Hessian-force}. 
Denote by $z_2$ the unique symmetric Kepler orbit from the previous
subsection and by $q_2(t)=z_2(\tau(t))^2$ its Levi-Civita transform.
We denote its average using~\eqref{bov1} by
$$
   a := \ol q_2 = \frac{\|z_2^2\|}{\|z_2\|^2} >0. 
$$
Using this, we rewrite the functions $a_1$ and $b_1$ from
Section~\ref{ss:Fredholm-homotopy} as
\begin{eqnarray*}
  a_1(z_1,z_2)
  &=& \frac{||z_1'||^2}{||z_1||^2}-\frac{1}{||z_1||^6}-\frac{||z_2||^4
    \cdot ||z_1^2||^2}{2||z_1||^2 \cdot\big(||z_1^2||^2\cdot ||z_2||^2 
  -||z_2^2||^2 \cdot ||z_1||^2\big)^2},\\
  &=& \frac{||z_1'||^2}{||z_1||^2}-\frac{1}{||z_1||^6}-\frac{
  ||z_1^2||^2}{2||z_1||^2 \cdot\big(||z_1^2||^2 
  -a \cdot ||z_1||^2\big)^2},\\
  b_1(z_1,z_2)
  &=& \frac{||z_2||^4}{\big(||z_1^2||^2\cdot ||z_2||^2-||z_2^2||^2 \cdot ||z_1||^2\big)^2}\\
  &=& \frac{1}{\big(||z_1^2||^2-a \cdot ||z_1||^2\big)^2}.
\end{eqnarray*}
Renaming $z_1$ to $z$, we thus consider the map 
\begin{equation}
  F: W\to H_{sym}^0(S^1,\R),\qquad
  F(z) = -z''+a_1(z)z+b_1(z)z^3
\end{equation}
with 
\begin{eqnarray*}
  a_1(z)
  &:=& \frac{||z'||^2}{||z||^2}-\frac{1}{||z||^6}-\frac{
  ||z^2||^2}{2||z||^2 \cdot\big(||z^2||^2 
  -a \cdot ||z||^2\big)^2},\\
  b_1(z)
  &:=& \frac{1}{\big(||z^2||^2-a \cdot ||z||^2\big)^2} > 0
\end{eqnarray*}
defined on the space
$$
   W := \{z\in H^2(S^1,\R)\mid z(\tau)>0\text{ for all }\tau\in S^1,\;
   \|z^2\|^2>a\|z\|^2\}. 
$$
So $F$ corresponds to the map $z_1\mapsto F_1(z_1,z_2)$ in
Proposition~\ref{prop:Kepler-Hessian-force}. 
By the discussion in Section~\ref{ss:decoupling}, the unique zero of
$F$ is the constant function
$$
   z(\tau)\equiv\ol z, 
$$
where $\ol z>0$ is uniquely determined by the equation $F(z)=0$, or equivalently
\begin{equation}\label{eq:a1b1}
   a_1(\ol z) + b_1(\ol z)\ol z^2 = 0. 
\end{equation}
We need to show that the derivative of $F$ at its zero $\ol z$ has trivial
kernel. In direction $v\in H_{sym}^2(S^1,\R)$ it is given by
$$
  DF(\ol z)v
  = -v'' + (a_1+3b_1\ol z^2)v + c_1
$$
with the constants $a_1=a_1(\ol z)$, $b_2=b_1(\ol z)$ and
$$
   c_1 := \la\nabla a_1(\ol z),v\ra\ol z+\la\nabla b_1(\ol z),v\ra\ol z^3.
$$
Using~\eqref{eq:a1b1} the equation $DF(\ol z)v=0$ thus becomes
\begin{equation}\label{eq:v-Kepler-force}
  v'' = 2b_1\ol z^2v + c_1. 
\end{equation} 
As in the proof of Lemma~\ref{lem:q-estimates} it follows that $v$ is
constant: It attains its maximum at some time $\tau_{\rm
  max}\in S^1$ satisfying
$$
  v''(\tau_{\rm max}) = 2b_1\ol z^2v(\tau_{\rm max}) + c_1 \leq 0.
$$
From $b_1>0$ it follows that for all $\tau\in S^1$ we have
$$
  v''(\tau) = 2b_1\ol z^2v(\tau) + c_1 
  \leq 2b_1\ol z^2v(\tau_{\rm max}) + c_1 \leq 0,
$$
which by periodicity implies that $v(\tau)\equiv\ol v$ is constant.
Plugging this into~\eqref{eq:v-Kepler-force} yields
\begin{equation}\label{eq:lambda-v}
  \lambda\,\ol v=0
\end{equation}
with the constant
$$
  \lambda := 2b_1(\ol z)\ol z^2 + \nabla a_1(\ol z)\ol z + \nabla
  b_1(\ol z)\ol z^3.
$$
To compute $\lambda$, we plug $z=\ol z$ into $a_1$ and $b_1$ to get
\begin{eqnarray*}
  a_1(\ol z)
  &=& -\frac{1}{\ol z^6} - \frac{1}{2\ol z^2(\ol z^2-a)^2},\\
  b_1(\ol z)
  &=& \frac{1}{\ol z^4(\ol z^2-a)^2}
\end{eqnarray*}
and compute their derivatives (as fuctions $\R\to\R$)
\begin{eqnarray*}
  \nabla a_1(\ol z)
  &=& \frac{6}{\ol z^7}
  + \frac{2\ol z(\ol z^2-a)^2+\ol z^2\cdot 2(\ol z^2-a)\cdot 2\ol z}
  {2\ol z^4(\ol z^2-a)^4}\\
  &=& \frac{6}{\ol z^7}
  + \frac{3\ol z^2-a}{\ol z^3(\ol z^2-a)^3},\\
  \nabla b_1(\ol z)
  &=& -\frac{4\ol z^3(\ol z^2-a)^2+\ol z^4\cdot 2(\ol z^2-a)\cdot 2\ol z}
  {\ol z^8(\ol z^2-a)^4}\\
  &=& \frac{4a-8\ol z^2}{\ol z^5(\ol z^2-a)^3}.
\end{eqnarray*}
From equation~\eqref{eq:a1b1} we obtain
$$
  0 = a_1(\ol z) + b_1(\ol z)\ol z^2 
  = -\frac{1}{\ol z^6} - \frac{1}{2\ol z^2(\ol z^2-a)^2}
  + \frac{1}{\ol z^2(\ol z^2-a)^2},
$$
and therefore
$$
  \frac{1}{\ol z^6} = \frac{1}{2\ol z^2(\ol z^2-a)^2}.
$$
Using this and the preceding formulae we compute
\begin{eqnarray*}
  \lambda
  &=& 2b_1(\ol z)\ol z^2 + \nabla a_1(\ol z)\ol z + \nabla
  b_1(\ol z)\ol z^3\\
  &=& \frac{2}{\ol z^2(\ol z^2-a)^2} + \frac{6}{\ol z^6}
  + \frac{3\ol z^2-a+4a-8\ol z^2}{\ol z^3(\ol z^2-a)^3}\\
  &=& \frac{2+3}{\ol z^2(\ol z^2-a)^2} 
  + \frac{3a-5\ol z^2}{\ol z^3(\ol z^2-a)^3}\\
  &=& \frac{5\ol z^2-5a+3a-5\ol z^2}{\ol z^3(\ol z^2-a)^3}\\
  &=& \frac{-2a}{\ol z^3(\ol z^2-a)^3} < 0.\\
\end{eqnarray*}
Hence equation~\eqref{eq:lambda-v} implies $\ol v=0$. This concludes
the proof of Proposition~\ref{prop:Kepler-Hessian-force}, and
therefore of Theorem~\ref{thm:Euler-mean-intersection}.

%%%%%%%%%%%%%%%%%%%%%%%%%%%%%%%%%%%%%%%%%%%%%%%%%%%%%%%%%%%%%%%%%%%%%%%%%%
%%%%%%%%%%%%%%%%%%%%%%%%%%%%%%%%%%%%%%%%%%%%%%%%%%%%%%%%%%%%%%%%%%%%%%%%%%


\begin{thebibliography}{99}
%%%%%%%%%%%%%%%%%%%%%%%%%%%%%%%%%%%%%%%%%%%%%%%%%%%%%%%%%%%%%%%%%%%%%%%%%%
%%%%%%%%%%%%%%%%%%%%%%%%%%%%%%%%%%%%%%%%%%%%%%%%%%%%%%%%%%%%%%%%%%%%%%%%%%

\bibitem{barutello-ortega-verzini} V.\,Barutello, R.\,Ortega, G.\,Verzini, \emph{Regularized variational principles
for the perturbed Kepler problem}, Advances in Mathematics
{\bf 383} (2021), Article no.\,107694.
\bibitem{cieliebak-mundet-salamon} K.~Cieliebak, I.~Mundet i Riera,
  D.~Salamon, {\em Equivariant moduli problems, branched manifolds,
    and the Euler class}, Topology {\bf 42} (2003), no.~3, 641--700.
\bibitem{frauenfelder} U.\,Frauenfelder, \emph{Helium and Hamiltonian delay equations}, 
Isr.\,J.\,Math. (2021), https://doi.org/10.1007/s11856-021-2242-x
\bibitem{frauenfelder2} U.\,Frauenfelder, \emph{A compactness theorem for frozen planets}, 
Jour.\,Top.\,Anal. Online Ready, https://doi.org/10.1142/S1793525321500448 
\bibitem{frauenfelder-weber} U.~Frauenfelder, J.~Weber, 
\emph{The regularized free fall I: Index computations}, 
arXiv:2102.01688v2, to appear in Russ.\,Jour.\,Math.\,Phys.
\bibitem{hofer-wysocki-zehnder} H.\,Hofer, K.\,Wysocki, E.\,Zehnder, \emph{Polyfold and Fredholm
 Theory}, Ergebnisse der Mathematic und ihrer Grenzgebiete, 3.\,Folge. A Series of Modern
 Surveys in Mathematics \textbf{72}, Springer, Cham (2021).
\bibitem{levi-civita} T. Levi-Civita, \emph{Sur la r\'egularisation du probleme des trois corps}, Acta Math., \textbf{42}, 99-144, (1920).
\bibitem{salamon-weber} D.\,Salamon, J.\,Weber, \emph{Floer homology and the heat flow}, GAFA \textbf{16},
no.\,5, (2006), 1050--1138.
\bibitem{smale} S.~Smale, {\em An infinite dimensional version of
  Sard's theorem}, Amer. J.~Math.~87 (1965), 861--866. 
\bibitem{tanner-richter-rost} G.\,Tanner, K.\,Richter, J.\,Rost, \emph{The theory of two-electron atoms: Between ground state and complete fragmentation}, Review of Modern Physics \textbf{72}(2), 497--544 (2000).
\bibitem{weber} J.\,Weber, \emph{Perturbed closed geodesics are periodic orbits: index and
transversality}, Math.\,Z. \textbf{241} (2002), no.\,1, 45--82.
\bibitem{wintgen-richter-tanner} D.\,Wintgen, K.\,Richter, G.\,Tanner, \emph{The Semi-Classical Helium Atom}, in Proceedings of the International School of Physics ``Enrico Fermi'', Course CXIX, 113--143 (1993).
\end{thebibliography}
\end{document}